\documentclass[review]{siamart0516}
\usepackage{amsopn}
\usepackage{algorithmic}
\Crefname{ALC@unique}{Line}{Lines}
\usepackage{array}
\usepackage{amsfonts}
\usepackage{amsmath}
\usepackage{amssymb}
\usepackage{blkarray}
\usepackage{bm}
\usepackage{bbm}
\usepackage{caption}
\usepackage{color, colortbl}
\usepackage{enumitem}
\definecolor{Gray}{gray}{0.9}
\definecolor{Gray2}{gray}{0.95}
\definecolor{Gray3}{gray}{0.7}
\usepackage{epstopdf}
\usepackage{float}
\usepackage{graphicx}
\graphicspath{{./images/}}
\usepackage{lipsum}
\usepackage{longtable}
\ifpdf
  \DeclareGraphicsExtensions{.eps,.pdf,.png,.jpg}
\else
  \DeclareGraphicsExtensions{.eps}
\fi
\usepackage{mathtools}
\usepackage{multirow}
\usepackage[counterclockwise, figuresleft]{rotating}
\usepackage{ulem}
\usepackage{setspace}
\usepackage{subcaption}
\usepackage{tabu}
\usepackage{xr-hyper}
\usepackage{hyperref}

\newcolumntype{C}[1]{>{\centering\let\newline\\\arraybackslash\hspace{0pt}}m{#1}}

\newtheorem{thm}{Theorem}[section]

\newtheorem{corollary}[thm]{Corollary}

\newtheorem{lemma}[thm]{Lemma}
\newtheorem{definition}[thm]{Definition}
\newtheorem{remark}[thm]{Remark}

\newcommand{\Span}{\mathrm{span}}
\newcommand{\argmin}{\arg\min}

\newcommand{\TheTitle}{A Fast Hierarchically Preconditioned Eigensolver \\Based on Multiresolution Matrix Decomposition}
\newcommand{\TheAuthors}{Authors, Authors, Authors, Authors}



\usepackage{fancyhdr}

\title{{\TheTitle}}%\thanks{Submitted to the editors DATE.

\author{
  Thomas Y. Hou
  \and
  De Huang
  \and
  Ka Chun Lam
  \and
  Ziyun Zhang
}

\usepackage{amsopn}

\ifpdf
\hypersetup{
  pdftitle={\TheTitle},
  pdfauthor={\TheAuthors}
}
\fi

\begin{document}\sloppy

\maketitle

% REQURIRED
\begin{abstract}
In this paper we propose a new iterative method to hierarchically compute a relatively large number of leftmost eigenpairs of a sparse symmetric positive matrix under the multiresolution operator compression framework. We exploit the well-conditioned property of every decomposition components by integrating the multiresolution framework into the Implicitly Restarted Lanczos method. We achieve this combination by proposing an extension-refinement iterative scheme, in which the intrinsic idea is to decompose the target spectrum into several segments such that the corresponding eigenproblem in each segment is well-conditioned. Theoretical analysis and numerical illustration are also reported to illustrate the efficiency and effectiveness of this algorithm.
\end{abstract}

% REQUIRED


\section{Introduction}
\label{sec:introduction}

The computation of eigenpairs for large and sparse matrices is one of the most fundamental tasks in many scientific applications. For example, the leftmost eigenpairs (i.e., the $N$ smallest eigenpairs for some $N\in \mathbb{N}$) of a graph laplacian $L$ help revealing the topological information of the corresponding network from real data. One illustrative example is that the multiplicity of the smallest eigenvalue $\lambda_1$ of $L$ coincides with the number of the connected components of the corresponding graph $G$. In particular, the second-smallest eigenvalue of $L$ is well-known as the algebraic connectivity or the Fiedler value of the graph $G$, which is applied to develop algorithms for graph partitioning \cite{chung1997spectral, newman2010networks, ng2001spectral}. Another important example regarding the use of leftmost eigenpairs is the computation of betweenness centrality of graphs as mentioned in \cite{bozzo2012effective, bozzo2013resistance, bergamaschi2015computing}. Computing the leftmost eigenpairs of large and sparse Symmetric Positive Definite (SPD) matrices is also stemmed from the problem of predicting electronic properties in complex structural systems \cite{goedecker1995low}. Such prediction is achieved by solving the Schr\"odinger equation $\mathcal{H} \Psi = \mathcal{E} \Psi$, where $\mathcal{H}$ is the Hamiltonian operator for the system, $\mathcal{E}$ corresponds to the total energy and $|\Psi(r)|^2$ represents the charge density at location $r$. Solving this equation using the Self Consistent Field (SCF) requires computing the eigenpairs of $\mathcal{H}$ repeatedly, which dominates the overall computation cost of the overall iterations. Thus, an efficient algorithm to solve the eigenproblem is indispensable. Usage of leftmost eigenpairs can also be found in vibrational analysis in mechanical engineering \cite{meirovitch1975elements}. In \cite{cocco2013principal}, authors also suggest that the leftmost eigenpairs of the covariance matrix between residues are important to extract functional and structural information about protein families. Efficient algorithms for computing $p$ smallest eigenpairs for relatively large $p$ are therefore crucial in various applications.

As most of the linear systems from engineering problems or networks are typically large and sparse in nature, iterative methods are preferred. Recently, several efficient algorithms have been developed to obtain leftmost eigenpairs of $A$. These include the Jacobi-Davidson (JD) method \cite{sleijpen2000jacobi}, implicit restarted Arnoldi/Lanczos method \cite{calvetti1994implicitly,sorensen1997implicitly,lehoucq1998arpack},  and the Deflation-accelerated Newton method (DACG) \cite{bergamaschi1997asymptotic}. All these methods give promising results \cite{bergamaschi2015computing,martinez2016tuned}, especially for finding a small amount of leftmost eigenpairs. However, as reported in \cite{martinez2016tuned}, the Implicit Restarted Lanczos Method (IRLM) is still the most performing algorithm when a large amount of smallest eigenpairs are required. Therefore, it is highly desirable to develop a new algorithm, based on the architecture of the IRLM, that can further optimize the performance.

The main purpose of this paper is to explore the possibility of exploiting the advantageous energy decomposition framework under the architecture of the IRLM. In particular, we propose a new spectrum-preserving preconditioned hierarchical eigensolver for computing a large amount of smallest eigenpairs. This eigensolver takes full advantage of the intrinsic structure of the given matrix, the nice spectral property in the Lanczos procedure and also the preconditioning characteristics of the Conjugate Gradient method. Given a sparse symmetric positive matrix $A$ which is assumed to be energy decomposable (See \ref{subsec:energy_decomposable} or \Cref{sec:preliminaries} for details), we integrate the well-behaved matrix properties that are inherited from the Multiresolution Matrix Decomposition (MMD) with IRLM. The preconditioner we propose for the Conjugate Gradient method can also preserve the narrowed residual spectrum of $A$ during the Lanzcos procedure. Throughout this paper, theoretical performance of our proposed algorithm is analyzed rigorously and we conduct a number of numerical experiments to verify the efficacy and effectiveness of the algorithm in practice. To summarize, our contributions are three-fold:
\begin{itemize}
\item We propose a hierarchical framework to compute a relatively large number of leftmost eigenpairs of a sparse symmetric positive matrix. This framework employs the MMD algorithm to further optimize the performance of IRLM. In particular, a specially designed spectrum-preserving preconditioner is introduced for the Conjugate Gradient method to solve for $A^{-1}$.
\item The proposed framework improves the running time of finding $m_{tar}$ smallest eigenpairs of a matrix $A \in \mathbb{R}^{n\times n}$ from $O(m_{tar} \cdot \kappa(A)\cdot nnz(A) \log \frac{1}{\varepsilon})$ (which is achieved by the classical IRLM) to $O\left(m_{tar} \cdot nnz(A)\cdot (\log\frac{1}{\varepsilon} + \log n)^{C}\right)$, where $\kappa(A)$ is the condition number of $A$, $nnz(\cdot)$ is the number of nonzero entries and $C$ is some small constant independent of $m_{tar}, nnz(A)$ and $\kappa(A)$.

%$O(N \cdot \log\kappa(A) m \log \frac{1}{\epsilon}) + O(n \log n)$, where $\kappa(A)$ is the condition number of $A$, $m$ is the number of nonzero entries and $C$ is some constant independent of $m, n$ and $\kappa(A)$.
\item We also provide a rigorous analysis on both the accuracy and the asymptotic computational complexity of our proposed algorithm. This ensures the correctness and efficiency of the algorithm even in large-scale, ill-conditioned scenarios.
\end{itemize}

\subsection{Overview of the algorithm}
In this paper, we propose and develop an iterative scheme under the framework of energy decomposition introduced in \cite{hou2017adaptive}. Under this framework, we can decompose $A^{-1} \in \mathbb{R}^{n \times n}$ into
\[
A^{-1} = P^A_{\mathcal{U}} A^{-1} + P^A_{\Psi} A^{-1} := P^A_{\mathcal{U}} A^{-1} + \Theta,
\]
where $[\mathcal{U}, \Psi]$ corresponds to a basis of $\mathbb{R}^n$; $P^A_{\mathcal{U}}$ and $P^A_{\Psi}$ are the corresponding subspace projections. Recursively, we can also consider $\Theta$ as a ``new'' $A^{-1}$ and decompose $\Theta$ in the same manner. This will give a MMD of $A^{-1} = \sum^K_{k=1} P^A_{\mathcal{U}^{(k)}} A^{-1} + \Theta^{(K)}$. To illustrate, we first consider a 1-level decomposition, i.e., $K = 1$. One important observation regarding this decomposition is that the spectrum of the original operator $A^{-1}$ resembles that of the compressed operator $\Theta$. In particular, if $\lambda_{i,\Theta}$ is the $i^{\text{th}}$ smallest eigenvalue of $\Theta$ and $\zeta_{i,\Theta}$ is the corresponding eigenvector, then $( \lambda_{i,\Theta}^{-1}, \zeta_{i,\Theta})$ is a good approximation of $(\lambda_i^{-1}, q_i)$ for small $\lambda_i$, where $(\lambda_i^{-1}, q_i)$ denotes the $i^{\text{th}}$ eigenpair of $A^{-1}$. These approximate eigenpairs $(\lambda_{i,\Theta}^{-1}, \zeta_{i,\Theta})$ can then used as the initial approximation of the required eigenpairs. Notice that compression errors are introduced into these eigenpairs by the matrix decomposition. Therefore, a refinement procedure should be carried out to diminish these errors up to the prescribed accuracy. Once we obtain the refined eigenpairs, we may extend the spectrum in order to obtain the required amount of eigenpairs. As observed in \cite{martinez2016tuned}, the Implicit Restarted Lanczos Method (IRLM) is the most performing algorithm when large eigenpairs are considered, we therefore employ the Krylov subspace extension technique to extend spectrum up to some prescribed control of the well-posedness. Intuitively, the MMD decomposes the spectrum of $A^{-1}$ into different segments of different scales. Using a subset of the decomposed components to approximate $A^{-1}$ yields a great reduction of the relative condition number. Thus, we can further trim down the complexity of the IRLM by approximating $A^{-1}$ during the shifting process.

To generalize, we propose a hierarchical scheme to compute the leftmost eigenpairs of an energy decomposable matrix. Given the $K$-level multiresolution decomposition $\{\Theta^{(k)}\}^K_{k=1}$ of an energy decomposable matrix $A$, we first compute the eigen decomposition $[V_{ex}^{(K)}, D_{ex}^{(K)}]$ of $\Theta^{(K)}$ (with dimension $N^{(K)}$) corresponding to the coarsest level by using some standard direct method. Then we propose an compatible refinement scheme for both $V_{ex}^{(K)}$ and $D_{ex}^{(K)}$ to obtain $V_{ini}^{(K-1)}$ and $D_{ini}^{(K-1)}$, which will then be the initial spectrum in the consecutive finer level. The efficiency of the cross-level refinement is achieved by a modified version of the orthogonal iteration with the Ritz Acceleration, where we exploit the proximity of the eigenspace across levels to accelerate the Conjugate gradient (CG) method within the refinement step. Using this refined initial spectrum, our second stage is to extend spectrum up to some prescribed control of the well-posedness using the Implicit Restarted Lanczos architecture. Recall that a shifting approach is introduced to reduce the iteration number for the extension, which again requires solving $A^{(K-1)} x = w$ with the CG method in each iteration. However, the preconditioner for CG when we are solving for $A^{(K-1)}w$ must be chosen carefully. Otherwise the orthogonal property brought about by the Krylov subspace methods may not be utilized and a large CG iteration number will be observed (See \Cref{sec:compare_IRLM}). In view of this, we propose a spectrum-preserving hierarchical preconditioner $M^{(K-1)} := (\mathrm{\Psi}^{(K-1)})^T \mathrm{\Psi}^{(K-1)}$ for accelerating the CG iteration during the Lanczos iteration. In particular, we can show that using the preconditioner $M^{(K-1)}$, the number of Preconditioned Conjugate gradient (PCG) iteration to achieve a relative $\varepsilon$ in $A^{(K-1)}$-norm can be controlled in terms of the condition factor $\delta(\mathcal{P})$ (from the energy decomposition of the matrix) and an extension threshold $\mu_{ex}^{(K-1)}$.

This process then repeats hierarchically until we reach the finest level. Under this framework, the condition number of every engaged operators is controlled. The overall accuracy of our proposed algorithm is also determined by the prescribed compression error at the highest level.

\subsection{Previous Works}
Several important iterative methods have been proposed to tackle the eigenproblems of SPD matrices. One of the well established algorithms is the Implicitly Restarted Lanczos Method (IRLM) (or the Implicitly Restarted Arnoldi Method (IRAM) for unsymmetric sparse matrices), which has been implemented in various popular scientific computing packages like MATLAB, R and ARPACK. The IRLM combines both the techniques of the implicitly shifted QR method and the shifting of the operators to avoid the difficulties for obtaining the leftmost eigenpairs. Another popular algorithm for finding leftmost eigenpairs is the Jacobi-Davidson method. The main idea is to minimize the Rayleigh Quotient $q(x) = \frac{x^TAx}{x^Tx}$ using a Newton-type methodology. Efficacy and stability of the algorithm are then achieved by using a projected simplification of the Hessian of the Rayleigh Quotient namely, $\tilde{J}(x_k) := (I - x_k x_k^T)(A - q(x_k)I)(I - x_k x_k^T)$ with the update of $x_k$ to be
\begin{equation}
\label{eqt:JD}
x_{k+1} = x_k - \tilde{J}(x_k)^{-1}(Ax_k - q(x_k)x_k).
\end{equation}
Notice that the advantage of such approach is the low accuracy requirement for solving \cref{eqt:JD}. A parallelization was also proposed \cite{romero2010parallel}. In \cite{bergamaschi1997asymptotic}, the authors proposed the Deflation Accelerated Conjugate Gradient (DACG) method designed for solving the eigenproblem of SPD matrices. The main idea is to replace the Newton's minimization procedure of the Rayleigh quotient $r(x)$ by the nonlinear Conjugate Gradient method which avoids solving linear systems within the algorithm. A comprehensive numerical comparison between the three algorithms was reported in \cite{bergamaschi2015computing}. Recently, Mart\'inez \cite{martinez2016tuned} studied a class of tuned preconditioners for accelerating both the DACG and the IRLM for the computation of the smallest set of eigenpairs of large and sparse SPD matrices. However, as reported in \cite{martinez2016tuned}, the IRLM still outperforms the others when a relatively large number of leftmost eigenpairs is desired. By virtue of this, we are motivated to develop a more efficient algorithm particularly designed for computing a considerable amount of leftmost eigenpairs.

Another class of methods related to localized spectrum is the compression of the eigenmodes. One of the representative pioneer works is proposed by Ozoliņ\v{s} et al. in \cite{ozolicnvs2013compressed}. The goal of this work is to obtain a spatially localized solution of a class of problems in mathematical physics by constructing the compressed modes. In particular, finding these localized modes can be formulated as an optimization problem
\[
\Psi_N = \text{arg}\min_{\hat{\Psi}^N } \frac{1}{\mu} \Vert \Psi_N \Vert_{1} + \text{Tr}(\hat{\Psi}^T_N H \hat{\Psi}_N) \quad \text{such that} \quad \hat{\Psi}_N^T\hat{\Psi}_N = I.
\]
The authors in \cite{ozolicnvs2013compressed} proposed an algorithm based on the split Bregman iteration to solve the $L_1$ minimization problem. By replacing the discrete operator $H$ by the graph Laplacian matrix $A$, one obtains the $L_1$ regularized Principal component analysis (PCA). In particular, if there is no $L_1$ regularization term in the optimization problem, the optimal $\Psi_N$ will be the first $m_{tar}$ eigenvectors of $A$. In other words, this procedure provides an effective way to obtain $N$ (where $N \geq m_{tar}$) localized basis functions that can approximately span the $m_{tar}$ leftmost eigenspace (i.e., eigenspace spanned by the $m_{tar}$ eigenvectors corresponding to the leftmost eigenvalues). Similarly, the MMD framework provides us the hierarchical and sparse/localized basis $\Psi$. These localized basis functions capture the compressed modes and eventually provide us a convenient way to control the complexity of the Eigensolver. %In particular, our proposed algorithm also provides a hierarchical formulation to construct localized modes.

Stiffness matrices discretizing heterogeneous and rough elliptic operators, or graph Laplacians representing general sparse networks are commonly found in practice. Recently, the problem of compressing these SPD matrices has been tackled in different perspectives. M{\aa}lqvist and Petersein \cite{maalqvist2014localization} proposed the use of modified coarse space in order to handle roughness of the coefficients when solving elliptic equations with Finite Element Methods. They construct localized multiscale basis functions from the modified coarse space $V_{H}^{ms}=V_{H}-\mathfrak{F}V_{H}$, where $V_{H}$ is the original coarse space spanned by nodal basis, and $\mathfrak{F}$ is the energy projection onto the space $(V_{H})^\perp$. The exponential decaying property of these modified basis functions has been shown both theoretically and numerically. In \cite{owhadi2017multigrid}, Owhadi reformulated the problem from the decision theory perspective using the idea of {\it Gamblets} as the modified basis. In particular, a coarse space $\Phi$ of measurement functions is constructed from the Bayesian perspective, and the gamblet space is explicitly given as $\Psi=A^{-1}(\Phi)$, which turns out to be a counterpart of the modified coarse space in \cite{maalqvist2014localization}. The exponential decaying property of these localized basis functions is also proved independently using the idea of gamblets. Hou and Zhang in \cite{hou2016sparse} further extended these works and constructed localized basis functions for higher order strongly elliptic operators. To further promote the operator compression for situations where the physical domain is unknown or is embedded in some nontrivial high dimensional manifolds, Hou et. al. propose to exploit the local spectrum information of a general class of SPD matrices to by-pass the needs of adopting knowledge of computational domain during the construction of local basis. Recently, Sch\"afer et. al \cite{schafer2017compression} proposed a near-linear running time algorithm to compress a large class of dense kernel matrices $\Theta \in \mathbb{R}^{n\times n}$. The authors also provided rigorous complexity analyses and showed that the complexity of the proposed algorithm is $O(n \log(n) \log^d (n/\epsilon))$ in space and $O(n \log^2(n) \log^{2d}(n/\epsilon))$ in time, where $d$ is the intrinsic dimension of the problem.

\subsection{Outline}
The layout of the rest of this paper is as follows: In \Cref{sec:preliminaries} we review the Energy Decomposition framework for symmetric positive definite matrices proposed in \cite{hou2017adaptive} and in particular, a brief review of the operator compression and multiresolution matrix decomposition is summarized. This is then followed by the review of the implicitly restarted Arnoldi iteration procedure. Some error analysis and perturbation theories subject to our operator compression framework are discussed. Theoretical developments and algorithms of the hierarchical spectrum extension/compression and the eigenpair refinement are then proposed in \Cref{sec:spectrum_extension} and \Cref{sec:eigenpair_refinement} respectively. Combining these two methods, we propose our hierarchical eigensolver in \Cref{sec:overall_algorithm}, where details of the choice of parameters are discussed. \Cref{sec:experiment} is devoted to experimental results to justify the effectiveness of our proposed algorithm. In \cref{sec:compare_IRLM}, we provide a quantitative numerical comparison with the IRLM. The numerical results show that our proposed algorithm gives a promising results in terms of runtime complexity. Discussion of future works and conclusion are drawn in \Cref{sec:conclusion}.

\section{Preliminaries}
\label{sec:preliminaries}
The purpose of this section is to provide a general summary of the {\bf Energy Decomposition} framework for operator compression and multiresolution matrix decomposition. One may refer to \cite{hou2017adaptive} for detailed numerical analysis and experimental results.

\subsection{Energy Decomposition}
\label{subsec:energy_decomposable}
Let $A$ be a $n \times n$ symmetric positive definite (SPD) matrix. We call $\mathcal{E} = \{ E_k \}_{k=1}^m$ an \textbf{energy decomposition} of $A$ and $E_k$ to be an \textbf{energy element} of $A$ if we can express $A = \sum_{k=1}^m E_k$, where $E_k \succeq 0\ \forall k = 1,\ldots, m$. For the ease of discussion, we always assume that the given $\mathcal{E} = \{ E_k \}^m_{k=1}$ is the finest underlying energy decomposition of $A$, meaning that no $E_k \in \mathcal{E}$ can be further decomposed as $E_k = E_{k,1} + E_{k,2}$.

Let $\mathcal{V}$ be a basis of $\mathbb{R}^n$. For any subset $\mathcal{S} \subset \mathcal{V}$, we denote $P_{\mathcal{S}}$ as the orthogonal projection onto $\mathcal{S}$. Following the notations in \cite{hou2017adaptive}, we also denote $A_\mathcal{S}$, $\underline{A}_\mathcal{S}$ and $\overline{A}_\mathcal{S}$ as the \textbf{restricted}, \textbf{interior} and \textbf{closed energy} of $\mathcal{S}$ with respect to $A$ and $\mathcal{E}$.

\subsection{Operator Compression}
The procedures of compressing the solver $A^{-1}$ with broad-banded spectrum are: (i) construct a partition of the computational basis using local information of $A$; (ii) construct the coarse space $\Phi$ that is locally computable and has good interpolation property; (iii) construct the modified coarse space $\Psi = A^{-1}(\Phi)$ of $\mathbb{R}^n$ as proposed in \cite{hou2016sparse,maalqvist2014localization,owhadi2017multigrid}. If an appropriate partitioning is given, we have the following error estimate for operator compression.
 
% \begin{thm}
% Let $\mathcal{P} = \{ P_j \}^M_{j=1}$ be a partition of $\mathcal{V}$. For each $1 \leq j \leq M$, let $\Phi_j$ be some subspace of $\text{span}\{P_j\}$ such that
% \begin{equation}
% \|x-P_{\Phi_j}x\|_2\leq \sqrt{\epsilon} \|x\|_{\underline{A}_{P_j}},\quad \forall x\in \mathrm{span}\{P_j\},
% \label{eqt:localphi_constraint}
% \end{equation}
% where $P_{\Phi_j}$ is the orthogonal projection onto $\Phi_j$. Take $\Phi=\bigoplus_{j}\Phi_j$ and let $\Psi$ be a subspace of $\mathbb{R}^n$ given by $\Psi = A^{-1}(\Phi)$. Denote $P_{\Psi}^A$ to be the orthogonal projection onto $\Psi$ with respect to $\langle \cdot, \cdot \rangle_A$, i.e.,
% \begin{equation}
% P_\Psi^A=\Psi(\Psi^TA\Psi)^{-1}\Psi^TA=A^{-1}\Phi(\Phi^TA^{-1}\Phi)^{-1}\Phi^T
% \label{eqt:projection}
% \end{equation}
% and $\Theta = P^A_{\Psi} A^{-1}$ to be the rank-$N$ compressed approximation of $A^{-1}$. Then for any $x \in \mathbb{R}^n$, and $b=Ax$, we have
% \begin{equation}
% \|x-P_\Psi^A x\|_A\leq \sqrt{\epsilon}\|b\|_2 \quad \text{and} \quad \|x-P_\Psi^A x\|_2\leq \epsilon\|b\|_2,
% \end{equation}
% and thus
% \begin{equation}
% \|A^{-1}-\Theta\|_2\leq \epsilon.
% \end{equation}
% \label{thm:compression_thm1}
% \end{thm}

\begin{thm}
Let $\Phi$ be a $N$ dimensional subspace of $\mathbb{R}^n$ such that for some $\epsilon>0$, 
\begin{equation}
\|x-P_{\Phi}x\|_2\leq \sqrt{\epsilon} \|x\|_{A},\quad \forall x\in \mathbb{R}^n,
\label{eqt:localphi_constraint}
\end{equation}
where $P_{\Phi}$ is the orthogonal projection onto $\Phi$. Let $\Psi$ be a subspace of $\mathbb{R}^n$ given by $\Psi = A^{-1}(\Phi)$. Denote $P_{\Psi}^A$ as the orthogonal projection onto $\Psi$ with respect to $\langle \cdot, \cdot \rangle_A$, and $\Theta = P^A_{\Psi} A^{-1}$ as the rank-$N$ compressed approximation of $A^{-1}$. Then for any $x \in \mathbb{R}^n$, and $b=Ax$, we have
\begin{equation}
\|x-P_\Psi^A x\|_A\leq \sqrt{\epsilon}\|b\|_2 \quad \text{and} \quad \|x-P_\Psi^A x\|_2\leq \epsilon\|b\|_2,
\end{equation}
and thus
\begin{equation}
\|A^{-1}-\Theta\|_2\leq \epsilon.
\end{equation}
\label{thm:compression_thm1}
\end{thm}

As discussed in \cite{hou2017adaptive}, to satisfy \cref{eqt:localphi_constraint}, $\Phi$ can be constructed by choosing some optimal local basis $\Phi_j$ on each patch $P_j$, where $\mathcal{P} = \{ P_j \}^M_{j=1}$ is a partition of $\mathcal{V}$. To minimize $\dim \Phi$, the local basis $\Phi_j$ is chosen to be the eigenvectors corresponding to the smallest interior eigenvalues (i.e., eigenvalues of $\underline{A}_{P_j}$) $\lambda_1(P_j) \leq \lambda_2(P_j) \leq \cdots \leq \lambda_{q_j(\epsilon)}(P_j)$, where $q_j(\epsilon)$ is the smallest integer such that $\frac{1}{\epsilon} \leq \lambda_{q_j(\epsilon)}(P_j)$. By reversing the statement, we introduce the {\bf error factor} $\varepsilon(\mathcal{P}) = \max_j (\lambda_{q+1}(P_j))^{-1}$ of partition $\mathcal{P}$, where $q$ is some prescribed uniform integer for all patches. Then locally on each patch we have $\|x-P_{\Phi_j}x\|_2\leq \sqrt{\varepsilon(\mathcal{P})} \|x\|_{\underline{A}_{P_j}}, \forall x\in \mathrm{span}\{P_j\}$, and by collecting $\Phi=\bigoplus_{j}\Phi_j$ we have globally $\|x-P_{\Phi}x\|_2\leq \sqrt{\varepsilon(\mathcal{P})} \|x\|_{A},\forall x\in \mathbb{R}^n$. In the following, we assume that $q=1$ in all cases. Under this setting, the problem of minimizing $\dim \Phi$ subject to \cref{eqt:localphi_constraint} is transformed into finding a partition $\mathcal{P} = \{ P_j \}_{j=1}^N$ with minimal patch number and satisfies $\varepsilon(\mathcal{P}) \leq \epsilon$.

Following the notations in \cite{hou2017adaptive}, we also use $\Phi,\Psi$ to denote the matrices whose columns are the basis vectors of the subspaces $\Phi,\Psi$ respectively. We remark that using the matrix form, the $A$-orthogonal projection $P_\Psi^A$ can be written as 
\begin{equation}
P_\Psi^A=\Psi(\Psi^TA\Psi)^{-1}\Psi^TA=A^{-1}\Phi(\Phi^TA^{-1}\Phi)^{-1}\Phi^T,
\label{eqt:projection}
\end{equation}
and the rank-$N$ compressed approximation is explicitly $\Theta = P^A_{\Psi} A^{-1}=\Psi A_{st}^{-1}\Psi^T$, where
\begin{equation}
A_{st} = \Psi^TA\Psi
\label{eqt:stiffness}
\end{equation}
is the stiffness matrix in the basis $\Psi$. Once the coarse space/basis $\Phi$ is constructed, the next step is to find $\Psi=[\psi_1,\psi_2,\cdots,\psi_N]= A^{-1}(\Phi)$ such that (i) the stiffness matrix $A_{\text{st}}$ has a relatively small condition number, or the condition number can be bounded by some local information; (ii) each $\psi_i$ is locally computable, or can be approximated by some $\widetilde{\psi}_i$ that is locally computable. To achieve these two requirements, we impose the correlation condition $\Phi^T\Psi=I_{N}$, which is equivalent to choosing $\Psi = [\psi_1,\psi_2,\ldots, \psi_N]$ to be
\begin{equation}
\Psi= A^{-1} \Phi (\Phi A^{-1} \Phi)^{-1}
\label{eqt:psi_closedform}
\end{equation}
and we have the following theorem for the well-posedness of $A_{\text{st}}$:

\begin{thm} Let $A_{st}$ be the stiffness matrix given by \cref{eqt:stiffness}. Let $\lambda_{\min}(A_{\text{st}})$ and $\lambda_{\max}(A_{\text{st}})$ denote the smallest and largest eigenvalues of $A_{\text{st}}$ respectively, then we have
\begin{equation}
\lambda_{\min}(A_{\text{st}})\geq \lambda_{\min}(A),\qquad \lambda_{\max}(A_{\text{st}})\leq \delta(\mathcal{P}),
\label{eqt:lambdamin_control} 
\end{equation}
with
$$ \delta(\mathcal{P}) = \delta(\mathcal{P},\Phi) = \underset{P_j \in \mathcal{P}}{\max} \ \delta(P_j,\Phi_j) \quad \text{ and } \quad \delta(P_j,\Phi_j) = \underset{x \in \Phi_j}{\max} \ \frac{x^T x}{x^T \overline{A}_{P_j}^{-1}x},$$
where $\delta(\mathcal{P})$ is called the {\bf condition factor} of the partition $\mathcal{P}$.
\label{thm:condition_number}
\end{thm}

In other words, by defining $\Psi$ as in \cref{eqt:psi_closedform}, the first requirement can be satisfied. Moreover, such choice of $\Psi$ also satisfies the second requirement. In fact, we can prove the spatial exponential decaying property of every basis function $\psi_i$ (See \cite{hou2017adaptive}, \cite{owhadi2017multigrid} for details). This fast decay feature makes it possible to approximate $\Psi$ by some localized basis $\widetilde{\Psi}$ that preserves the good properties of $\Psi$.  In particular, we can construct a basis $\widetilde{\Psi}=[\widetilde{\psi}_1,\widetilde{\psi}_2,\cdots,\widetilde{\psi}_N]$ such that each $\widetilde{\psi}_i$ satisfies $\|\psi_i-\widetilde{\psi}_i\|_A\leq C\sqrt{\frac{\epsilon}{N}}$ for some constant $C$, and has support size $O((\log\frac{1}{\epsilon} + \log N)^d)$, where $d$ is the intrinsic dimension of the problem that characterizes its connectivity. For this localized $\widetilde{\Psi}$, we have an analogy of \cref{thm:compression_thm1} stating that the operator compression error can be bounded by $\|A^{-1}-\widetilde{\Theta}\|_2\leq (1+C\|A^{-1}\|_2)^2\varepsilon(\mathcal{P})$ (where $\widetilde{\Theta} := P_{\widetilde{\Psi}}^A A^{-1} = \widetilde{\Psi}(\widetilde{\Psi}^T A \widetilde{\Psi})^{-1} \widetilde{\Psi}^T$), and the condition bound of the localized stiffness matrix can be estimated by
\begin{equation}
\kappa(\widetilde{A}_{\text{st}}) = \frac{\lambda_{\max}(\widetilde{A}_{\text{st}})}{\lambda_{\min}(\widetilde{A}_{\text{st}})} \leq \left(1+C\sqrt{\frac{\epsilon}{\delta(\mathcal{P})}}\right)^2\delta(\mathcal{P}) \| A^{-1} \|_2,
\label{eqt:lambdatildemin_control} 
\end{equation}
where $\kappa(\widetilde{A}_{\text{st}})$ is the condition number of $\widetilde{A}_{\text{st}} := \widetilde{\Psi}^T A \widetilde{\Psi}$. Therefore the burden of controlling the accuracy, sparsity and well-posedness of the compressed operator $A_{\text{st}}$ falls into the procedure of partitioning. We then propose a nearly-linear time algorithm using the indicators {\bf error factor} and {\bf condition factor} to obtain an appropriate partition $\mathcal{P}$ subject to $\varepsilon(\mathcal{P}) \delta(\mathcal{P}) \leq c$ for some prescribed upper bound $c$. For details of the notations and the algorithm, please refer to \cite{hou2017adaptive}.

\subsection{Multiresolution Matrix Decomposition}
Recall that the main purpose of decomposing $A^{-1}$ into hierarchical resolutions is to resolve the difficulty of large condition number $\kappa(A)$ when solving the linear system $Ax = b$. Through decomposition, the relative condition number in each scale/level can be bounded by some prescribed value. Using the notation as in the previous subsections, we denote $U = [U_1, U_2, \cdots, U_M]$ and therefore $[U, \Psi]$ forms a basis of $\mathbb{R}^n$. We also have $U^T A \Psi = U^T \Phi (\Phi^T A^{-1} \Phi)^{-1} = 0$. Thus the inverse of $A$ can be written as
\begin{equation}
\begin{split}
A^{-1}=&\ \Big(\left[\begin{array}{c} U^T\\ \Psi^T\end{array}\right]^{-1}\left[\begin{array}{c} U^T\\ \Psi^T\end{array}\right]A\left[\begin{array}{cc} U & \Psi\end{array}\right]\left[\begin{array}{cc} U & \Psi\end{array}\right]^{-1}\Big)^{-1} \\
=&\ U(\underbrace{U^TAU}_{B_{\text{st}}})^{-1}U^T +\Psi(\underbrace{\Psi^TA\Psi}_{A_{\text{st}}})^{-1}\Psi^T.
\end{split}
\label{eqt:Ainv_representation}
\end{equation}
Therefore, solving $A^{-1}b$ is equivalent to solving $A^{-1}_{\text{st}}(\Psi^T b)$ and $B^{-1}_{\text{st}}(U^T b)$ separately. For $B_{\text{st}}$, since the sparsity of $U$ will be inherited to $B_{\text{st}}$, it will be efficient to solve $B^{-1}_{\text{st}}b$ if $\kappa(B_{\text{st}})$ is bounded. The following lemma estimates such upper bound.
\begin{lemma}
If $\Phi$ satisfies the condition as in \cref{thm:compression_thm1} with $\epsilon(\mathcal{P})$ and $B_{\text{st}} = U^T A U$, then 
\begin{equation}
\lambda_{\max}(B_{\text{st}})\leq\lambda_{\max}(A)\cdot\lambda_{\max}(U^T U),\qquad \lambda_{\min}(B_{\text{st}})\geq \frac{1}{\varepsilon(\mathcal{P})}\cdot\lambda_{\min}(U^T U),
\end{equation}
and thus 
\begin{equation}
\kappa(B_{\text{st}})\leq \varepsilon(\mathcal{P})\cdot\lambda_{\max}(A)\cdot\kappa(U^TU).
\end{equation}
\label{lemma:Bst_eigenvalue}
\end{lemma}
Notice that $U^T U$ is block-diagonal with blocks $U_j^T U_j$, therefore
\begin{equation}
\kappa(U^TU) = \frac{\lambda_{\max}(U^TU)}{\lambda_{\min}(U^TU)} = \frac{\max_{1\leq j\leq M}{\lambda_{\max}(U_j^T U_j)}}{\min_{1\leq j\leq M}{\lambda_{\min}(U_j^T U_j)}}.
\end{equation}
In particular, if we extend $\Phi_j$ to an orthonormal basis of $\text{span}\{P_j\}$ to get $U_j$ using the QR factorization, we have $\kappa(U^T U) = 1$. So if the condition number of $A$ is huge, we can first set a small enough $\varepsilon$ to sufficiently bound $\kappa(B_{\text{st}})$; if $\kappa(A_{\text{st}})$ is still large, we apply the decomposition to $A_{\text{st}}^{-1}$ again to further decompose $\kappa(A_{\text{st}})$. In order to further decompose the stiffness matrix $A_{\text{st}}$, we need to construct the corresponding energy decomposition of $A_{\text{st}}$.

\begin{definition}[Inherited energy decomposition]
Let $\mathcal{E} = \{E_k\}_{k=1}^m$ be the energy decomposition of $A$, then the inherited energy decomposition of $A_{\text{st}} = \Psi^T A \Psi$ with respect to $\mathcal{E}$ is simply given by $\mathcal{E}^{\Psi} = \{ E_k^{\Psi}\}^m_{k=1}$, where $E^{\Psi}_{k} = \Psi^T E_k \Psi,\quad k = 1,2,\cdots, m.$
\end{definition}

Once we have the underlying energy decomposition of $A_{\text{st}} $, we can repeat the procedure to decompose $A_{\text{st}}^{-1}$ in $\mathbb{R}^N$ as what we have done to $A^{-1}$ in $\mathbb{R}^n$, and furthermore to obtain a multi-level decomposition of $A^{-1}$. In particular, at level $k$, we construct the partition $\mathcal{P}^{(k)}$ and the basis $\Phi^{(k)},U^{(k)},\Psi^{(k)}$ accordingly, and decompose $(A^{(k)})^{-1}$ as
\[
(A^{(k)})^{-1} = U^{(k+1)}\big((U^{(k+1)})^TA^{(k)}U^{(k+1)}\big)^{-1}(U^{(k+1)})^T + \Psi^{(k+1)}\big((\Psi^{(k+1)})^TA^{(k)}\Psi^{(k+1)}\big)^{-1}(\Psi^{(k+1)})^T,
\]
and then define $A^{(k+1)}=(\Psi^{(k+1)})^TA^{(k)}\Psi^{(k+1)}$ and $B^{(k+1)}=(U^{(k+1)})^TA^{(k)}U^{(k+1)}$. We also recall the following notations
\begin{subequations}
\begin{align}
\bm{\Phi}^{(k)}=\Phi^{(1)}\cdots\Phi^{(k-1)}\Phi^{(k)},\quad k\geq 1, \label{eqt:mutliresolution_phi&U&psi_a}\\
\mathcal{U}^{(k)}=\Psi^{(1)}\cdots\Psi^{(k-1)}U^{(k)},\quad k\geq1, \label{eqt:mutliresolution_phi&U&psi_b}\\
\bm{\Psi}^{(k)}=\Psi^{(1)}\cdots\Psi^{(k-1)}\Psi^{(k)},\quad k\geq 1.
\label{eqt:mutliresolution_phi&U&psi_c}
\end{align}
\label{eqt:mutliresolution_phi&U&psi}
\end{subequations}
Using these notations and noticing that $(\Phi^{(k)})^T\Phi^{(k)}=(\Phi^{(k)})^T\Psi^{(k)}=I_{N^{(k)}}$, we have 
\[A^{(k)}=(\bm{\Psi}^{(k)})^TA\bm{\Psi}^{(k)}=\big((\bm{\Phi}^{(k)})^TA^{-1}\bm{\Phi}^{(k)}\big)^{-1},\quad B^{(k)}=(\mathcal{U}^{(k)})^TA\mathcal{U}^{(k)},\]
\[(\bm{\Phi}^{(k)})^T\bm{\Phi}^{(k)}=(\bm{\Phi}^{(k)})^T\bm{\Psi}^{(k)}=I_{N^{(k)}},\quad \bm{\Psi}^{(k)}=A^{-1}\bm{\Phi}^{(k)}\big((\bm{\Phi}^{(k)})^TA^{-1}\bm{\Phi}^{(k)}\big)^{-1},\]
and for any integer $K$, 
\begin{equation}
A^{-1} = (A^{(0)})^{-1} = \sum_{k=1}^{K}\mathcal{U}^{(k)}\big((\mathcal{U}^{(k)})^TA\mathcal{U}^{(k)}\big)^{-1}(\mathcal{U}^{(k)})^T+\bm{\Psi}^{(K)}\big((\bm{\Psi}^{(K)})^TA\bm{\Psi}^{(K)}\big)^{-1}(\bm{\Psi}^{(K)})^T.
\label{eqt:MMD}
\end{equation}
We call \cref{eqt:MMD} the Multiresolution Matrix Decomposition (MMD) of $A^{-1}$. We remark that as $k$ increases, the compressed dimension $N^{(k)}$ decreases, and the scale of the subspace spanned by $\bm{\Psi}^{(k)}$ becomes coarser. In the subspace spanned by $\bm{\Psi}^{(k-1)}$, the basis $\mathcal{U}^{(k)}$ represents the features that are finer than $\bm{\Psi}^{(k)}$. This decomposition helps separate $A$ that has a large condition number into a sequence of matrices with more controllable conditioned numbers. This is stated in the following corollary.
\begin{corollary} We have
\[\kappa(A^{(k)})\leq \delta(\mathcal{P}^{(k)})\|A^{-1}\|_2,\]
\[\kappa(B^{(k)})\leq \varepsilon(\mathcal{P}^{(k)})\delta(\mathcal{P}^{(k-1)})\kappa\big((U^{(k)})^TU^{(k)}\big).\]
For consistency, we write $\delta(\mathcal{P}^{(0)})=\lambda_{\max}(A^{(0)})=\lambda_{\max}(A)$.
\label{cor:Bconditonnumber}
\end{corollary}

The following theorem provides an estimation of the total compression error under $K$ levels of matrix decomposition.
\begin{thm}
\label{thm:error_accumulate}
Assume we have constructed $\Phi^{(k)}, k=1,2,\cdots,K$ on each level accordingly, then we have
\begin{equation}
\|x-P_{\bm{\Phi}^{(k)}}x\|_2^2\leq \varepsilon_k\|x\|^2_A\quad \forall x\in \mathbb{R}^n,\quad \text{where } \varepsilon_k = \sum_{k'=1}^{k}\varepsilon(\mathcal{P}^{(k')}),
\label{eqt:error_accumulate}
\end{equation}
and thus for any $x\in \mathbb{R}^n$ and $b=Ax$, we have
\[\|x-P_{\bm{\Psi}^{(k)}}^Ax\|^2_A\leq \varepsilon_k\|b\|^2_2, \quad \|x-P_{\bm{\Psi}^{(k)}}^Ax\|_2\leq \varepsilon_k \|b\|_2, \quad \text{and} \quad \|A^{-1}-P_{\bm{\Psi}^{(k)}}^AA^{-1}\|_2\leq \varepsilon_k.
\]
\end{thm}

Notice that the compression error $\varepsilon_k$ is in a cumulative form. However, we can restrict $\varepsilon(\mathcal{P}^{(k)})$ to increase with $k$ at certain rate, i.e. $\frac{\varepsilon(\mathcal{P}^{(k+1)})}{\varepsilon(\mathcal{P}^{(k)})} = \frac{1}{\eta}$ for some $\eta \in (0,1)$, which gives
\begin{equation}
\varepsilon_k \leq \frac{1}{1-\eta} \varepsilon(\mathcal{P}^{(k)}).
\end{equation}

With the above framework for the MMD, the original matrix $A$ can be decomposed into bounded pieces, such that the condition number $\kappa(B^{(k)})$ is controlled by choosing an appropriating partition $\mathcal{P}$ with $\varepsilon(\mathcal{P}^{(k)}) \delta(\mathcal{P}^{(k)}) \leq c$ for some constant $c$. Therefore, we can apply the MMD to solve a linear system. Notice that the difference between $\varepsilon_k$ and $\varepsilon(\mathcal{P}^{(k)})$ is very small and can be neglected, in this manuscript, we will treat $\varepsilon(\mathcal{P}^{(k)})$ as $\varepsilon_k$ and denote them simply by $\varepsilon_k$. To be coherent, we also replace the notation of $\delta(\mathcal{P}^{(k)})$ by $\delta_k$ to avoid confusion that may arise due to various notations. 

In practice, we also introduce a local approximator $\widetilde{\Psi}^{(k)}$, with which the sparsity of $\widetilde{A}^{(k)}$ and $\widetilde{B}^{(k)}$ can be preserved. In particular, we require $nnz(\tilde{A}^{(k)}) = O(nnz(A))$, where $nnz$ denotes the number of nonzero entries. We remark that, since $\widetilde{B}^{(k)}=(U^{(k)})^T\widetilde{A}^{(k-1)}U^{(k)}$, any multiplication operation concerning $\widetilde{B}^{(k)}$ only requires the applying of $(U^{(k)})^T,U^{(k)}$ and $\widetilde{A}^{(k-1)}$ separately. The applying of $(U^{(k)})^T,U^{(k)}$ can be done implicitly by performing local Householder transform with cost linear in $n$. So only the sparsity of $\widetilde{A}^{(k)}$ matters. From the estimates for the multiresolution matrix decomposition in \cite{hou2017adaptive}, we can preserve the sparsity of $\widetilde{A}^{(k)}$ by choosing the scale ratio $\eta^{-1}$ to be 
\begin{equation}
\eta^{-1} = (\log\frac{1}{\varepsilon}+\log n)^p,
\end{equation}
where we remark that $p=1$ for graph Laplacian cases. Such choice of $\eta$ also gives us the estimate of the total level number as
\begin{equation}
K=O\left(\frac{\log n}{\log(\log\frac{1}{\varepsilon}+\log n)}\right).
\label{eqt:levelnumber}
\end{equation}
Moreover, the uniform condition bound $\kappa(\mathcal{P}^{(k)},q^{(k)})\leq c$ can be imposed directly through the MMD Algorithm. For more details, please refer to Section 6 of \cite{hou2017adaptive}. For the ease of discussion in this paper, we presume using the localized decomposition to control the sparsity throughout levels and simply write $\widetilde{\psi}^{(k)}$, $\widetilde{A}^{(k)}$ and $\widetilde{B}^{(k)}$ as $\Psi^{(k)}$, $A^{(k)}$ and $B^{(k)}$.

\subsection{Implicitly Restarted Lanczos Method (IRLM)}
The Arnoldi iteration is a widely used method to find eigenvalues of unsymmetric sparse matrices. It belongs to the family of Krylov subspace methods. For symmetric case, we can further simplify it as the Lanczos iteration. A direct application of Lanczos iteration gives the largest eigenvalues of an operator by calculating the eigenvalues of its projection on a Krylov subspace. In each step the algorithm expands the Krylov subspace and finds an orthogonal basis of the space. Namely, after $k$ steps, the factorization is 
\begin{equation}
AV_k = V_k T_k + f_k e_k^T.
\label{eqt:lanczos_vector}
\end{equation}
where we recall that $T_k$ is a tridiagonal matrix when $A$ is symmetric. Denote $(\theta, y)$ as an eigenpair of $T_k$. Let $x = V_k y$. Then we have
\begin{alignat}{2}
\|Ax - x\theta\|_2 &= \|A V_k y - V_k y \theta\|_2 &= \|f_k\|_2 |e_k^T y|.
\end{alignat}
Therefore $\theta$ is a good approximation of the eigenvalue of $A$ if and only if $\|f_k\|_2 |e_k^T y|$ is small. The latter is called the Ritz residual. An analogy to the power method shows that, to compute the largest $m$ eigenvalues, the convergence rate of the largest $m$ eigenvalues of $A$ is $(\lambda_{m+1}/\lambda_m)^{k}$ where $\lambda_i$ is the $i$th largest eigenvalue of $A$. 
 
The direct Lanczos method is not practical due to the fact that $\|f_k\|_2$ rarely becomes small enough until the size of $T_k$ approaches that of $A$. An improvement is the implicitly restarted Lanczos Method (IRLM) \cite{sorensen1992implicit, lehoucq1996deflation}. The IRLM employs the idea analogous to the implicitly shifted QR-iteration \cite{francis1961}. With this approach, the ``unwanted'' eigenvalues (in this case the leftmost ones) are shifted away implicitly in each round of implicit restart, and $T_k$ is kept with a small size equal to the number of desired eigenvalues. This is one of the state-of-the-art algorithms for large-scale partial eigenproblems.

Yet, it is still complicated if we want to find the leftmost eigenvalues. One possible approach is to use a shifted IRLM. Namely, to find eigenvalues nearest to $\sigma$, we can replace $A$ with $(A-\sigma I)^{-1}$ as the target operator. By taking $\sigma=0$ we get the eigenvalues with smallest magnitude. Such approach usually converges with a few iterations, but it requires solving $A^{-1}$ in every iteration. For large sparse problems, $A^{-1}$ is usually solved by the Conjugate Gradient (CG) method. The complexity of CG is the complexity of matrix-vector product times the number of CG iterations. The former is equal to the number of nonzero entries of $A$ (denoted as $nnz(A)$), while the latter is controlled by the condition number $\kappa(A)$. Therefore, the total complexity of the shifted IRLM for solving $m_{tar}$ smallest eigenvalues is 
\begin{equation}
O(R_{\text{IRLM}} \cdot m_{tar} \cdot nnz(A) \cdot \kappa(A)),
\end{equation}
where $R_{\text{IRLM}}$ is the number of IRLM rounds. In the following, we will develop the extension-refinement algorithm to integrate the MMD framework with the shifted IRLM which gives considerable improvement in terms of iteration numbers of CG and PCG throughout the algorithm.

\begin{algorithm}[h!]
\caption{Lanczos Iteration ($p$-step extension)}
\label{alg:arnoldi_iteration}
\begin{algorithmic}[1]
\REQUIRE{$V$, $T$, $f$, target operator $op(\cdot)$, $p$.}
\ENSURE{$V$, $T$, $f$.}
\STATE{$k$ = \text{column number of }$V$;}
\FOR{$i = 1:p$}
\STATE{$\beta = \|f\|_2$;}
\IF{$\beta < \epsilon$}
\STATE{generate a new random $f$, $\beta = \|f\|_2$;}
\ENDIF
\STATE{$T \leftarrow \left(\begin{smallmatrix} T \\ \beta e^T_{k+i-1} \end{smallmatrix}\right), \quad v = f/\beta,\quad V \leftarrow [V,v]$;}
\STATE{$w = op(v)$;}
\STATE{$h = V^T w,\quad T \leftarrow [T, h]$;}
\STATE{$f = w - V h$;}
\STATE{Re-orthogonalize to adjust $f$;}
\ENDFOR
\end{algorithmic}
\end{algorithm}

\begin{algorithm}[h!]
\caption{Inner Iteration of the Implicitly Restarted Lanczos Method (IRLM)}
\label{alg:IRAM}
\begin{algorithmic}[1]
\REQUIRE{$V$, $T$, $f$.}
\ENSURE{$V$, $T$, $f$.}
\STATE{$k$ = \text{column number of }$V$;}
\STATE{Set $Q = I_{k+p}$ and $\{ \sigma_j\}$ to be the $p$ smallest eigenvalues;}
\STATE{Perform \cref{alg:arnoldi_iteration} on $V$, $T$ and $f$ for $p$ steps;}
\FOR{$j = 1:p$}
\STATE{$T - \sigma_j I = Q_j R_j$;}
\STATE{$T = Q^T_j T Q_j, Q \leftarrow QQ_j$;}
\ENDFOR
\STATE{$V \leftarrow V \cdot Q(:,1:k)$, $T \leftarrow T(1:k, 1:k)$;}
\STATE{$f \leftarrow V \cdot Q(:,k+1) \cdot T(k+1, k) + f\cdot Q(k+p,k)$;}
%\STATE{Extract $k$ eigenpairs $(\theta_i, y_i)$ and compute the Ritz residuals $\| f\|_2$;}
\end{algorithmic}
\end{algorithm}

\section{The Compressed Eigen Problem}
\label{sec:compressed_eigenproblem}
In the previous section, we introduced an effective compression technique for a SPD matrix $A$ subject to a prescribed compression error $\epsilon$. The compressed operator is also being symmetric positive definite. Therefore, by the well-known eigenvalue perturbation theory, we know that the eigenparis of the compressed operator can be used as good approximations for the eigenpairs of the original matrix. In particular, we have the following estimate:
\begin{lemma}
Let $\Theta=\Psi(\Psi^TA\Psi)^{-1}\Psi^T$ be the rank-$N$ compressed approximation of $A^{-1}$ introduced in \Cref{thm:compression_thm1} such that $\|A^{-1}-\Theta\|_2\leq \varepsilon$. Let $\mu_1\geq \mu_2\geq \cdots\geq \mu_n> 0$ be the eigenvalues of $A^{-1}$ in a descending order, and $\tilde{\mu}_1\geq \tilde{\mu}_2\geq \cdots\geq \tilde{\mu}_N> 0$ be the non-zero eigenvalues of $\Theta$ in a descending order. Then we have
\[|\mu_{i}-\tilde{\mu}_{i}|\leq \varepsilon, \quad 1\leq i\leq N; \qquad \mu_{i}\leq \varepsilon,\quad N<i\leq n.\]
Moreover, let $\tilde{v}_i$, $i=1,\cdots, N$, be the corresponding normalized eigenvectors of $\Theta$ such that $\Theta \tilde{v}_i=\tilde{\mu}_i \tilde{v}_i$, then we have 
\[\|A^{-1}\tilde{v}_i-\mu_i \tilde{v}_i\|_2\leq 2\varepsilon, \quad 1\leq i\leq N.\]
\label{lemma:approximate_eigenpair}
\end{lemma}

Since the non-zero eigenvalues of $\Theta$ and the corresponding eigenvectors actually result from the non-singular stiffness matrix $A_{st}=\Psi^TA\Psi$, we will call these eigenpairs the \textbf{essential eigenpairs} of $\Theta$ in what follows. We will also need the following lemma for developing our algorithms.
\begin{lemma}
Let $(\tilde{\mu}_i, \tilde{v}_i)$, $i=1,\cdots,N$, be the $N$ essential eigenpairs of $\Theta$ given in \Cref{lemma:approximate_eigenpair}.
\label{lemma:generaleigenproblem}
\begin{itemize}
\item[(i)] Let $w_i=\Psi^T \tilde{v}_i$, then
\[\Psi^T\Psi A_{st}^{-1} w_i=\tilde{\mu}_iw_i, \quad 1\leq i\leq N.\]
\item[(ii)] Let $z_i=\Psi^\dagger \tilde{v}_i=(\Psi^T\Psi)^{-1}\Psi^T \tilde{v}_i$, then 
\[A_{st}^{-1}\Psi^T\Psi z_i=\tilde{\mu}_iz_i, \quad 1\leq i\leq N,\]
\end{itemize}
where $A_{st}=\Psi^TA\Psi$ is the stiffness matrix. Conversely, if either (i) or (ii) is true, then $(\tilde{\mu}_i, \tilde{v}_i)$, $i=1,\cdots,N$, are eigenpairs of $\Theta$.
\end{lemma}

Similar to \Cref{lemma:approximate_eigenpair}, we have the following estimates for multiresolution decomposition.

\begin{lemma} 
Given an integer $K$, let $\Theta^{(k)}=\bm{\Psi}^{(k)}\big((\bm{\Psi}^{(k)})^TA\bm{\Psi}^{(k)}\big)^{-1}(\bm{\Psi}^{(k)})^T$, $k=1,2,\cdots,K$, with $\bm{\Psi}^{(k)}$ given in \Cref{eqt:mutliresolution_phi&U&psi}. Write $A^{-1}=\Theta^{(0)}$. Let $(\mu_i^{(k)},v_i^{(k)})$, $i=1,2,\cdots,N^{(k)}$, be the essential eigenpairs of $\Theta^{(k)}$ where $\mu_1^{(k)}\geq \mu_2^{(k)}\geq \cdots\geq \mu_{N^{(k)}}^{(k)}> 0$. Then for any $0\leq k'<k\leq K$, we have
\[|\mu_i^{(k')}-\mu_i^{(k)}|\leq \varepsilon_k,\quad 1\leq i \leq N^{(k)};\qquad |\mu_i^{(k')}|\leq \varepsilon_k,\quad N^{(k)}< i \leq N^{(k')},\]
and
\[\|\Theta^{(k')}v_i^{(k)}-\mu_i^{(k')}v_i^{(k)}\|_2\leq 2\varepsilon_k,\quad 1\leq i \leq N^{(k)}.\]
\label{lemma:approximate_eigenpair_multiresolution}
\end{lemma}
\begin{proof}
By \Cref{thm:error_accumulate} we have that $\|\Theta^{(0)}-\Theta^{(k)}\|_2=\|A^{-1}-\Theta^{(k)}\|_2\leq \varepsilon_k$, $k=1,2,\cdots,K$. From the definition of $\Theta^{(k)}$ and the decomposition $\cref{eqt:MMD}$, one can easily check that 
\[A^{-1}=\Theta^{(0)}\succeq\Theta^{(1)}\succeq\cdots \succeq\Theta^{(K-1)}\succeq\Theta^{(K)}.\]
Then the results follow immediately.
\end{proof}

\subsection*{On Compressed Eigenproblems} 
We should remark that the efficiency of constructing the compressed operator we propose relies on the exponential decay property of the basis $\Psi$. This spacial exponential decay feature allows us to localize $\Psi$ and to construct sparse stiffness matrix $A_{st}=\Psi^TA\Psi$ without compromising compression accuracy $\varepsilon$ in $O(nnz(A) \cdot (\log(\frac{1}{\varepsilon})+\log n)^c)$ time. In fact, the problem of using spatially localized/compact basis to compress high dimensional operator and to approximate eigenspace of smallest eigenvalues has long been studied in different ways. A representative pioneer work is the method of compressed modes proposed by Ozoliņ\v{s} et al. \cite{ozolicnvs2013compressed}, intended originally for Schr\"odinger’s equation in quantum physics. By adding a $L_1$ regularization to the variational form of an eigenproblem, they obtained spatially compressed basis modes that well span the desired eigenspace. Though the way they obtain sparsity is quite different from what we do, both methods obtain interestingly similar results for some model problems. It can be inspiring to make comparison between their method and ours, so that readers can have better understanding of our approach. We leave the detailed comparison to the Appendix.

\section{Hierarchical Spectrum Completion} 
\label{sec:spectrum_extension}
Now that we have a sequence of compressed approximations, we next seek to use this decomposition to compute the dominant spectrum of $A^{-1}$ down to a prescribed value in a hierarchical manner. In particular, we propose to decompose the target spectrum into several segments of different scales, and then allocate the computation of each segment to a certain level of the compressing sequence so that the problem on each level is well-conditioned.

To implement this idea, we first go back to the one-level compression settings. Suppose that we have accurately obtained the first $m$ essential eigenpairs $(\mu_i,v_i)$, $i=1,\cdots,m$, of $\Theta=\Psi(\Psi^TA\Psi)^{-1}\Psi^T=\Psi A_{st}^{-1}\Psi^T$, and our aim is to compute the following $\hat{m}-m$ eigenpairs (namely extend to the first $\hat{m}$ eigenpairs) using the Lanczos method. Define $V_m=\Span\{v_i:1\leq i\leq m\}$ and $V_{m^+}=\Span\{v_i:m<i\leq N\}=V_m^{\perp}\cap \Span\{\Psi\}$. Then to perform the Lanczos method to compute the next segment of eigenpairs of $\Theta$, we need to repeatedly apply the operator $\Psi A_{st}^{-1}\Psi^T$ to vectors in $V_{m^+}$, which requires to compute $A_{st}^{-1}w$ for $w\in W_{m^+}=\Psi^T(V_{m^+})$.

Ideally we want the computation of the following $\hat{m}-m$ eigenpairs to be restricted to a problem with bounded spectrum width that is proportional to $\mu_{m}/\mu_{\hat{m}}$. This is possible since we assume that we have accurately obtained the span space $V_m$ of the first $m$ eigenvectors, and thus we can consider our problem in the reduced space orthogonal to $V_m$. In this case, the CG method will be efficient for computing inverse matrix operations.

\begin{definition}
Let $A$ be a symmetric, positive definite matrix, and $V$ be an invariant subspace of $A$. We define the condition number of $A$ with respect to $V$ as
\[\kappa(A,V)=\frac{\lambda_{max}(A,V)}{\lambda_{min}(A,V)},\]
where 
\[\lambda_{max}(A,V)=\max_{v\in V, v\neq0}\frac{v^TAv}{v^Tv},\qquad \lambda_{min}(A,V)=\min_{v\in V, v\neq0}\frac{v^TAv}{v^Tv}.\]
\end{definition}

\begin{thm}
Let $A$ be a symmetric, positive definite matrix, and $V$ be an invariant subspace of $A$. When using the conjugate gradient method to solve $Ax=b$ with initial guess $x_0$ such that $r_0=b-Ax_0\in V$, we have the following estimate 
\[\|x_k-x_*\|_A\leq 2\Big(\frac{\sqrt{\kappa(A,V)}-1}{\sqrt{\kappa(A,V)}+1}\Big)^k\|x_0-x_*\|_A,\]
and
\[\|x_k-x_*\|_2\leq 2\sqrt{\kappa(A,V)}\Big(\frac{\sqrt{\kappa(A,V)}-1}{\sqrt{\kappa(A,V)}+1}\Big)^k\|x_0-x_*\|_2,\]
where $x_*$ is the exact solution, and $x_k\in x_*+V$ is the solution at the $k^{th}$ step of CG iteration. Thus it takes $k=O(\kappa(A,V)\cdot\log\frac{1}{\epsilon})$ steps (or $k=O\big(\kappa(A,V)\cdot(\log\kappa(A,V)+\log\frac{1}{\epsilon})\big)$ steps) to obtain a solution subject to relative error $\epsilon$ in the energy norm (or $l_2$ norm). 
\label{thm:CG_estimate}
\end{thm}
\begin{proof} We only need to notice that the $k$-order Krylov subspace $\mathcal{K}(A, r_0, k)$ generated by $A$ and $r_0$ satisfies
\[\mathcal{K}(A,r_0,k)\subset V,\quad \forall k\in \mathbb{Z}.\]
\end{proof}

Notice that, for any $i=m+1,\cdots, N$, though $v_i\in V_{m^+}$ is an eigenvector of $\Theta=\Psi A_{st}^{-1}\Psi^T$, $w_i=\Psi^Tv_i$ is not an eigenvector of $A_{st}^{-1}$(but an eigenvector of $\Psi^T\Psi A_{st}^{-1}$) since we do not require $\Psi$ to be orthonormal. Therefore the space $W_{m^+}$ is not an invariant space of $A_{st}$, and if we directly use the CG method to solve $A_{st}x=w$, the convergence rate will depend on $\kappa(A_{st})$, instead of $\kappa(A_{st})/\mu_{m}$ as intended. Though we bound $\lambda_{max}(A_{st})$ from above by $\delta(\mathcal{P})$ and $\lambda_{min}(A_{st})$ from below by $\lambda_{min}(A)$ (See \Cref{thm:condition_number}), $\kappa(A_{st})$ can be still large since we prescribe a bounded compression rate in practice to ensure the efficiency of the compression algorithm.

Therefore, we need to find a proper invariant space, so that we can make use of the knowledge of the space $V_m$ and restrict the computation of $A_{st}^{-1}w$ to a problem of narrower spectrum.

\begin{lemma} 
Let $(\mu_i,v_i)$, $i=1,\cdots, N$, be the essential eigenpairs of $\Theta=\Psi(\Psi^TA\Psi)^{-1}\Psi^T=\Psi A_{st}^{-1}\Psi^T$, such that $\mu_1\geq \mu_2\geq \cdots\geq \mu_N>0$. Let $(\Psi^T\Psi)^\frac{1}{2}$ be the square root of the symmetric, positive definite matrix $\Psi^T\Psi$. Then $(\mu_i,z_i)$, $i=1,\cdots, N$, are all eigenpairs of $(\Psi^T\Psi)^{\frac{1}{2}}A_{st}^{-1}(\Psi^T\Psi)^{\frac{1}{2}}$, where
\[z_i=(\Psi^T\Psi)^{-\frac{1}{2}}\Psi^Tv_i, \quad 1\leq i\leq N.\]
Moreover, for any subset $S\subset \{1,2,\cdots,N\}$, and $Z_S=\Span\{z_i: i\in S\}$, we have
\[\mathcal{K}(A_{\Psi}, z, k)\subset Z_S, \quad \forall z\in Z_S, \ \forall k\in \mathbb{Z},\]
where $A_{\Psi}=(\Psi^T\Psi)^{-\frac{1}{2}}A_{st}(\Psi^T\Psi)^{-\frac{1}{2}}$.
\label{lemma:eigenpreconditioner}
\end{lemma}

\begin{lemma}
Let $\Psi$ be given in \Cref{eqt:psi_closedform}, then we have
\[\lambda_{min}(\Psi^T\Psi)\geq 1,\quad \lambda_{max}(\Psi^T\Psi)\leq 1+\varepsilon(\mathcal{P}) \delta(\mathcal{P}),\]
and thus
\[\kappa(\Psi^T\Psi)\leq 1+\varepsilon(\mathcal{P}) \delta(\mathcal{P}).\]
\label{lemma:psipsi_conditionnumb}
\end{lemma}

\begin{proof}
Let $U$ be the orthogonal complement basis of $\Phi$ given in \Cref{eqt:Ainv_representation}, so $[\Phi,U]$ is an orthonormal basis of $\mathbb{R}^n$, and we have $\Phi\Phi^T+UU^T = I_n$. Since $\Phi^T\Psi = \Phi^TA^{-1}\Phi(\Phi^TA^{-1}\Phi)^{-1} = I_{N}$, we have
\[\Psi^T\Psi = \Psi^T\Phi\Phi^T\Psi + \Psi^TUU^T\Psi = I_N + \Psi^TUU^T\Psi. \]
We then immediately obtain $\Psi^T\Psi\succeq I_N$, and thus $\lambda_{min}(\Psi^T\Psi)\geq 1$. To obtain an upper bound of $\lambda_{max}(\Psi^T\Psi)$, we notice that from the construction of $\Phi$ we have
\[\|x-P_{\Phi}x\|_2^2 \leq \varepsilon(\mathcal{P}) x^TAx,\quad \forall x\in \mathbb{R}^n\quad  \Longrightarrow\quad (I_n-P_{\Phi})^2 \preceq \varepsilon(\mathcal{P})A, \]
where $P_{\Phi}=\Phi\Phi^T$ denotes the orthogonal projection into $\Span\{\Phi\}$. Since $\Phi\Phi^T+UU^T = I_n$, we have 
\[UU^T = I_n-\Phi\Phi^T = (I_n-\Phi\Phi^T)^2 \preceq \varepsilon(\mathcal{P})A.\]
Therefore we have
\[\Psi^T\Psi= I_N + \Psi^TUU^T\Psi\preceq I_N + \varepsilon(\mathcal{P})\Psi^TA\Psi= I_N+\varepsilon(\mathcal{P})A_{st},\]
and by \Cref{thm:condition_number} we obtain 
\[\lambda_{max}(\Psi^T\Psi)\leq 1 + \varepsilon(\mathcal{P})\lambda_{max}(A_{st})\leq 1+\varepsilon(\mathcal{P}) \delta(\mathcal{P}).\]
\end{proof} 

\begin{thm}
Let $A_{\Psi}$ and $(\mu_i,z_i)$ be defined as in \Cref{lemma:eigenpreconditioner}. Let $Z_{m^+}=\Span\{ z_i: m<i\leq N\}$, then $Z_{m^+}$ is an invariant space of $A_{\Psi}$, and we have
\[\kappa(A_\Psi,Z_{m^+})\leq\mu_{m+1}\delta(\mathcal{P}).\]
\label{thm:Apsi_conditionbound}
\end{thm}
\begin{proof}
By \Cref{lemma:psipsi_conditionnumb}, we have
\[\lambda_{max}(A_\Psi,Z_{m^+})\leq \lambda_{max}(A_\Psi)=\|(\Psi^T\Psi)^{-\frac{1}{2}}A_{st}(\Psi^T\Psi)^{-\frac{1}{2}}\|_2\leq \|A_{st}\|_2\|(\Psi^T\Psi)^{-1}\|_2\leq \delta(\mathcal{P)}.\]
And by the definition of $Z_{m^+}$, we have
\[\lambda_{min}(A_\Psi,Z_{m^+})=\frac{1}{\lambda_{max}(A_\Psi^{-1},Z_{m^+})}=\frac{1}{\lambda_{max}((\Psi^T\Psi)^{\frac{1}{2}}A_{st}^{-1}(\Psi^T\Psi)^{\frac{1}{2}},Z_{m^+})}=\frac{1}{\mu_{m+1}}.\]
\end{proof}

Inspired by \Cref{lemma:psipsi_conditionnumb} and \Cref{thm:Apsi_conditionbound}, we now consider to solve $A_{st}x=w$ efficiently for $w\in W_{m^+}=\Psi^T(V_{m^+})=(\Psi^T\Psi)^{\frac{1}{2}}(Z_{m^+})$ by making use of the controlled condition number $\kappa(A_\Psi,Z_{m^+})$ and $\kappa(\Psi^T\Psi)$. Theoretically, we can compute $x=A_{st}^{-1}w$ by the following steps:
\begin{itemize}
\item[(i)] Compute $b=(\Psi^T\Psi)^{-\frac{1}{2}}w\in Z_{m^+}$;
\item[(ii)] Use the CG method to compute $y=A_{\Psi}^{-1}b$ with initial guess $y_0$ such that $b-A_\Psi y_0\in Z_{m^+}$;
\item[(iii)] Compute $x=(\Psi^T\Psi)^{-\frac{1}{2}}y$.
\end{itemize}
Notice that this procedure is exactly solving $A_{st}x=w$ using the preconditioned CG method with preconditioner $\Psi^T\Psi$, which only involves applying $A_{st}$ and $(\Psi^T\Psi)^{-1}$ to vectors, but still enjoys the good conditioning property of $A_\Psi$ restricted to $Z_{m^+}$. Therefore we have the following estimate:\\

\begin{corollary} Consider using the PCG method to solve $A_{st}x=w$ for $w\in W_{m^+}$ with preconditioner $\Psi^T\Psi$ and initial guess $x_0$ such that $r_0=w-A_{st}x_0\in W_{m^+}$. Let $x_*$ be the exact solution, and $x_k$ be the solution at the $k^{th}$ step of the PCG iteration. Then we have
\[\|x_k-x_*\|_{A_{st}}\leq 2\Big(\frac{\sqrt{\kappa(A_\Psi,Z_{m^+})}-1}{\sqrt{\kappa(A_\Psi,Z_{m^+})}+1}\Big)^k\|x_0-x_*\|_{A_{st}},\]
and 
\[\|x_k-x_*\|_2\leq 2\sqrt{\kappa(\Psi^T\Psi)\kappa(A_\Psi,Z_{m^+})}\Big(\frac{\sqrt{\kappa(A_\Psi,Z_{m^+})}-1}{\sqrt{\kappa(A_\Psi,Z_{m^+})}+1}\Big)^k\|x_0-x_*\|_2.\]
\label{cor:PCG_convergence}
\end{corollary}

\begin{proof} Let $y_k=(\Psi^T\Psi)^{\frac{1}{2}}x_k$ and $y_*=(\Psi^T\Psi)^{\frac{1}{2}}x_*$, then we have
\[\|y_k-y_*\|_2^2=(x_k-x_*)^T\Psi^T\Psi(x_k-x_*),\]
and 
\[\|y_k-y_*\|_{A_\Psi}^2=(y_k-y_*)^TA_\Psi(y_k-y_*)=\|x_k-x_*\|_{A_{st}}^2.\] 
Noticing that $(\Psi^T\Psi)^{-\frac{1}{2}}r_0\in Z_{m^+}$ and $\mathcal{K}(A_\Psi,(\Psi^T\Psi)^{-\frac{1}{2}}r_0,k)\subset Z_{m^+}\ \forall k$, the results follow from \Cref{thm:CG_estimate}.
\end{proof}

By \Cref{cor:PCG_convergence}, to compute a solution of $A_{st}x=w$ subject to a relative error $\epsilon$ in the $A_{st}$-norm, the number of needed PCG iterations is 
\[O\big(\kappa(A_\Psi,Z_{m^+})\cdot\log\frac{1}{\epsilon}\big)= O\big(\mu_{m+1}\delta(\mathcal{P})\cdot\log\frac{1}{\epsilon}\big).\] 
This is also an estimate of the number of needed PCG iterations for a relative error $\epsilon$ in the $l_2$-norm, if we assume that $\kappa(\Psi^T\Psi),\kappa(A_\Psi,Z_{m^+}) \leq \frac{1}{\epsilon}$.

In what follows we will denote $M=\Psi^T\Psi$. Notice that the nonzero entries of $M$ are due to the overlapping support of column basis vectors of $\Psi$, while the nonzero entries of $A_{st}=\Psi^TA\Psi$ are results of interactions between column basis vectors of $\Psi$ with respect to $A$. Thus we can reasonably assume that $nnz(M)\leq nnz(A_{st})$. Suppose that in each iteration of the whole PCG procedure, we also use the CG method to solve for $M^{-1}$ subject to a relatively higher precision $\hat{\epsilon}$, which requires a cost of $O(nnz(M)\cdot\kappa(M)\cdot\log\frac{1}{\hat{\epsilon}})$. In practice it is sufficient to take $\hat \epsilon$ smaller than but comparable to $\epsilon$ (e.g. $\hat \epsilon =0.1\epsilon$), so $\log(\frac{1}{\hat\epsilon}) = O(\log\frac{1}{\epsilon})$. By \Cref{lemma:psipsi_conditionnumb} we have $\kappa(M)=O(\varepsilon(\mathcal{P})\delta(\mathcal{P}))$. Then the computational complexity of each single iteration can be bounded by 
\[O\big(nnz(A_{st})\big)+O\big(nnz(M)\cdot\kappa(M)\cdot\log\frac{1}{\epsilon}\big)= O\big(nnz(A_{st})\cdot \varepsilon(\mathcal{P})\delta(\mathcal{P}) \cdot\log\frac{1}{\epsilon}\big),\]
and the total cost of computing a solution of $A_{st}x=w$ subject to a relative error $\epsilon$ is 
\begin{equation}
O\big(\mu_{m+1}\delta(\mathcal{P})\cdot nnz(A_{st})\cdot \varepsilon(\mathcal{P})\delta(\mathcal{P}) \cdot(\log\frac{1}{\epsilon})^2\big).
\label{eqt:extend_complexity}
\end{equation}

We remark that when the original size of $A\in\mathbb{R}^{n\times n}$ is large, the eigenvectors $V$ are long and dense. It would be expensive to compute inner products with these long vectors over and over again. In fact, in the previous discussions the operator $\Theta=\Psi A_{st}^{-1}\Psi^T$(of the same size as $A$) and the eigenvectors $V$ are only for purpose of analysis use to explain the idea of our method. In practical, for a long vector $v=\Psi \hat{v}$, we don't need to keep track of the whole vector, but only need to store its much shorter coefficients $\hat{v}$ of compressed dimension $N$ instead. When we compute $v_2=\Theta v_1=\Psi A_{st}^{-1}\Psi^T v_1$, it is equivalent to computing $\hat{v}_2=A_{st}^{-1}M \hat{v}_1$, where $v_j=\Psi \hat{v}_j,\ j=1,2$, and $M=\Psi^T\Psi$. One can check that the analysis presented above still applies. So in the implementation of our method, we only deal with operator $A_{st}^{-1}M$ and short vectors $\widehat{V}$, and the long eigenvectors $V$ and $\Psi$ will not appear until in the very end when we recover $V=\Psi\widehat{V}$. We remark that since the eigenvectors of $\Theta$ are orthogonal, their coefficient vectors $\widehat{V}$ are $M$-orthogonal, i.e. $\widehat{V}^TM\widehat{V}=I$. We use $\|x\|_M$ to denote the norm $\sqrt{x^TMx}$.

Recall that in the Lanczos method with respect to operator $\Theta$, the upper-Hessenberg matrix $T$ in the Arnoldi relation  
\[\Theta V = VT+fe^T\]
is indeed tridiagonal, since $\Theta$ is symmetric, and $V^T[V,f]=[I,\bf 0]$. This upper-Hessenberg matrix $T$ being tridiagonal is the reason why the implicit restarting process(\Cref{alg:IRAM}) is efficient. Now since we are actually dealing with the operator $A_{st}^{-1}M$ and the coefficient vectors $\widehat{V}^T=M^{-1}\Psi^T V$, the Arnoldi relation becomes 
\[A_{st}^{-1}M \widehat{V} = \widehat{V} T +\hat{f}e^T,\]
where $\hat{f}=M^{-1}\Psi^T f$. So as long as we keep $\widehat{V}$ $M$-orthogonal and $f$ $M$-orthogonal to $\widehat{V}$, $T$ will still be tridiagonal since
\[T=\widehat{V}^TM\widehat{V}T=\widehat{V}^TM(\widehat{V}T+\hat{f}e^T )=\widehat{V}^TMA_{st}^{-1}M \widehat{V}\]
is symmetric. We therefore modified \Cref{alg:arnoldi_iteration} to \Cref{alg:arnoldi_iteration_M} to take $M$-orthogonality into consideration.

Summarizing the analysis above, we propose \Cref{alg:eigenpair_extension} for extending a given collection of eigenpairs using the Lanczos type method. The operator $OP(\ \cdot\ ;A_{st},M,\epsilon_{op})$ exploits our key idea that uses $M=\Psi^T\Psi$ as the preconditioner to effectively reduce the number of PCG iterations in every operation of $A_{st}^{-1}M$. For convenience, we will use ``$x=pcg(A,b,M,x_0,\epsilon)$'' to represent the operation of computing $x=A^{-1}b$ using the PCG method with preconditioner $M$ and initial guess $x_0$, subject to relative error $\epsilon$. ``$x=pcg(A,b,-,x_0,\epsilon)$'' means no preconditioner is used (i.e. the normal CG method), and ``$x=pcg(A,b,M,-,\epsilon)$'' means an all zero vector is used as the initial guess.
% \begin{table}[H]
% \centering
% \renewcommand{\arraystretch}{1.3}
%     \begin{tabular}{|ll|}
%     \hline
%        &\textbf{Operator $OP(\ \cdot\ ;A_{st},M,\Psi,\epsilon_{op})$}\\ \hline
%         & \textbf{input:} $v$\\
%         & \textbf{output:} $y$ \\ 
%         1. & $w=\Psi^Tv$ \\
%         2. & $\tilde{w}=pcg(A_{st},w,M,-,\epsilon_{op})$\\
%         3. & $y=\Psi\tilde{w}$\\
%     \hline 
%     \end{tabular}
% \end{table}

\begin{algorithm}[h!] 
\caption{General Lanczos Iteration ($p$-step extension)}
\label{alg:arnoldi_iteration_M}
\begin{algorithmic}[1]
\REQUIRE{$\widehat{V}$, $T$, $\hat{f}$, target operator $op(\cdot)$, $p$, inner product matrix $M$}
\ENSURE{$\widehat{V}$, $T$, $\hat{f}$}
\STATE{$k$ = \text{column number of }$\widehat{V}$;}
\FOR{$i = 1:p$}
\STATE{$\beta = \|\hat{f}\|_M$;}
\IF{$\beta < \epsilon$}
\STATE{generate a new random $\hat{f}$, $\beta = \|\hat{f}\|_M$;}
\ENDIF
\STATE{$T \leftarrow \left(\begin{smallmatrix} T \\ \beta e^T_{k+i-1} \end{smallmatrix}\right), \quad \hat{v} = \hat{f}/\beta,\quad \widehat{V} \leftarrow [\widehat{V},\hat{v}]$;}
\STATE{$w = op(\hat{v})$;}
\STATE{$h = \widehat{V}^T Mw,\quad T \leftarrow [T, h]$;}
\STATE{$\hat{f} = w - \widehat{V} h$;}
\STATE{Re-orthogonalize to adjust $f$(with respect to $M$-orthogonality);}
\ENDFOR
\end{algorithmic}
\end{algorithm} 

\begin{algorithm}
\renewcommand{\thealgorithm}{}
\floatname{algorithm}{}
\caption{Function $y$ = \textbf{Operator} $OP(x;A_{st},M,\epsilon_{op})$}
\begin{algorithmic}[1]
\STATE{$w=Mx$;}
\STATE{$y=pcg(A_{st},w,M,-,\epsilon_{op})$;}
\end{algorithmic}
\end{algorithm}

\addtocounter{algorithm}{-1}
\begin{algorithm}[h!]
\caption{Eigenpair Extension}
\label{alg:eigenpair_extension}
\begin{algorithmic}[1]
\REQUIRE{$\widehat{V}_{ini}$, $D_{ini}$, $OP(\ \cdot\ ; A_{st},M,\epsilon_{op})$, target number $m_{tar}$,\\prescribed accuracy $\epsilon$, eigenvalue threshold $\mu$, searching step $d$.}
\ENSURE{$\widehat{V}_{ex}$, $D_{ex}$.}
\STATE{Generate random initial vector $\widehat{V}=\hat{v}$ that is $M$-orthogonal to $\widehat{V}_{ini}$;}
\REPEAT
\STATE{perform $d$ steps of general Lanczos iteration (\Cref{alg:arnoldi_iteration_M}) with operator $OP$ to extend $\widehat{V},T$;}
\WHILE{Lanczos residual $> \epsilon$,}
\STATE{Perform $c\cdot d$ steps of shifts to restart Lanczos (\Cref{alg:IRAM}) and renew $\widehat{V},T$;}
\ENDWHILE
\STATE{Find the $d^{th}$ smallest eigenvalue of $T$ as $\hat{\mu}$;}
\UNTIL{$\hat{\mu} < \mu$ \text{ or } $\dim(\widehat{V}) \geq m_{tar}-\dim(\widehat{V}_{ini})$.}
\STATE{$m_{new}=\dim(\widehat{V})$;}
\WHILE{Lanczos residual $> \epsilon$,}
\STATE{Perform $c\cdot m_{new}$ steps of shifts to restart Lanczos (\Cref{alg:IRAM}) and renew $\widehat{V},T$;}
\ENDWHILE
\STATE{$PSP^T = T$ (Schur Decomposition);}
\STATE{$\widehat{V}_{ex}=[\widehat{V}_{ini},\widehat{V}P]$,\quad $D_{ex}=\left[\begin{array}{cc} D_{ini} & \\ & S\end{array}\right]$;}
\end{algorithmic}
\end{algorithm}

% \begin{table}[H]
% \centering
% \renewcommand{\arraystretch}{1.3}
%     \begin{tabular}{|ll|}
%     \hline
%        &\textbf{Algorithm: Eigenpair Extension}\\ \hline
%         & \textbf{input:} $V_{ini}$, $D_{ini}$, $OP(\ \cdot\ ; A_{st},M,\Psi,\epsilon_{op})$, target number $m_{tar}$\\ 
%         & \qquad\quad precribed accuracy $\epsilon$, eigenvalue threshold $\mu$, searching step $d$\\
%         & \textbf{output:} $V_{ex}$, $D_{ex}$ \\ 
%        1.&  generate random initial vector $V=v$ in the space $V_{ini}^\perp$ \\ 
%        2.&  \textbf{do} \\
%        3.&  \qquad perform $d$ steps of Alnordi iteration with operator $OP$ to extend $V,H$\\
%        4.&  \qquad \textbf{while} Alnordi error $>\epsilon$ \\ 
%        5.&  \qquad\qquad perform $c\cdot d$ steps of shifts to restart Alnordi and renew $V,H$\\
%        6.&  \qquad \textbf{end while} \\
%        7.&  \qquad find the $d^{th}$ smallest eigenvalue of $H$ as $\hat{\mu}$ \\
%        8.&  \textbf{until} $\hat{\mu}<\mu$ \textbf{or} $\dim(V)\geq m_{tar}-\dim(V_{ini})$\\
%        9.&  $m_{new}=\dim(V)$ \\
%        10.&  \textbf{while} Alnordi error $>\epsilon$ \\
%        11.& \qquad perform $c\cdot m_{new}$ steps of shifts to restart Alnordi and renew $V,H$\\
%        12.& \textbf{end while} \\
%        13.& $PSP^T=H$\quad (Schur decomposition)\\
%        14.& $V_{ex}=[V_{ini},VP]$,\quad $D_{ex}=\left[\begin{array}{cc} D_{ini} & \\ & S\end{array}\right]$\\ 
%     \hline 
%     \end{tabular}
% \end{table}
%
Given an existing eigenspace $V_{ini}=\Psi\widehat{V}_{ini}$, \Cref{alg:eigenpair_extension} basically uses the Lanczos method to find the following eigenpairs of $\Theta$ in the space $V_{ini}^\perp$. Notice that the output $\widehat{V}_{ex}$ gives the coefficients of the desired eigenvectors $V_{ex}$ in the basis $\Psi$. However, different from the classical Lanczos method, we do not prescribe a specific number for the output eigenpairs. Instead, we set a threshold $\mu$ to bound the last output eigenvalue. As we will develop our idea into a multi-level algorithm that pursues a number of target eigenpairs hierarchically, the output of the current level will be used to generate the initial eigenspace for the higher level. Therefore, the purpose of setting a threshold $\mu$ on the current level is to bound the restricted condition number on the higher level, as the initial eigenspace $V_{ini}$ from the lower level helps to bound the restricted condition number on the current level.

The choice of the threshold $\mu$ will be discussed in detail after we introduce the refinement procedure. Here, to develop a hierarchical spectrum completion method using the analysis above, we state the hierarchical versions of \Cref{lemma:psipsi_conditionnumb} and \Cref{thm:Apsi_conditionbound}.

\begin{lemma}
Let $\bm{\Psi}^{(k)}$ be given in \Cref{eqt:mutliresolution_phi&U&psi}, and $M^{(k)}=(\bm{\Psi}^{(k)})^T\bm{\Psi}^{(k)}$. Then we have 
\[\lambda_{min}(M^{(k)})\geq 1,\quad \lambda_{max}(M^{(k)})\leq 1+\varepsilon_k\delta_k,\]
and thus
\[\kappa(M^{(k)})\leq 1+\varepsilon_k\delta_k.\]
\label{lemma:multi_psipsi_conditionnumb}
\end{lemma}
\begin{proof}
The proof is similar to the proof of \Cref{lemma:psipsi_conditionnumb}. Let $\bm{U}^{(k)}=(\bm{\Phi}^{(k)})^{\perp}$ be the orthogonal complement basis of $\bm{\Phi}^{(k)}$. According to \Cref{thm:error_accumulate}, we have 
\[\|x-P_{\bm{\Phi}^{(k)}}x\|_2^2\leq\varepsilon_k\|x\|_A^2, \]
which implies that
\[\bm{U}^{(k)}(\bm{U}^{(k)})^T = (I_n - \bm{\Phi}^{(k)}(\bm{\Phi}^{(k)})^T) \leq \varepsilon_k A.\] 
Notice that $(\bm{\Phi}^{(k)})^T\bm{\Psi}^{(k)}=I_{N^{(k)}}$, $\bm{\Phi}^{(k)}(\bm{\Phi}^{(k)})^T+\bm{U}^{(k)}(\bm{U}^{(k)})^T=I_n$, we thus have 
\[M^{(k)} = (\bm{\Psi}^{(k)})^T\bm{\Phi}^{(k)}(\bm{\Phi}^{(k)})^T\bm{\Psi}^{(k)}+(\bm{\Psi}^{(k)})^T\bm{U}^{(k)}(\bm{U}^{(k)})^T\bm{\Psi}^{(k)}= I_{N^{(k)}}+(\bm{\Psi}^{(k)})^T\bm{U}^{(k)}(\bm{U}^{(k)})^T\bm{\Psi}^{(k)},\]
\[\Longrightarrow\qquad  I_{N^{(k)}} \preceq M^{(k)} \preceq I_{N^{(k)}} + \varepsilon_k (\bm{\Psi}^{(k)})^TA\bm{\Psi}^{(k)}=I_{N^{(k)}} + \varepsilon_k A^{(k)}.\]
Therefore we have $\lambda_{min}(M^{(k)})\geq 1$, and by \Cref{cor:Bconditonnumber} we have
\[\lambda_{max}(M^{(k)})\leq 1+\varepsilon_k\lambda_{max}(A^{(k)})\leq 1+\varepsilon_k\delta_k.\]
\end{proof}

\begin{thm}
Let $A^{(k)}$ and $\bm{\Psi}^{(k)}$ be given in \Cref{eqt:mutliresolution_phi&U&psi}, and $M^{(k)}=(\bm{\Psi}^{(k)})^T\bm{\Psi}^{(k)}$. Let $(\mu_i^{(k)},v_i^{(k)})$, $i=1,\cdots,N^{(k)}$, be the essential eigenpairs of $\Theta^{(k)}=\bm{\Psi}^{(k)}(A^{(k)})^{-1}(\bm{\Psi}^{(k)})^T$. Define 
\[z_i^{(k)}=(M^{(k)})^{-\frac{1}{2}}(\bm{\Psi}^{(k)})^Tv_i^{(k)},\quad 1\leq i\leq N^{(k)}.\]
Given an integer $m_{k}$, let $Z_{m_k^+}^{(k)}=\Span\{ z_i^{(k)}: m_k<i\leq N^{(k)}\}$, then $Z_{m_k^+}^{(k)}$ is an invariant space of $A^{(k)}_{\bm{\Psi}}=(M^{(k)})^{-\frac{1}{2}}A^{(k)}(M^{(k)})^{-\frac{1}{2}}$, and we have
\[\kappa(A^{(k)}_{\bm{\Psi}},Z_{m^+}^{(k)})\leq\mu_{m_k+1}^{(k)}\delta_k.\]
Moreover, consider using the PCG method to solve $A^{(k)}x=w$ for $w\in W_{m_k^+}^{(k)}$ with preconditioner $M^{(k)}$ and initial guess $x_0$ such that $r_0=w-A^{(k)}x_0\in W_{m_k^+}^{(k)}$, where $W_{m_k^+}^{(k)}=\Span\{ (\bm{\Psi}^{(k)})^Tv_i^{(k)}: m_k<i\leq N^{(k)}\}$. Let $x_*$ be the exact solution, and $x_t$ be the solution at the $t^{th}$ step of the PCG iteration. Then we have
\[\|x_t-x_*\|_{A^{(k)}}\leq 2\Big(\frac{\sqrt{\mu_{m_k+1}^{(k)}\delta_k}-1}{\sqrt{\mu_{m_k+1}^{(k)}\delta_k}+1}\Big)^t\|x_0-x_*\|_{A^{(k)}},\]
and 
\[\|x_t-x_*\|_2\leq 2\sqrt{\varepsilon_k\mu_{m_k+1}^{(k)}\delta_k^2}\Big(\frac{\sqrt{\mu_{m_k+1}^{(k)}\delta_k}-1}{\sqrt{\mu_{m_k+1}^{(k)}\delta_k}+1}\Big)^t\|x_0-x_*\|_2.\]
\label{thm:multi_spectrum_completion_overall}
\end{thm}

Recall that we will use the CG method to implement Lanczos iteration on each level $k$ to complete the target spectrum. To ensure the efficiency of the CG method, namely to bound the restricted condition number $\kappa(A^{(k)}_{\bm{\Psi}},Z_{m^+}^{(k)})$ on each level, we need a priori knowledge of the spectrum $\{(\mu_i^{(k)},v_i^{(k)}): 1\leq i\leq m_k\}$ such that $\mu_{m_k+1}^{(k)}\delta_k$ is uniformly bounded. This given spectrum should be inductively computed on the lower level $k+1$. But notice that there is a compression error between each two neighbour levels, which will compromise the orthogonality and thus the theoretical bound for restricted condition number, if we directly use the spectrum of the lower level as a priori spectrum of the current level. Therefore we introduce a refinement method in \Cref{sec:eigenpair_refinement} to overcome this difficulty.

\subsection*{Preconditioning In Eigenproblems} Before we proceed, we would like to have some discussions on the critical choice of the preconditioner $M= \Psi^T\Psi$ for inverting $A_{st}$ in the Lanzcos method. Though the preconditioner $M$ comes naturally from the derivation of our method, it reveals an important phenomenon that arises when we use CG type methods to handle matrix inversion in eigenproblems. 

Given a symmetric matrix $A$ with large condition number, we know that choosing a good preconditioner $C$ is critical for improving the performance of using the CG method to solve linear system $Ax=f$. Generally, such improvement is ``uniformly'' good for all right hand side $f$,  which may become a ``curse'' in eigenproblems. In an extreme case, suppose the right hand side $f$ is an eigenvector of $A$, then the CG method without any preconditioner actually converges exactly in one iteration. However, if $C$ does not preserve the eigenvectors of $A$, it will still take some ``uniform'' number of iterations to converge for the PCG method with preconditioner $C$. This happens, for example, when we choose $C=LL^T$ as the incomplete Cholesky decomposition of $A$, which is a common choice of preconditioner. 

Now consider computing the smallest eigenvalues of $A$ using the Lanzcos method. Suppose we have already computed some eigenspace $V$, then the next step would be computing $A^{-1}f$ for some $f\in V^{\perp}$. Therefore, the efficiency of using the CG method is subject to the restricted condition number $\kappa(A,V^{\perp})$. As $V$ gets larger, $\kappa(A,V^{\perp})$ gets smaller, and it takes less iterations for the CG method to converge (subject to some prescribed tolerance). However, using the PCG method with incomplete Cholesky preconditioning cannot benefit from what we have computed, since the preconditioner $C=LL^T$ compromises the spectral property that the right hand side $f$ is in a smaller and smaller invariant space of $A$. So it can be more efficient to use the CG method than to use the PCG method when we are computing a relative large number of partial eigenpairs of $A$ with the Lanzcos method. We will verify this phenomenon in numerical experiments in \Cref{sec:experiment}.

Inspired by this observation, we seek to combine the nice spectral property in the Lanzcos procedure and the advantage of preconditioning in the CG method. So in our method, we not only apply the multiresolution matrix decomposition to resolve the large condition number of $A$, but also use a proper choice of preconditoners with good spectral property so that we can take advantage from the narrowing down residual spectrum of $A$ in the Lanzcos procedure.

\section{Cross-level Refinement Of Eigenspace}
\label{sec:eigenpair_refinement}
In the previous section we have established a one level spectrum extension method, given that a partial accurate spectrum is provided. To develop this method into an inductive hierarchical spectrum completion procedure, a natural idea is to use the spectrum computed at the lower level as the initial spectrum to be used in the higher level. However, such initial spectrum is not actually good enough since there is a compression error between each two neighboring levels. Thus we need to use a compatible refinement technique to refine the initial spectrum.

Now consider the cross-level spectrum refinement between the two consecutive levels, the $h$-level and the $l$-level. The two operators are $\Theta^h=\bm{\Psi}^h\big((\bm{\Psi}^h)^TA\bm{\Psi}^h\big)^{-1}(\bm{\Psi}^h)^T$ and $\Theta^l=\bm{\Psi}^l\big((\bm{\Psi}^l)^TA\bm{\Psi}^l\big)^{-1}(\bm{\Psi}^l)^T$ respectively. We have the relations 
\[\bm{\Psi}^l=\bm{\Psi}^h\Psi^l,\quad \mathcal{U}^l=\bm{\Psi}^hU^l,\]
\[A_{st}^l=(\bm{\Psi}^l)^TA\bm{\Psi}^l=(\Psi^l)^T(\bm{\Psi}^h)^TA\bm{\Psi}^h\Psi^l=(\Psi^l)^TA_{st}^h\Psi^l,\]
\[B_{st}^l=(\mathcal{U}^l)^TA\mathcal{U}^l=(U^l)^T(\bm{\Psi}^h)^TA\bm{\Psi}^hU^l=(U^l)^TA_{st}^hU^l,\]
\[(A_{st}^h)^{-1}=\Psi^l(A_{st}^l)^{-1}(\Psi^l)^T+U^l(B_{st}^l)^{-1}(U^l)^T,\]
\begin{equation}
\Theta^h=\bm{\Psi}^h\big(\Psi^l(A_{st}^l)^{-1}(\Psi^l)^T+U^l(B_{st}^l)^{-1}(U^l)^T\big)(\bm{\Psi}^h)^T=\Theta^l+\mathcal{U}^l(B_{st}^l)^{-1}(\mathcal{U}^l)^T.
\label{eqt:Theta_decomposition}
\end{equation}

Now suppose that we have obtained the first $m_l$ essential eigenpairs $(\mu_{l,i},v_{l,i})$, $i=1,\cdots,m_l$, of $\Theta^l$. We want to use these eigenpairs as initial guess to obtain the first $m_h$ essential eigenpairs of $\Theta^h$. Recall that we have the estimates 
\[|\mu_{h,i}-\mu_{l,i}|\leq \varepsilon_l,\quad 1\leq i\leq m_l,\]
and 
\[\|\Theta^hv_{l,i}-\mu_{h,i}v_{l,i}\|_2\leq2\varepsilon_l,\quad 1\leq i\leq m_l,\]
where $\varepsilon_l$ is the compression error bound. These estimates give us confidence that we can obtain $(\mu_{h,i},v_{h,i})$, $i=1,\cdots,m_h$, efficiently from $(\mu_{l,i},v_{l,i})$, $i=1,\cdots,m_l$, by using some refinement technique.

Indeed, we will use the Orthogonal Iteration with Ritz Acceleration as our refinement method. Consider an initial guess $Q^{(0)}$ of the first $m$ eigenvectors of a SPD operator $\Theta$. To obtain more accurate eigenvalues and eigenspace, the Orthogonal Iteration with Ritz Acceleration runs as follows: 
\begin{align*}
&Q^{(0)}\in \mathbb{R}^{n\times m}\ \text{given with}\ (Q^{(0)})^TQ^{(0)}=I_m\\
&F^{(0)} = \Theta Q^{(0)}\\
&\textbf{for}\ k = 1,2,\cdots\\
&\qquad Q^{(k)}R^{(k)}=F^{(k-1)}\qquad (\text{QR factorization})\\ \tag{$\ast$} \label{eqt:orthogonaliteration}
&\qquad F^{(k)} = \Theta Q^{(k)}\\
&\qquad S^{(k)}=(Q^{(k)})^TF^{(k)}\\
&\qquad P^{(k)}D^{(k)}(P^{(k)})^T = S^{(k)}\qquad (\text{Schur decomposition})\\
&\qquad Q^{(k)} \leftarrow Q^{(k)}P^{(k)}\\
&\qquad F^{(k)} \leftarrow F^{(k)}P^{(k)}\\
&\textbf{end}
\end{align*}

To state the convergence property of the Orthogonal Iteration with Ritz Acceleration, we first define the distance between two spaces. Let $V_1,V_2\subset \mathbb{R}^n$ be two linear spaces, and $\bm{P}_{V_1},\bm{P}_{V_2}$ be the orthogonal projections onto $V_1,V_2$ respectively. We define the distance between $V_1$ and $V_2$ as
\[\mathrm{dist}(V_1,V_2)=\|\bm{P}_{V_1}-\bm{P}_{V_2}\|_2.\]
We also use the same notation $\mathrm{dist}(V_1,V_2)$ when $V_1,V_2$ are matrices of column vectors. In this case $\mathrm{dist}(V_1,V_2)$ means $\mathrm{dist}(\Span\{V_1\},\Span\{V_2\})$.

Suppose that the diagonal entries $\mu_i^{(k)}$, $i=1,\cdots,m$, of $D^{(k)}$ are in a decreasing order, then $\mu_i^{(k)}$ is a good approximation of the $i^{th}$ eigenvalue of $\Theta$, and $\Span\{Q^{(k)}_i\}$ is a good approximation of the eigenspace spanned by the first $i$ eigenvectors of $\Theta$, where $Q_i^{(k)}$ denotes the first $i$ columns of $Q^{(k)}$. We would like to emphasize that the meaning of the superscript $(k)$ of $\mu_i^{(k)}$ is different from those in \Cref{sec:spectrum_extension}. More precisely, we have the following convergence estimate:\\

\textbf{Theorem (Stewart, 1968):}\cite{stewart1969accelerating} Let $(\mu_i,v_i)$, $i=1,\cdots,N$, be the ordered (essential) eigenpairs of $\Theta$, and let $\mu_i^{(k)}$, $i=1,\cdots,m$, be the ordered eigenvalues of $D^{(k)}=(Q^{(k)})^T\Theta Q^{(k)}$ given in the Orthogonal Iteration with Ritz Acceleration \cref{eqt:orthogonaliteration}. Let $V_m=[v_1,v_2,\cdots,v_m]$, and $d^{(0)}=\mathrm{dist}(V_m,Q^{(0)})$. Then we have
\[|\mu_i-\mu_i^{(k)}|\leq O\Big(\big(\frac{\mu_{m+1}}{\mu_i}\big)^{2k}\cdot\|\Theta\|_2\cdot\frac{(d^{(0)})^2}{1-(d^{(0)})^2}\Big),\quad 1\leq i\leq m.\]
Moreover, we have 
\[\mathrm{dist}(V_m,Q^{(k)})\leq O\Big(\big(\frac{\mu_{m+1}}{\mu_m}\big)^{k}\cdot\frac{d^{(0)}}{\sqrt{1-(d^{(0)})^2}}\Big),\]
and for $i=1,\cdots,m-1,$ if we further assume that $\alpha_i=\mu_i-\mu_{i+1}>0$, then we have
\[\mathrm{dist}(V_i,Q^{(k)}_i)\leq O\Big(\big(\frac{\mu_{m+1}}{\mu_i}\big)^{k}\cdot\frac{d^{(0)}}{\sqrt{1-(d^{(0)})^2}}\Big)+O\Big(\frac{\sqrt{i}}{\alpha_i}\cdot\big(\frac{\mu_{m+1}^2}{\mu_m\mu_i}\big)^{k}\cdot \|\Theta\|_2\cdot\frac{(d^{(0)})^2}{1-(d^{(0)})^2}\Big),\] 
where $V_i$ and $Q_i^{(k)}$ are the first $i$ columns of $V_m$ and $Q^{(k)}$ respectively.\\

Now we go back to our problem, where we have $\Theta=\Theta^h$, $m=m_l$, and $Q^{(0)}=V_{m_l}^l=[v_{l,1},\cdots,v_{l,m_l}]$. We next consider the efficiency of this refinement technique in our problem. As long as the initial distance $d^{(0)}=\mathrm{dist}(V_{m_l}^h,V_{m_l}^l)<1$, the first $m_h$ eigenvalues and the eigenspace of the first $m_h$ eigenvectors of $\Theta^h$ converges exponentially fast at a rate $(\frac{\mu_{h,m_l+1}}{\mu_{h,m_h}})^k$. We can expect that a few iterations of refinement will be sufficient to give an accurate eigenspace for narrowing down the residual spectrum of $\Theta^h$, if we can ensure that the ratio $\frac{\mu_{h,m_l+1}}{\mu_{h,m_h}}$ is small enough. This will be verified in our numerical examples to be presented in \cref{sec:experiment}. In particular, to refine the first $m_h$ eigenpairs subject to a prescribed accuracy $\epsilon$, we need $K = O(\log(\frac{1}{\epsilon})/\log(\frac{\mu_{h,m_h}}{\mu_{h,m_l+1}}))$ refinement iterations. 

The main cost of the refinement procedure comes from the computation of $\Theta^hQ^{(0)}$ and the computation of $\Theta^hQ^{(k)}$ in each iteration. We will reduce the computational cost by using the fact that $Q^{(k)}$ is a good approximation of eigenvectors of $\Theta^h$. We first consider how to compute $\Theta^hQ^{(0)}$ efficiently.

Notice that in our problem, we take $Q^{(0)}=V^l_{m_l}$, whose columns are the first $m_l$ eigenvectors of $\Theta^l$. Therefore by \Cref{eqt:Theta_decomposition}, we have
\[\Theta^hQ^{(0)}=\Theta^hV^l_{m_l}=\Theta^lV^l_{m_l}+\mathcal{U}^l(B_{st}^l)^{-1}(\mathcal{U}^l)^TV^l_{m_l}=V^l_{m_l}D^l_{m_l}+\mathcal{U}^l(B_{st}^l)^{-1}(\mathcal{U}^l)^TV^l_{m_l},\]
where $D^l_{m_l}$ is a diagonal matrix whose diagonal entries are $\mu_{l,1},\mu_{l,2},\cdots,\mu_{l,m_l}$. Recall that by \Cref{lemma:Bst_eigenvalue} and \Cref{cor:Bconditonnumber}, $\kappa(B^{l}_{st})$ is bounded by $\varepsilon_l\delta_h$ that can be well controlled in the decomposition procedure. Thus it is efficient to solve $(B_{st}^l)^{-1}$ using the CG method. As we have mentioned before, applying $(U^l)^T$ or $U^l$ from the left is performed by doing patch-wise Householder transformations that involve only one local Householder vector on each patch, which takes $O(N^h)$ computational cost, where $N^h$ is the compressed dimension on level $h$ or the size of $A_{st}^h$. Therefore in the CG method, the cost of matrix multiplication of $B^l_{st}=(U^l)^TA_{st}^hU^l$ mainly comes from the number of nonzero entries of $A_{st}^h$. Then the total computational cost of computing $\Theta^hQ^{(0)}$ subject to a relative error $\epsilon$ can be bounded by 
\[O\left(m_l\cdot nnz(A^h_{st})\cdot \varepsilon_l \delta_h\cdot \log(\frac{1}{\epsilon})\right).\]

Next, we consider how to compute $\Theta^hQ^{(k)}$. To do so, we first compute $w^{(k)}_i=(\bm{\Psi}^h)^Tq^{(k)}_i$, where $q^{(k)}_i$ is the $i^{th}$ column of $Q^{(k)}$, then compute $(A_{st}^h)^{-1}w^{(k)}_i$, and apply $\bm{\Psi}^h$. Again we will use the PCG method with predictioner $M^h=(\bm{\Psi}^h)^T\bm{\Psi}^h$ to compute $(A_{st}^h)^{-1}w^{(k)}_i$. As we have discussed in \Cref{sec:spectrum_extension}, this is equivalent to using the CG method to compute $(A^h_{\bm{\Psi}})^{-1}z_i^{(k)}$, where $A^h_{\bm{\Psi}} = (M^h)^{-\frac{1}{2}}A_{st}^h(M^h)^{-\frac{1}{2}}$, and $z^{(k)}_i = (M^h)^{-\frac{1}{2}}w^{(k)}_i = (M^h)^{-\frac{1}{2}}(\bm{\Psi}^h)^Tq^{(k)}_i$. Inspired by \Cref{cor:PCG_convergence}, we seek to provide a good initial guess for the CG method to ensure efficiency. In the Orthogonal Iteration with Ritz Acceleration \cref{eqt:orthogonaliteration}, one can check that $(Q^{(k)})^T(\Theta^h Q^{(k)} - Q^{(k)}D^{(k)}) = \bm{0}$, where $D^{(k)}$ is a diagonal matrix with diagonal entries $\mu^{(k)}_1,\mu^{(k)}_2,\cdots,\mu^{(k)}_{m_l}$, and therefore 
\begin{align*}
&\ (Z^{(k)})^T\big((A^h_{\bm{\Psi}})^{-1}Z^{(k)}-Z^{(k)}D^{(k)}\big) \\
= &\ (Q^{(k)})^T \bm{\Psi}^h(M^h)^{-\frac{1}{2}}\left((A^h_{\bm{\Psi}})^{-1}(M^h)^{-\frac{1}{2}}(\bm{\Psi}^h)^TQ^{(k)} -(M^h)^{-\frac{1}{2}}(\bm{\Psi}^h)^TQ^{(k)}D^{(k)}\right)\\
=&\ (Q^{(k)})^T \left(\bm{\Psi}^h(A_{st}^h)^{-1}(\bm{\Psi}^h)^TQ^{(k)} - \bm{\Psi}^h(M^h)^{-1}(\bm{\Psi}^h)^TQ^{(k)}D^{(k)}\right) \\
=&\ (Q^{(k)})^T \left(\Theta^hQ^{(k)} - Q^{(k)}D^{(k)}\right)\\
=&\ \bm{0},
\end{align*}
where we have used that $Q^{(k)}\in \Span \{\bm{\Psi}^h\}$ and so $\bm{\Psi}^h(M^h)^{-1}(\bm{\Psi}^h)^TQ^{(k)}=Q^{(k)}$. This observation implies that if we use $\mu^{(k)}_iz_i^{(k)}$ as the initial guess for computing $(A^h_{\bm{\Psi}})^{-1}z_i^{(k)}$ using the CG method, the initial residual $z_i^{(k)}-(A^h_{\bm{\Psi}})(\mu^{(k)}_iz_i^{(k)})$ is orthogonal to $(A^h_{\bm{\Psi}})^{-1}Z^{(k)}$. Since $Q^{(k)}$ are already good approximate essential eigenvectors of $\Theta^h$, $Z^{(k)}$ are good approximate eigenvectors of $(A_{\bm{\Psi}}^h)^{-1}$, we can expect that the target eigenspace $Z_{m_h}$, namely the eigenspace of the first $m_h$ eigenvectors of $(A_{\bm{\Psi}}^h)^{-1}$, can be well spanned in $\Span\{(A^h_{\bm{\Psi}})^{-1}Z^{(k)}\}$. Therefore we can reasonably assume that $z_i^{(k)}-(A^h_{\bm{\Psi}})(\mu^{(k)}_iz_i^{(k)})\in Z_{m_h^+}=Z_{m_h}^\perp$, and so again we can benefit from the restricted condition number $\kappa(A^h_{\bm{\Psi}},Z_{m_h^+})\leq \mu_{h,m_h+1}\delta_h$ as introduced in \Cref{sec:spectrum_extension}. Moreover, we notice that the spectral residual $\|\Theta^hq_i^{(k)}-\mu_i^{(k)}q_i^{(k)}\|_2$ is bounded by $2\varepsilon_l$ by \Cref{lemma:approximate_eigenpair_multiresolution}, and we have 
\begin{equation}
\|(A^{h}_{st})^{-1}w^{(k)}_i-\mu_i^{(k)}(M^h)^{-1}w^{(k)}_i\|_2\leq \|(A^h_{\bm{\Psi}})^{-1}z_i^{(k)}-\mu^{(k)}_iz_i^{(k)}\|_2 = \|\Theta^hq_i^{(k)}-\mu_i^{(k)}q_i^{(k)}\|_2,
\label{eqt:estimate_in_refine}
\end{equation}
where we have used $\lambda_{min}(M^h)\geq 1$ (\Cref{lemma:multi_psipsi_conditionnumb}). Thus if we use $\mu^{(k)}_iz_i^{(k)}$ as the initial guess, the initial error will be bounded by $2\varepsilon_l$ at most, and the CG procedure will only need 
\[O\left(\kappa(A^h_{\bm{\Psi}},Z_{m_h^+})\cdot \log(\frac{\varepsilon_l}{\epsilon})\right) = O\left(\mu_{h,m_h+1}\delta_h\cdot \log(\frac{\varepsilon_l}{\epsilon})\right)\]
iterations to achieve a relative accuracy $\epsilon$, instead of $O(\kappa(A^h_{\bm{\Psi}},Z_{m_h^+})\cdot \log(\frac{1}{\epsilon}))$. Notice that using the initial guess $\mu^{(k)}_iz_i^{(k)}$ for $(A^h_{\bm{\Psi}})^{-1}z_i^{(k)}$ is equivalent to using the initial guess $\mu_i^{(k)}(M^h)^{-1}w^{(k)}_i$ for $(A^{h}_{st})^{-1}w^{(k)}_i$. 

Supported by the analysis above, we will compute $(A_{st}^h)^{-1}w^{(k)}_i$ using the preconditioned CG method with preconditioner $M^h$ and initial guess $\mu^{(k)}_i(M^h)^{-1}w_i^{(k)}$. Again suppose that in each PCG iteration, we also use the CG method to apply $(M^h)^{-1}$ subject to a higher relative accuracy $\hat\epsilon$, which takes $O(nnz(M^h)\cdot \kappa(M^h)\cdot \log(\frac{1}{\hat\epsilon}))$ computational cost. In practice, it is sufficient to take $\hat \epsilon$ comparable to $\epsilon$. Recall that $nnz(M^h)\leq nnz(A^h_{st})$, and $\kappa(M^h)\leq O(\varepsilon_h\delta_h)$ (\Cref{lemma:multi_psipsi_conditionnumb}), the cost of computing $\Theta^hQ^{(k)}$ subject to a relative error $\epsilon$ is then bounded by  
\[O\left(m_l\cdot \mu_{h,m_h+1}\delta_h\cdot \log(\frac{\varepsilon_l}{\epsilon})\cdot nnz(A^h_{st})\cdot \varepsilon_h\delta_h\cdot \log(\frac{1}{\epsilon})\right).\] 

Notice that in each refinement iteration we also need to perform one QR factorization and one Schur decomposition, which together cost $O(N^h\cdot m_l^2)$. However, as we have mentioned in the introduction, we only consider the asymptotic complexity of our method when the original $A$ becomes super large. In this case, the number $m_{tar}$ of the target eigenpairs is considered as a fixed constant, and so the term $O(N^h\cdot m_l^2)\leq O(N^h m_{tar}^2)$ is considered to be minor and will be omitted in our complexity analysis. Therefore, the total cost of refining the first $m_h$ eigenpairs subject to a prescribed accuracy $\epsilon$ can be bounded by 
\begin{align}
& O\left(m_l\cdot nnz(A^h_{st})\cdot \varepsilon_l \delta_h\cdot \log(\frac{1}{\epsilon})\right) \\ &\quad +O\left(m_l\cdot \mu_{h,m_h+1}\delta_h\cdot \log(\frac{\varepsilon_l}{\epsilon})\cdot nnz(A^h_{st})\cdot \varepsilon_h\delta_h\cdot \log(\frac{1}{\epsilon})\cdot \log(\frac{1}{\epsilon})/\log(\frac{\mu_{h,m_h}}{\mu_{h,m_l+1}})\right).\nonumber
\label{eqt:refine_complexity}
\end{align}

Again we remark that the operator $\Theta^h$, the long vectors $Q^{(k)}$, $F^{(k)}$, $V^l$ and $V^h$ are only for analysis use. Operations on long vectors of size $n$ will be very expensive and unnecessary, especially on lower levels where the compression dimension $N^{h}$(the size of $A_{st}^h$) is small. Notice that all long vectors on the $h$-level are in $\Span\{\bm{\Psi}^h\}$ as
\[Q^{(k)}=\bm{\Psi}^h\widehat{Q}^{(k)},\quad F^{(k)}=\bm{\Psi}^h\widehat{F}^{(k)},\quad V^l_{m_l}=\bm{\Psi}^h\widehat{V}^l_{m_l},\quad V^h_{m_h}=\bm{\Psi}^h\widehat{V}^h_{m_h},\]
we thus only operate on their coefficients in the basis $\bm{\Psi}^h$. Correspondingly, whenever we need to consider orthogonality of long vectors, we replace it by the $M^h$-orthogonality of their coefficient vectors. One can check that all discussions above still apply. Also another advantage of using the coefficient vectors is that in the previous discussions, the good initial guess $\mu_i^{(k)}(M^h)^{-1}w^{(k)}_i=\mu_i^{(k)}(M^h)^{-1}(\bm{\Psi}^h)^Tq^{(k)}_i=\mu_i^{(k)}\hat{q}^{(k)}$ is obtained explicitly.\\

Summarizing the analysis above, we propose the following \Cref{alg:eigenpair_refinement} as our refinement method. Since we want the eigenspace spanned by the first $m_h$ eigenvectors of $\Theta^h$ to be computed accurately, the refinement stops when $\mathrm{dist}(Q^{(k-1)}_{m_h},Q^{(k)}_{m_h})<\epsilon$ for some prescribed accuracy $\epsilon$, where $Q^{(k)}_{m_h}$ denotes the first $m_h$ columns of $Q^{(k)}$. Since $Q^{(k)}$ is orthogonal, one can check that 
\begin{align*}
\mathrm{dist}(Q^{(k-1)}_{m_h},Q^{(k)}_{m_h})=&\ \|Q^{(k)}_{m_h}-Q^{(k-1)}_{m_h}(Q^{(k-1)}_{m_h})^TQ^{(k)}_{m_h}\|_2\\
=&\ \|\widehat{Q}^{(k)}_{m_h}-\widehat{Q}^{(k-1)}_{m_h}(\widehat{Q}^{(k-1)}_{m_h})^TM^h\widehat{Q}^{(k)}_{m_h}\|_{M^h}\\
\leq&\ \sqrt{\lambda_{max}(M^h)}\|\widehat{Q}^{(k)}_{m_h}-\widehat{Q}^{(k-1)}_{m_h}(\widehat{Q}^{(k-1)}_{m_h})^TM^h\widehat{Q}^{(k)}_{m_h}\|_2\\
\leq&\ \sqrt{1+\varepsilon_h\delta_h}\|\widehat{Q}^{(k)}_{m_h}-\widehat{Q}^{(k-1)}_{m_h}(\widehat{Q}^{(k-1)}_{m_h})^TM^h\widehat{Q}^{(k)}_{m_h}\|_F.
\end{align*}
In practical, we use $\|\widehat{Q}^{(k)}_{m_h}-\widehat{Q}^{(k-1)}_{m_h}(\widehat{Q}^{(k-1)}_{m_h})^TM^h\widehat{Q}^{(k)}_{m_h}\|_F<\frac{\epsilon}{\sqrt{1+\varepsilon_h\delta_h}}$ as the stopping criterion since it is easy to check. We have used \Cref{lemma:multi_psipsi_conditionnumb} to bound $\lambda_{max}(M^h)$. 

\begin{algorithm}[h!]
\caption{Eigenpair Refinement}
\label{alg:eigenpair_refinement}
\begin{algorithmic}[1]
\REQUIRE{$\widehat{V}_{m_l}^l$, $D_{m_l}^l$, prescribed accuracy $\epsilon$, target eigenvalue threshold $\mu_h$.}
\ENSURE{$\widehat{V}_{m_h}^h$, $D_{m_h}^h$.}
\STATE{Set $\widehat{Q}^{(0)}=V_{m_l}^l$,\quad $D^{(0)}=D_{m_l}^l$,\quad $k=0$;}
\FOR{$i=1:m_l$}
\STATE{$g_i=pcg(B_{st}^l,(U^l)^TM^h\hat{q}_i^{(0)},-,-,\epsilon)$;\quad ($\widehat{Q}=[\hat{q}_1,\cdots,\hat{q}_{m_l}]$)}
\ENDFOR
\STATE{$\widehat{F}^{(0)}=\widehat{Q}^{(0)}D^{(0)}+U^lG$;\quad ($G=[g_1,\cdots,g_{m_l}]$)}
\REPEAT
\STATE{$k \leftarrow k + 1$;}
\STATE{$\widehat{Q}^{(k)}R^{(k)}=\widehat{F}^{(k-1)}$;\quad (QR factorization with respect to $M^h$ orthogonality, i.e. $(\widehat{Q}^{(k)})^TM^h\widehat{Q}^{(k)}=I$)}
\STATE{$W^{(k)}=M^h\widehat{Q}^{(k)}$;}
\FOR{$i = 1:m_l$}
% \STATE{$\bar{w}_i^{(k)}=pcg(M^h,w_i^{(k)},-,-,0.1*\epsilon)$;}\label{line:initialguess}
\STATE{$\hat{f}_i^{(k)}=pcg(A^h_{st},w_i^{(k)},M^h,\mu_i^{(k-1)}\hat{q}_i^{(k)},\epsilon)$;\quad ($\widehat{F}=[\hat{f}_1,\cdots,\hat{f}_{m_l}]$)}
\ENDFOR
\STATE{$S^{(k)}=(W^{(k)})^T\widehat{F}^{(k)}$;}
\STATE{$P^{(k)}D^{(k)}(P^{k)})^T=S^{(k)}$\quad (Schur decomposition, diagonals of $D^{(k)}$ in decreasing order);}
\STATE{renew $m_h$ so that $\mu_{m_h}^{(k)}\geq \mu_h>\mu_{m_h+1}^{(k)}$;}
\STATE{$\widehat{Q}^{(k)}\leftarrow\widehat{Q}^{(k)}P^{(k)}$,\quad $\widehat{F}^{(k)}\leftarrow\widehat{F}^{(k)}P^{(k)}$;}
\UNTIL{$\|\widehat{Q}^{(k)}_{m_h}-\widehat{Q}^{(k-1)}_{m_h}(\widehat{Q}^{(k-1)}_{m_h})^TM^h\widehat{Q}^{(k)}_{m_h}\|_F<\epsilon$.}
\STATE{$\widehat{V}_{m_h}^h=\widehat{Q}^{(k)}_{m_h}$,\quad $D_{m_h}^h=D^{(k)}_{m_h}$. \quad ($D^{(k)}_{m_h}$ denotes the first $m_h$-size block of $D^{(k)}$)}
\end{algorithmic}
\end{algorithm}

\section{Overall Algorithms}
\label{sec:overall_algorithm}
Combining the refinement method and the extension method, we now propose our overall \Cref{alg:hierarch_eigencompute} for computing partial eigenpairs of a SPD matrix $A$. It utilizes the {\it a priori} multiresolution decomposition of $A$ to compute the first $m_{tar}$ eigenpairs of $A^{-1}$, by passing approximate eigenpairs from lower levels to higher levels to finally reach a prescribed accuracy. In particular, this algorithm starts with the eigen decomposition of the lowest level (whose dimension is small enough), refines and extends the approximate eigenpairs on each level, and stops at the highest level. The overall accuracy is achieved by the prescribed compression error of the highest level.

Recall that the output $\widehat{V}^{(k)}_{ex}$ of the extension process and the initializing process are the coefficients of $V^{(k)}_{ex}$ in the basis $\bm{\Psi}^{(k)}$. When passing these results from level $k$ to level $k-1$, we need to recover the coefficients of $V^{(k)}_{ex}$ in the basis $\bm{\Psi}^{(k-1)}$. This can be done by simply reforming $\widehat{V}^{(k)}_{ex}\leftarrow\Psi^{(k)}\widehat{V}^{(k)}_{ex}$(Line \ref{line:reform} in \Cref{alg:hierarch_eigencompute}), since $V^{(k)}_{ex}=\bm{\Psi}^{(k)}\widehat{V}^{(k)}_{ex}=\bm{\Psi}^{(k-1)}\Psi^{(k)}\widehat{V}^{(k)}_{ex}$.

In \Cref{alg:hierarch_eigencompute}, the parameters should be chosen carefully to ensure computational efficiency, by using the analysis in the previous sections. We shall discuss the choice of each parameter separately. To be consistent, we first clarify some notations. Let $\hat{m}_k,m_k$ be the numbers of output eigenpairs of the refinement process and the extension process respectively on level $k$. Ignoring numerical errors, let $(\mu_i^{(k)},v_i^{(k)})$, $i=1,\cdots,N^{(k)}$, be the essential eigenpairs of the operator $\Theta^{(k)}$ as in \Cref{sec:spectrum_extension}. Let $(\mu_i^{(k)},v_i^{(k)})$, $i=1,\cdots,m_k$, denote the output eigenpairs on level $k$. Notice that $(\mu_i^{(k)},v_i^{(k)})$, $i=1,\cdots,\hat{m}_k$, are the output of the refinement process, and $(\mu_i^{(k)},v_i^{(k)})$, $i=\hat{m}_k+1,\cdots,m_k$, are the output of the extension process. We will use $(\tilde{\mu},\tilde{v})$ to denote the numerical output of $(\mu,v)$.

\textbf{Choice Of Multi-level Accuracies $\{\epsilon^{(k)}\}$:} Notice that there is a compression error $\varepsilon_k$ between level $k$ and level $k-1$. That is to say, no matter how accurately we compute the eigenpairs of $\Theta^{(k)}$, they are approximations of eigenpairs of $\Theta^{(k-1)}$ subject to accuracy no better that $\varepsilon_k$. Therefore, on the one hand, the choice of the algorithm accuracy $\epsilon^{(k)}$ for the eigenpairs of $\Theta^{(k)}$ on each level should not compromise the compression error. On the other hand, the accuracy should not be over-achieved due to the presence of the compression error. Therefore, we choose $\epsilon^{(k)}=0.1\times\varepsilon_k$ in practice.

\textbf{Choice Of Thresholds $\{(\mu^{(k)}_{re},\mu^{(k)}_{ex})\}_{k=1}^{K}$:} These thresholds provide control on the smallest eigenvalues of output eigenpairs of both the refinement process and the extension process in that 
\[\mu_{\hat{m}_k}^{(k)}\geq \mu_{re}^{(k)}>\mu_{\hat{m}_k+1}^{(k)},\quad \mu_{ex}^{(k)}\geq \mu_{m_k}^{(k)},\quad k=1,2,\cdots,K.\]
Recall that the outputs of the refinement process are the inputs of the extension process, and the outputs of the extension process are the inputs of the refinement process on the higher level. By \Cref{thm:multi_spectrum_completion_overall}, to ensure the efficiency of the extension process, we need to uniformly control the restricted condition number 
\begin{equation*}
\kappa(A^{(k)}_\Psi,Z_{\hat{m}^+_k}^{(k)})\leq \mu^{(k)}_{\hat{m}_k+1}\delta_k<\mu^{(k)}_{re}\delta_k.
\end{equation*}
Recall that in \Cref{sec:eigenpair_refinement} the convergence rate of the refinement process is given by $\frac{\mu_{h,m_l+1}}{\mu_{h,m_h}}$, where $l$ corresponds to $k+1$ and $h$ corresponds to $k$ on each level $k$. Thus to ensure the efficiency of the refinement process we need to uniformly control the ratio
\begin{equation*}
\frac{\mu_{m_{k+1}+1}^{(k)}}{\mu_{\hat{m}_{k}}^{(k)}}\leq \frac{\mu_{m_{k+1}}^{(k)}}{\mu^{(k)}_{re}}\leq \frac{\mu_{m_{k+1}}^{(k+1)}+\varepsilon_{k+1}}{\mu^{(k)}_{re}}\leq \frac{\mu_{ex}^{(k+1)}+\varepsilon_{k+1}}{\mu^{(k)}_{re}},
\end{equation*}
where $\varepsilon_{k+1}$ is the compression error between level $k+1$ and level $k$, and we have used \Cref{lemma:approximate_eigenpair_multiresolution}. Thus, more precisely, we need to choose thresholds $\{(\mu^{(k)}_{re},\mu^{(k)}_{ex})\}_{k=1}^{K}$ so that there exist uniform constants $\kappa>0,\gamma\in(0,1)$ so that 
\begin{equation}
\text{(i)}\ \mu^{(k)}_{re}\delta_k\leq \kappa,\qquad \text{(ii)}\ \frac{\mu_{ex}^{(k+1)}+\varepsilon_{k+1}}{\mu^{(k)}_{re}}\leq \gamma.
\label{efficiency_condition}
\end{equation}
Due to the existence of $\varepsilon_k$, condition (ii) implies that there is no need to choose $\mu_{ex}^{(k)}$ much smaller than $\varepsilon_k$, which suffers from over-computing but barely improves the efficiency of the refinement process. So one convenient way is to choose 
\begin{equation}
\mu_{re}^{(k)}=\alpha\varepsilon_{k+1},\qquad \mu_{ex}^{(k)}=\beta\varepsilon_k,
\label{eqt:parameter_choice}
\end{equation}
for some uniform constants $\alpha,\beta>0$ such that $\alpha>1+\beta$. Recall that when constructing the multiresolution decomposition, we impose conditions $\varepsilon_k\delta_k\leq c$ and $\varepsilon_{k}=\eta\varepsilon_{k+1}$ for some uniform constants $c>0$ and $\eta\in (0,1)$. Thus we have  
\[\mu^{(k)}_{re}\delta_k=\frac{\alpha}{\eta}\varepsilon_k\delta_k\leq\frac{\alpha c}{\eta}=\kappa,\qquad \frac{\mu_{ex}^{(k+1)}+\varepsilon_{k+1}}{\mu^{(k)}_{re}}=\frac{1+\beta}{\alpha}=\gamma<1.  \]

\textbf{Choice of Searching Step $d$:} In the first part of the extension algorithm, we explore the number $m_k$ so that $\mu_{m_k}^{(k)}\leq \mu_{ex}^{(k)}$, and we do this by setting an exploring step size $d$ and examining the last few eigenvalues every $d$ steps of the Lanczos iteration. The step size $d$ should neither be too large to avoid over computing, nor too small to ensure efficiency. In practical, we choose $d=\min\{\lfloor\frac{\dim\bm{\Psi}^{(k)}}{10}\rfloor,\lfloor\frac{m_{tar}}{10}\rfloor\}$.

\textbf{Complexity:} Now we summarize the complexity of \Cref{alg:hierarch_eigencompute} for computing the first $m_{tar}$ largest eigenpairs of $A^{-1}$ for a SPD matrix $A\in \mathbb{R}^{n\times n}$ subject to an error $\varepsilon$. Suppose we are provided a $K$-level multiresolution matrix decomposition of $A$ with $\varepsilon_k\delta_k\leq c$, $\varepsilon_k=\eta \varepsilon_{k+1}$, and $\varepsilon_1=\varepsilon$. In what follows, we will uniformly estimate $nnz(A_{st}^{(k)})\leq nnz(A)$, $\epsilon^{(k)}\geq \epsilon^{(1)} = 0.1 \varepsilon_1$ and $m_k\leq m_{tar}$.

We first consider the complexity of all refinement process. Notice that by our choice $\frac{\varepsilon_{k+1}}{\epsilon^{(k)}}=\frac{\varepsilon_{k+1}}{0.1\epsilon^{(k)}}=\frac{1}{0.1\eta}$, the factor $\log(\frac{\varepsilon_l}{\epsilon})$ in \Cref{eqt:refine_complexity}, which is now $\log(\frac{\varepsilon_{k+1}}{\epsilon^{(k)}})$, can be estimated as $O(\log(\frac{1}{\eta}))$. Since we can will make sure $\frac{\mu_{m_{k+1}+1}^{(k)}}{\mu_{\hat{m}_{k}}^{(k)}}\leq \gamma$ for some constant $\gamma<1$, the factor $\log(\frac{\mu_{h,m_h}}{\mu_{h,m_l+1}})$ in \Cref{eqt:refine_complexity}, which is now $\log(\frac{\mu_{\hat{m}_{k}}^{(k)}}{\mu_{m_{k+1}+1}^{(k)}})$, can be seen as a constant. Also using estimates $\mu_{h,m_h}\delta_h\leq \frac{\alpha c}{\eta}=O(\frac{c}{\eta})$, $\varepsilon_l\delta_h\leq \frac{c}{\eta}$, $\varepsilon_h\delta_h\leq c$ and $\log\frac{1}{\epsilon}=O(\log\frac{1}{\varepsilon})$, we modify \Cref{eqt:refine_complexity} to obtain the complexity of all $K$-level refinement process 
\begin{equation}
O\left(m_{tar}\cdot nnz(A)\cdot \frac{c^2}{\eta}\log(\frac{1}{\eta})\cdot (\log\frac{1}{\varepsilon})^2\cdot K \right).
\label{eqt:refine_complexity_total}
\end{equation}

Next we consider the complexity of all extension process. As we have discussed in \Cref{sec:spectrum_extension}, the major cost of the extension process comes from the operation of adding a new vector (the adding operation) to the Lanzcos vectors (Line 7 of \cref{alg:arnoldi_iteration_M} that happens in line 3 of \Cref{alg:eigenpair_extension}). Using estimates $\mu_m\delta(\mathcal{P})\leq \frac{\alpha c}{\eta} = O(\frac{c}{\eta})$, $\varepsilon(\mathcal{P})\delta(\mathcal{P})\leq c$, $\log\frac{1}{\epsilon}=O(\log\frac{1}{\varepsilon})$, we modify \cref{eqt:extend_complexity} to obtain the cost of every single call of the adding operation as
\[O\left(\frac{c^2}{\eta}\cdot nnz(A)\cdot (\log\frac{1}{\varepsilon})^2\right).\]  
On every level, the indexes contributing to adding operations go from $\hat m_k +1$ to $m_k$. Due to the refinement process, we have $\hat m_k\leq m_{k+1}$, and so every single index from $1$ to $m_{tar}$ may contribute more than one adding operations. But if we reasonably assume that $\mu_{ex}^{(k+1)}>\mu_{re}^{(k-1)}$, namely $\beta>\alpha\eta$ under parameter choice \Cref{eqt:parameter_choice}, we will have $m^{(k+1)}<\hat m^{(k-1)}$, and so every index from $1$ to $m_{tar}$ will contribute no more than two adding operations. Therefore the total cost of all extension process can be estimated as 
\begin{equation}
O\left(m_{tar}\cdot\frac{c^2}{\eta}\cdot nnz(A)\cdot (\log\frac{1}{\varepsilon})^2\right).
\label{eqt:extend_complexity_total}
\end{equation}
We remark that the cost of implicit restarting process is only a constant multiple of \Cref{eqt:extend_complexity_total}.
Combining \Cref{eqt:refine_complexity_total} and \Cref{eqt:extend_complexity_total}, we obtain the total complexity of our method
\begin{equation}
O\left(m_{tar}\cdot nnz(A)\cdot \frac{c^2}{\eta}\log(\frac{1}{\eta})\cdot (\log\frac{1}{\varepsilon})^2\cdot K \right).
\label{eqt:complexity_total}
\end{equation}
To further simplify \Cref{eqt:complexity_total}, we need to use estimates for the multiresolution matrix decomposition given in the previous work \cite{hou2017adaptive}. In particular, to preserve sparsity $nnz(A^{(k)}_{st})\leq nnz(A)$, we need to choose the scale ratio $\eta^{-1} = (\log\frac{1}{\varepsilon}+\log n)^p$ for some constant $p$. We remark that for graph Laplacian, $p=1$. The resulting level number is $K=O(\frac{\log n}{\log(\log\frac{1}{\varepsilon}+\log n)})$. The condition bound $c$ can be imposed to be uniform constant by the algorithm given in \cite{hou2017adaptive}. Then the overall complexity of \Cref{alg:hierarch_eigencompute} can be estimated as 
\begin{equation}
O\left(m_{tar}\cdot nnz(A)\cdot (\log\frac{1}{\varepsilon} + \log n)^p\cdot (\log\frac{1}{\varepsilon})^2 \cdot \log n \right) = O\left(m_{tar}\cdot nnz(A)\cdot (\log\frac{1}{\varepsilon} + \log n)^{p+3}\right).
\label{eqt:complexity_total_2}
\end{equation}

\begin{algorithm}[h!]
\caption{Hierarchical Eigenpair Computation}
\label{alg:hierarch_eigencompute}
\begin{algorithmic}[1]
\REQUIRE{$K$-level decomposition $\{\Theta^{(k)}\}_{k=1}^{K}$ of SPD matrix $A$, target number $m_{tar}$, searching step $d$, \\
prescribed multi-level accuracies $\{\epsilon^{(k)}\}$, extension thresholds $\{\mu^{(k)}_{ex}\}_{k=1}^{K}$, refinement thresholds $\{\mu^{(k)}_{re}\}_{k=1}^{K}$.}
\ENSURE{$V$, $D$.}
%\STATE{Find the eigen decomposition $[V^{(K)}_{ex},D^{(K)}_{ex}]$ of $\Theta^{(K)}$;}
\STATE{Find the eigen pairs $[\widehat{V}^{(K)}_{ex},D^{(K)}_{ex}]$ of the eigen problem $(A_{st}^{(K)})^{-1}M^{(K)}x=\mu x$;}
\FOR{$k = K-1:1$}
\STATE{$\widehat{V}^{(k+1)}_{ex}\leftarrow \Psi^{(k+1)}\widehat{V}^{(k+1)}_{ex}$} \label{line:reform}
\STATE{$[\widehat{V}^{(k)}_{ini},D^{(k)}_{ini}]=\text{Eigen\_Refine}([\widehat{V}^{(k+1)}_{ex},D^{(k+1)}_{ex}];\epsilon^{(k)},\mu_{re}^{(k)})$;}
\STATE{$op=OP(\ \cdot\ ;A^{(k)},M^{(k)},\epsilon^{(k)})$;}
\STATE{$[\widehat{V}^{(k)}_{ex},D^{(k)}_{ex}]=\text{Eigen\_Extend}([\widehat{V}^{(k)}_{ini},D^{(k)}_{ini}];op,\epsilon^{(k)},\mu_{ex}^{(k)},d,m_{tar})$;}
\ENDFOR
\STATE{$V=\Psi^{(1)}\widehat{V}^{(1)}_{ex}$\quad $D=D^{(1)}_{ex}$.}
\end{algorithmic}
\end{algorithm}

\section{Numerical Examples}
\label{sec:experiment}
In this section we present several numerical examples for the eigensolver. We will use \Cref{alg:hierarch_eigencompute} to compute a relative large number of eigenpairs of large matrices subject to prescribed accuracies.

\subsection{Dataset Description}
\begin{figure}[ht]
\centering
\includegraphics[width=.3\textwidth]{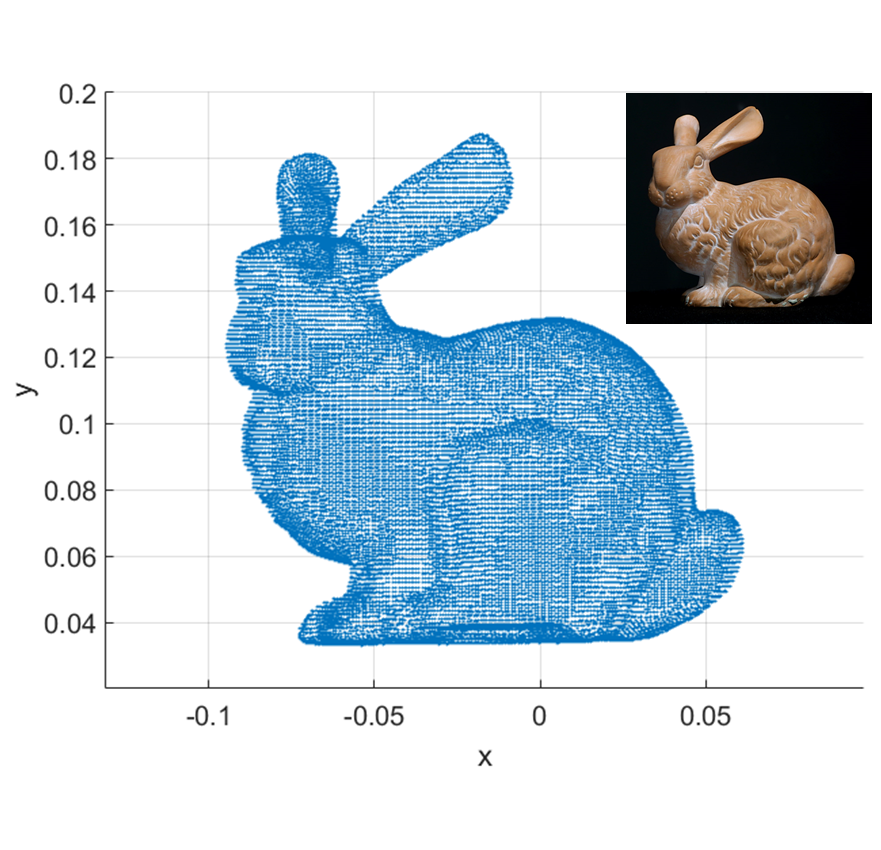} 
\includegraphics[width=.3\textwidth]{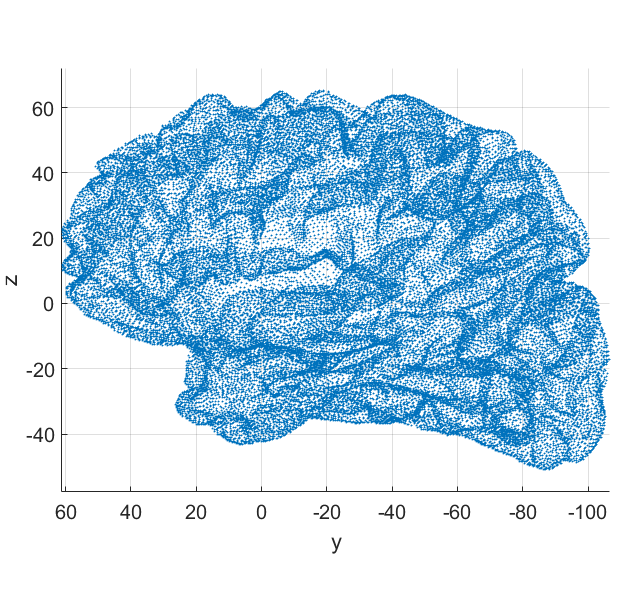} 
\includegraphics[width=.3\textwidth]{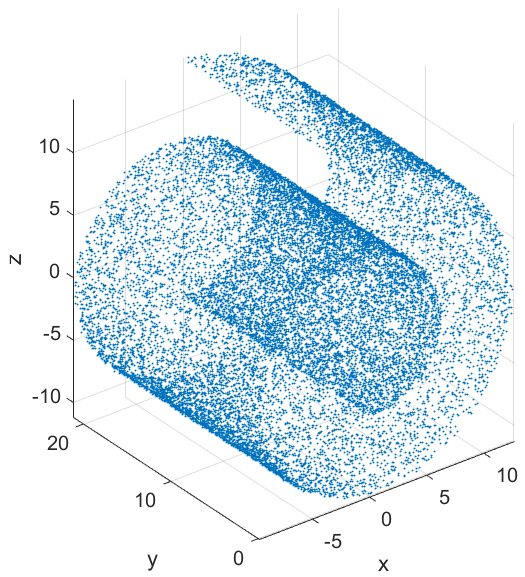} 
\caption{Datasets. From left to right: (1) bunny (point cloud and sculpture); (2) brain; (3) swiss roll.}
\label{fig:datasets}
\end{figure}

The datasets we use are drawn from different physical contexts. They are generated as 3D point clouds and  transformed into graphs by adding edges in the K-Nearest Neighbors (KNN) setting.
\begin{itemize}
\item The first dataset is the well-known ``Stanford Bunny'' from Stanford 3D Scanning Repository\footnote{http://graphics.stanford.edu/data/3Dscanrep/}. A reconstructed bunny has 35947 vertices that can be embedded into a surface in $\mathbb{R}^3$ with 5 holes in the bottom. 
\item The second dataset is a MRI data of brain from the Open Access Series of Imaging Sciences (OASIS)\footnote{http://www.oasis-brains.org/}. They use FreeSurfer to reconstruct the surface from MRI scan and obtain a point cloud with 48463 points. 
\item The third dataset is a ``SwissRoll'' model, which is popular in manifold learning. Vertices are generated by
\begin{equation}
(x_i, y_i, z_i) = (t_i cos(t_i), y_i, t_i sin(t_i)) + \boldsymbol{\eta}_i, \, i=1,2,...,n,
\end{equation}
where $t_i \stackrel{\text{i.i.d}}{\sim} \mathcal{U}[1.5\pi, 4.5\pi]$, $y_i \stackrel{\text{i.i.d}}{\sim} \mathcal{U}[0, 20]$, and $\eta_{i} \stackrel{\text{i.i.d}}{\sim} \mathcal{N}(\bm{0}, 0.05I_3)$. It can be viewed as a spiral of one and a half rounds plus random noise. In our examples the roll has $n=20000$ points.\\
\end{itemize}
With point clouds at hand, we apply the k-nearest neighbour (kNN) to construct graphs with $k_{bunny}=20$, $k_{brain}=20$ and $k_{swissroll}=10$. Each existing edge $e_{ij}$ is weighted as $e^{-r_{i,j}^2/\sigma}$, where $r_{i,j}$ is the Euclidean distance between vertices $v_i$ and $v_j$, and $\sigma$ is a parameter. We have $\sigma_{bunny} = 10^{-6}$, $\sigma_{brain} = 10^{-4}$ and $\sigma_{swiss} = 0.1$. Figure \ref{fig:datasets} shows the point clouds of datasets.\\\\
From the graphs given above, we construct their related graph laplacians $L$ in the general setting:
\[L_{ij}=\left\{\begin{array}{cc} 
 \sum_{k\sim i} w_{ik}, & i=j, \\\\
 -w_{ij}, & i\neq j.
 \end{array}\right.\]
Further, without loss of generality, we rescale all graph laplacians and add uniform selfloops of weight 1 to them, so that each of them satisfies (i)$\lambda_1=1$, (ii) $\lambda_2=O(1)$. Under this construction, we obtain three graph laplacian matrices $L_{bunny},L_{brain},L_{swissroll}$. $L_{bunny}$ has size $n=35947$, sparsity $nnz=714647$ and condition number $\kappa(L_{bunny})=1.86\times 10^4$; $L_{brain}$ has size $n=48463$, sparsity $nnz=1038065$ and condition number $\kappa(L_{bunny})=1.14\times 10^5$; $L_{swissroll}$ has size $n=20000$, sparsity $nnz=248010$ and condition number $\kappa(L_{bunny})=1.15\times 10^6$.

\subsection{Numerical Multiresolution Matrix Decomposition}

Before computing eigenpairs of graph laplacians from our datasets using \Cref{alg:hierarch_eigencompute}, we need to apply Algorithm 6 proposed in \cite{hou2017adaptive} to obtain the multiresolution decompositions. For each graph laplacian, we perform the decomposition with a prescribed condition bound $c$ and a series of multi-level resolutions (compression errors) $\{\varepsilon_k\}_{k=1}^{K}$. Note that we perform two decompositions with different multi-resolutions for the SwissRoll data.

\Cref{table:Decomposition} and \Cref{table:SurfaceDecomposition} give the detailed information of all decompositions we will use for eigenpair computation. In \Cref{table:Decomposition}, $K$ is the number of levels, $\varepsilon_1$ is the finest (prescribed) accuracy, $\eta$ is the ratio $\varepsilon_k/\varepsilon_{k+1}$ and $c$ is the condition bound such that $\varepsilon_k\delta_k\leq c$. By \Cref{lemma:multi_psipsi_conditionnumb}, the condition number of $M^{(k)}$ is bounded as $\kappa(M^{(k)})\leq 1+\varepsilon_k\delta_k\approx c$, and by \Cref{cor:Bconditonnumber}, the condition number of $B^{(k)}$ is bounded as $\kappa(B^{(k)})\leq \varepsilon_k\delta_{k-1}\leq c/\eta$. We can see in \Cref{table:SurfaceDecomposition} that these bounds are well satisfied. Recall that the bounded condition number of $M^{(k)}$ is essential for the efficiency of \Cref{alg:eigenpair_extension}, and the bounded condition number of $B^{(k)}$ is essential for the efficiency of \Cref{alg:eigenpair_refinement}.

\Cref{table:SurfaceDecomposition} also shows the detailed information for all four decompositions. The 2-norm of $A^{(k)}$, namely $\lambda_{max}(A^{(k)})$ decreases as $k$ increases, and well bounded as $\|A^{(k)}\|_2\leq \delta_k\leq c\varepsilon^{-1}_k$ as expected by \Cref{cor:Bconditonnumber} (we have normalized $\mu_1=\|L^{-1}\|_2$ to 1). And the sparsities of $A^{(k)}$ and $M^{(k)}$ are of the same order as the sparsity of $A^{(0)}=L$, i.e. $nnz(A^{(k)}),nnz(M^{(k)})=O(nnz(A^{(0)}))$ as we mentioned at the end of \Cref{subsec:energy_decomposable}. 
\begin{table}[!h]
\centering
%\scriptsize %\tiny % \footnotesize, \small, \normalsize
\footnotesize
\renewcommand{\arraystretch}{1.5}
    \begin{tabular}{|c|c|c|c|c|c|c|}
    \hline
    Data & $K$ & $\varepsilon_1$ & $\eta$ & $c$ & Bound on $\kappa(M^{(k)})$ & Bound on $\kappa(B^{(k)})$\\ \hline
    Bunny & 2 & $10^{-3}$ & 0.1 & 20 & 20 & 200 \\
    Brain & 4 & $10^{-4}$ & 0.2 & 20 & 20 & 100 \\
    SwissRoll & 3 & $10^{-5}$ & 0.1 & 20 & 20 & 200 \\
    SwissRoll & 4 & $10^{-5}$ & 0.2 & 20 & 20 & 100 \\
    \hline
    \end{tabular}
    \caption{Decomposition information}
    \label{table:Decomposition}
\end{table}

\begin{table}[!h]
\centering
%\scriptsize %\tiny % \footnotesize, \small, \normalsize
\footnotesize
\renewcommand{\arraystretch}{1.5}
    \begin{tabular}{|C{1cm}|C{1.6cm}|C{1.6cm}|C{2.5cm}|C{1.5cm}|C{2.5cm}|C{1.2cm}|C{1.1cm}|}
    \hline
    Level $k$ & $\varepsilon_k$ &Size of $A^{(k)}$ & $nnz(A^{(k)})$ & $\|A^{(k)}\|_2$  & $nnz(M^{(k)})$ & $\kappa(M^{(k)})$ &$\kappa(B^{(k)})$\\ 
    \hline
    \multicolumn{8}{|c|}{The 2-level decomposition of Bunny data.} \\
    \hline
    0 & - &$35947$ & $714647 = m$ & $1.86\times 10^4$ & - & - & -\\ 
    1 & $10^{-3}$ &$2641$ & $613571\approx0.86m$ & $1.05\times 10^4$ & $203445\approx0.28m$ & $1.45$ & $5.58$\\ 
    2 & $10^{-2}$ &$198$ & $27774\approx0.04m$ & $1.37\times 10^3$ & $10808\approx0.02m$ & $2.05$ & $45.03$\\
    \hline
    \hline
    \multicolumn{8}{|c|}{The 4-level decomposition of Brain data.} \\
    \hline
    0 & - &$48463$ & $1038065= m$ & $1.14\times 10^5$ & - & - & -\\ 
    1 & $10^{-4}$ & $11622$ & $2546246\approx2.45m$ & $7.82\times 10^4$ & $725328\approx0.70m$ & $1.29$ & $5.80$\\ 
    2 & $5\times10^{-4}$ & $1713$ & $431269\approx0.42m$ & $2.01\times 10^4$ & $189051\approx0.18m$ & $1.84$ & $18.34$\\
    3 & $2.5\times10^{-3}$ & $252$ & $37230\approx0.04m$ & $3.33\times 10^3$ & $20126\approx0.02m$ & $2.19$ & $28.23$\\
    4 & $1.25\times10^{-2}$ & $35$ & $1217<0.01m$ & $4.53\times 10^2$ & $1093<0.01m$ & $2.02$ & $35.08$\\
    \hline
    \hline
    \multicolumn{8}{|c|}{The 3-level decomposition of SwissRoll data.} \\
    \hline
    0 & - & $20000$ & $248010= m$ &  $1.15\times10^{6}$  & - & - & - \\ 
    1 & $10^{-5}$ &$5528$ & $689032\approx2.78m$ & $4.31\times10^{5}$  & $197020\approx0.79m$ & $1.45$ & $10.06$\\ 
    2 & $10^{-4}$ &$723$ & $108887\approx0.44m$ & $7.44\times10^{4}$ & $35213\approx0.14m$ & $2.30$ &$67.47$\\
    3 & $10^{-3}$ &$55$ & $2215<0.01m$ & $5.45\times10^{3}$ & $1365<0.01m$ & $3.92$ &$185.93$\\
    \hline
    \hline
    \multicolumn{8}{|c|}{The 4-level decomposition of SwissRoll data.} \\
    \hline
    0 & - & $20000$ & $248010= m$ &  $1.15\times10^{6}$  & - & - & - \\ 
    1 & $10^{-5}$ &$5528$ & $689032\approx2.78m$ & $4.31\times10^{5}$  & $197020\approx0.79m$ & $1.45$ & $10.06$\\ 
    2 & $5\times10^{-5}$ &$1347$ & $215811\approx0.87m$ & $9.36\times10^{4}$ & $65169\approx0.26m$ & $1.90$ &$26.52$\\
    3 & $2.5\times10^{-4}$ &$203$ & $18849\approx0.08m$ & $1.89\times10^{4}$ & $9063\approx0.04m$ & $3.06$ &$98.87$\\
    4 & $1.25\times10^{-3}$ &$53$ & $1939<0.01m$ & $3.72\times10^{3}$ & $1193<0.01m$ & $3.36$ & $51.14$\\
    \hline
    \end{tabular}
    \caption{Decomposition information of (i) Bunny (2-level) (ii) Brain (4-level) and (iii) SwissRoll (3, 4-level) data. $m\triangleq nnz(A^{(0)})$.}
    \label{table:SurfaceDecomposition}
\end{table}

\subsection{The Coarse Level Eigenpair Approximation}
We first use the decompositions given above to compute the first few eigenpairs of graph laplacians with relatively low accuracies. Even on the coarse levels, the compressed (low dimensional) operators show good spectral approximation properties with regard to the smallest eigenvalues of $L$ (or the largest eigenvalues of $L^{-1}$). Here we take the bunny data and the brain data as examples. For the bunny data, we use the lowest level $k=2$ with compression error $\varepsilon_2=0.01$; for the brain data, we use level $k=3$ with compression error $\varepsilon_3=0.0025$.

\begin{figure}[!h]
\centering
    \begin{subfigure}[b]{0.43\textwidth}
        \includegraphics[width=\textwidth]{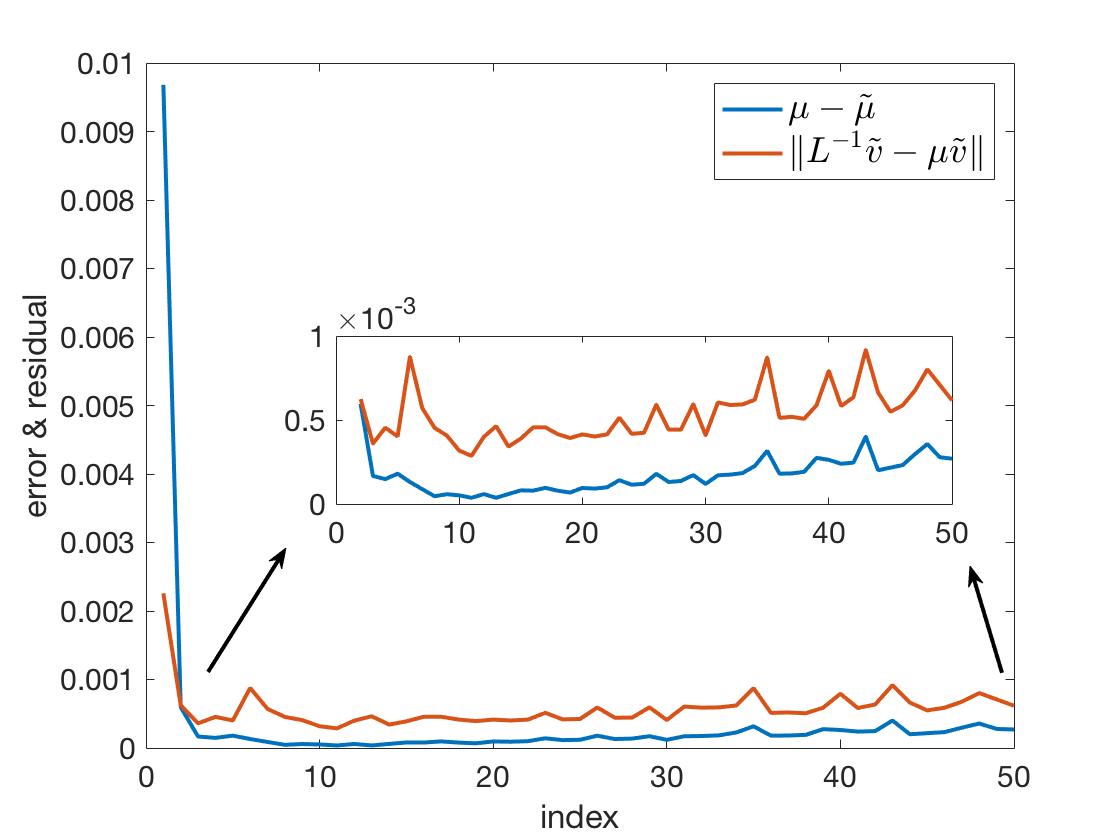}
        \caption{error $\mu_i-\tilde{\mu}_i$ and residual $\|L^{-1}\tilde{v}_i-\mu_i\tilde{v}_i\|_2$}
        \label{fig:bunny_coarseeigen}
    \end{subfigure}
    \begin{subfigure}[b]{0.43\textwidth}
        \includegraphics[width=\textwidth]{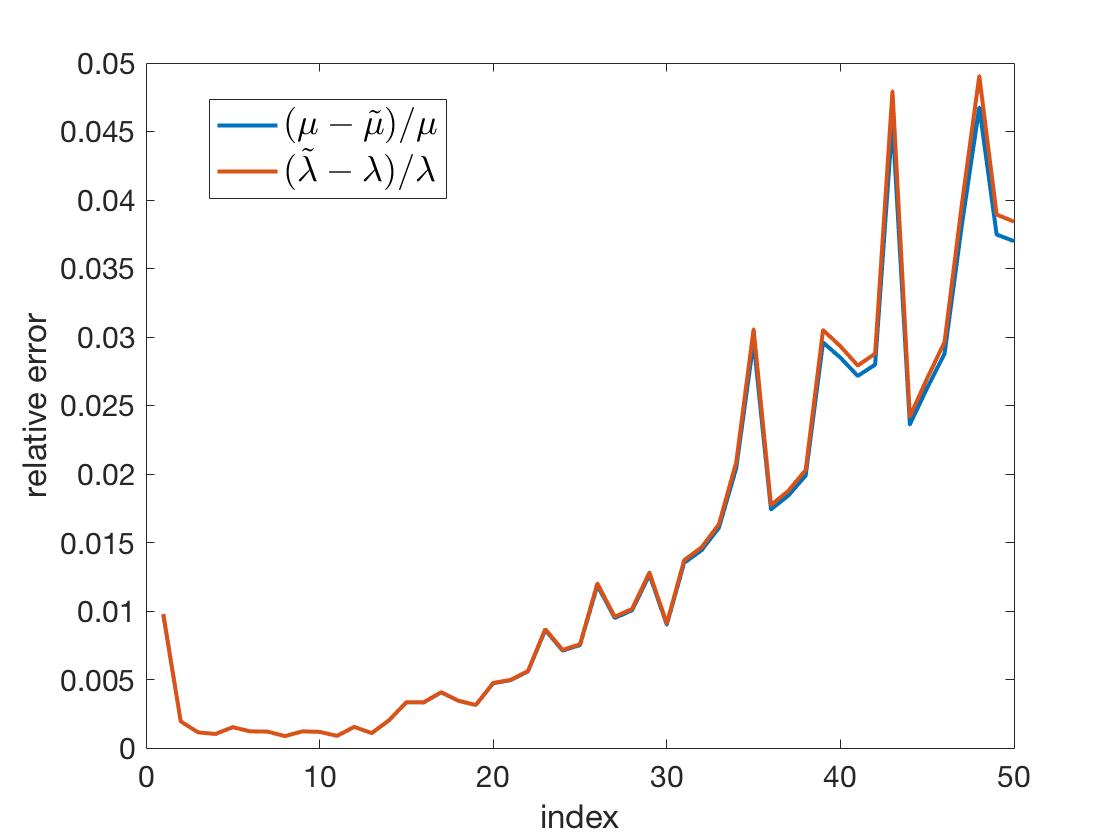}
        \caption{relative error of $\mu_i$ and $\lambda_i$}      
    \end{subfigure}
    \begin{subfigure}[b]{0.43\textwidth}
        \includegraphics[width=\textwidth]{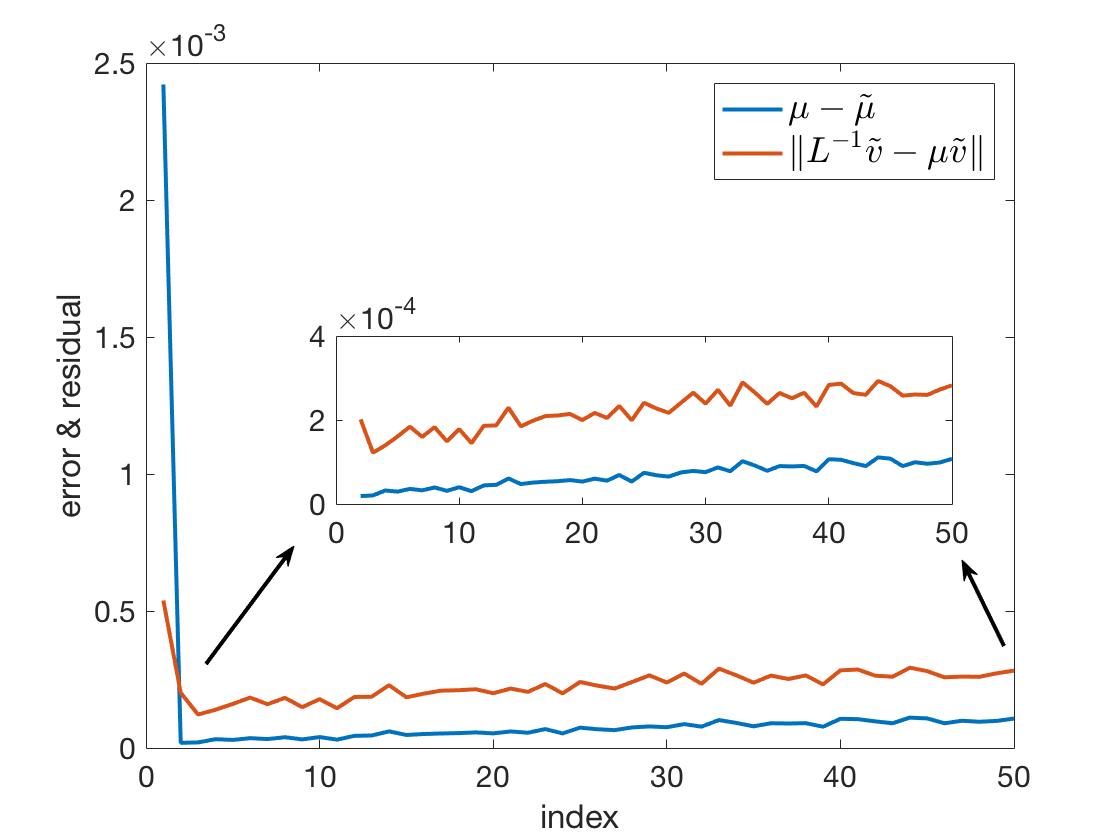}
        \caption{error $\mu_i-\tilde{\mu}_i$ and residual $\|L^{-1}\tilde{v}_i-\mu_i\tilde{v}_i\|_2$}  
        \label{fig:brain_coarseeigen}
    \end{subfigure}
	\begin{subfigure}[b]{0.43\textwidth}
        \includegraphics[width=\textwidth]{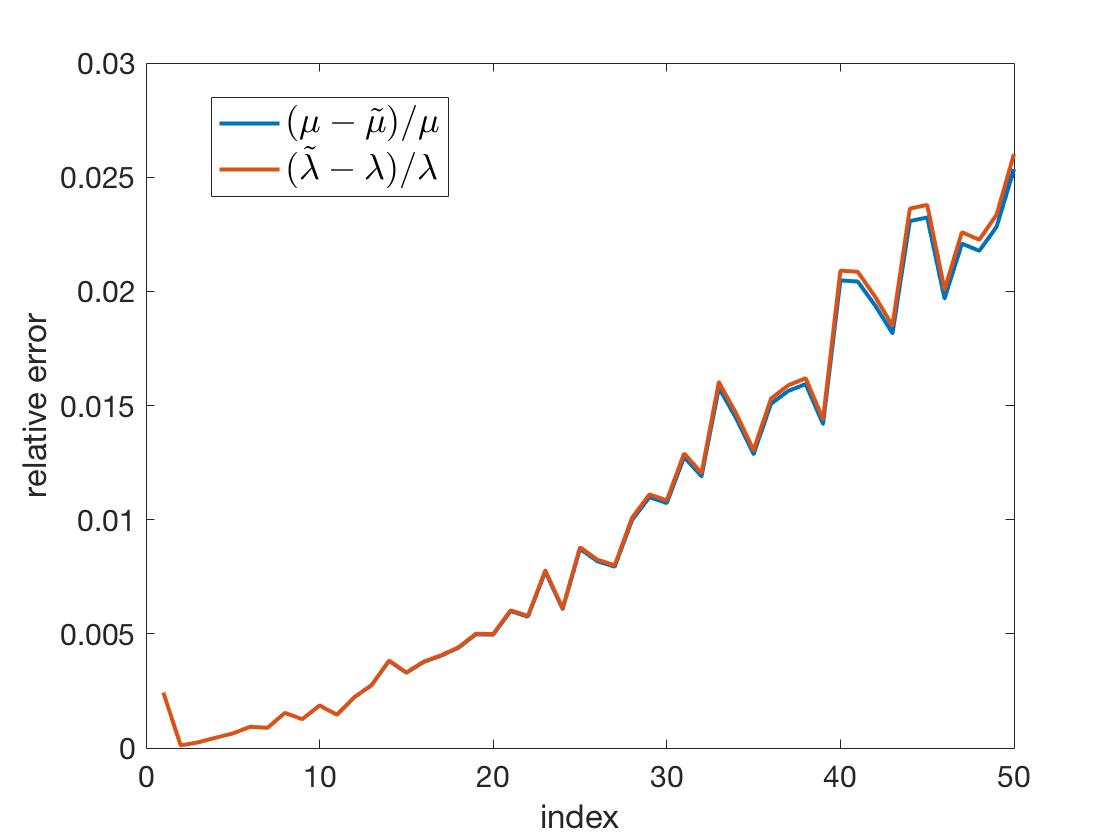}
        \caption{relative error of $\mu_i$ and $\lambda_i$}
    \end{subfigure}
    \caption{The error, the residual and the relative error. Top: Bunny data; bottom: Brain data}\label{fig:errors}
\end{figure}

We compute the first 50 eigenpairs $\{\tilde{v}_i,\tilde{\lambda}_i\}$ of the compressed operator by directly solving the general eigen problem (\Cref{lemma:generaleigenproblem})
\[A^{(k)}z_i=\tilde{\lambda}_iM^{(k)}z_i,\quad \tilde{v}_i=\bm{\Psi}^{(k)}z_i,\quad i=1,\cdots,50.\]
The computation of the coarse level eigenproblem is much more efficient due to the compressed dimension. To show the error of the approximate eigenvalues, the ground truth is obtained by using the Eigen C++ Library \footnotemark. \Cref{fig:errors} shows the absolute and relative errors of these eigenvalues. In both cases $\mu_i$ is the $i$th largest eigenvalue of $L^{-1}$ and $\lambda_i = 1/\mu_i$; $\tilde{\mu}_i$ is the $i$th largest eigenvalue of the compressed problem $\Theta^{(k)}$ and $\tilde{\lambda}_i = 1/\tilde{\mu}_i$. By \Cref{lemma:approximate_eigenpair_multiresolution}, $|\mu_i-\tilde{\mu}_i|$ is bounded by $\varepsilon_k$ and $\|L^{-1}\tilde{v}_i-\mu_i\tilde{v}_i\|_2$ is bounded $2\varepsilon_k$. We can see in \Cref{fig:errors} that both estimates are well satisfied. In particular, the error of the first eigenvalue is close to the bound of $\varepsilon_k$. However, the first eigenpair is already known. Therefore, we are only interested in the $2^{\text{nd}}$ up to $50^{\text{th}}$ eigenvalues and we embed the sub-plot of these eigenvalue errors as shown in \Cref{fig:bunny_coarseeigen} and \Cref{fig:brain_coarseeigen} respectively.

We can also qualitatively test the accuracy of the approximate eigenvectors of the compressed operators, by comparing their behaviors in image segmentation to those of the true eigenvectors of the original Laplacian operators. We will leave the detailed comparison to the Appendix.

\footnotetext{Eigen C++ Library is available at \url{http://eigen.tuxfamily.org/index.php?title=Main_Page}}

\subsection{The Multi-level Eigenpair Computation}
In this section, we use our main \Cref{alg:hierarch_eigencompute} to compute a relatively large number of eigenpairs of Laplacian matrices subject to the prescribed accuracy. For both the Brian data and the SwissRoll data, we compute the first $500$ eigenpairs of the graph Laplacian subject to prescribed accuracy $|\lambda^{-1}_i-\tilde{\lambda}^{-1}_i|=|\mu_i-\tilde{\mu}_i|\leq\epsilon=\varepsilon_1$.

The three decompositions of these two datasets are used in this section. For each decomposition, we apply \Cref{alg:hierarch_eigencompute} with two sets of parameters, $(\alpha,\beta)=(5,2)$ and $(\alpha,\beta)=(3,1)$. The details of the results that are obtained using \Cref{alg:hierarch_eigencompute} are summarized in \Cref{table:computationinformation}-\Cref{table:SwissRollEigencompute-2}. In \Cref{table:computationinformation}, parameters $\alpha,\beta,\kappa,\gamma$ are defined in \Cref{sec:overall_algorithm}. In \Cref{table:BrainEigencompute}-\Cref{table:SwissRollEigencompute-2}, we collect numerical results that reflect the efficiency of each single process (refinement or extension). Here we give a detailed description of the notations we use in these tables:
\begin{itemize}
\item \#I and \#O denote the numbers of input and output eigenpairs. To be consistent with the notations defined in \Cref{sec:overall_algorithm}, we use (\#I,\#O)$=(m_{k+1},\hat{m}_k)$ for refinement process on level $k$, and (\#I,\#O)$=(\hat{m}_k,m_k)$ for extension process on level $k$.
\item \#Iter denotes the number of orthogonal iterations in the refinement process. Note that this number is controlled by the ratio $\gamma$.
\item $\#_{\text{cg}}(B^{(k)})$ denotes number of CG calls concerning $B^{(k)}$ in the refinement process; $\#_{\text{pcg}}(A^{(k)})$ denotes the number of PCG calls concerning $A^{(k)}$ in the refinement process and the extension process. $\overline{\#}(B^{(k)})$ and $\overline{\#}(A^{(k)})$ denote the average numbers of matrix-vector multiplications concerning $B^{(k)},A^{(k)}$ respectively, namely the average numbers of iterations, in one single call of CG or PCG. Note that $\overline{\#}(B^{(k)})$ is controlled by $\log(1/\epsilon^{(k)})\kappa(B^{(k)})\leq \log(1/\epsilon^{(k)})c/\eta$, and $\overline{\#}(A^{(k)})$ by $\log(1/\epsilon^{(k)})\kappa(A_{\bm{\Psi}}^{(k)},Z^{(k)}_{\hat{m}_k^+})\leq \log(1/\epsilon^{(k)})\alpha c/\eta$.
\item As the extension process proceeds, the target spectrum to be computed on this level shrinks even more, and so does the restricted condition number of the operator. Thus the numbers of iterations in each PCG call get much smaller than its expected control $\log(1/\epsilon^{(k)})\alpha c/\eta$, which is a good thing in practice. So to study how the theoretical bound $\log(1/\epsilon^{(k)})\alpha c/\eta$ really affects the efficiency of PCG calls, it is more reasonable to investigate the maximal number of iterations in one PCG call on each level. We use $\widehat{\#}(A^{(k)})$ to denote the largest number of iterations in one single PCG call on level $k$.
\item $\overline{\#}(M^{(k)})$ denotes the average number of matrix-vector multiplications concerning $M^{(k)}$ in one single CG call concerning $M^{(k)}$. Such CG calls occur in the PCG calls concerning $A^{(k)}$ where $M^{(k)}$ acts as the preconditioner. Note that $\overline{\#}(M^{(k)})$ is controlled by $\log(1/\epsilon^{(k)})\kappa(M^{(k)})\leq \log(1/\epsilon^{(k)})(1+c)$.
\item ``Main Cost'' denotes the main computational cost contributed by matrix-vector multiplication flops. In the refinement process we have
\begin{equation*}
\begin{split}
\mathrm{Main\ Cost}=&\#_{\text{cg}}(B^{(k+1)})\cdot\overline{\#}(B^{(k+1)})\cdot nnz(A^{(k)})\\
&+\#_{\text{pcg}}(A^{(k)})\cdot\overline{\#}(A^{(k)})\cdot \big(nnz(A^{(k)})+\overline{\#}(M^{(k)})\cdot nnz(M^{(k)})\big),
\end{split}
\end{equation*}
while in the extension process we have
\[\mathrm{Main\ Cost}=\#_{\text{pcg}}(A^{(k)})\cdot\overline{\#}(A^{(k)})\cdot \big(nnz(A^{(k)})+\overline{\#}(M^{(k)})\cdot nnz(M^{(k)})\big).\]
\end{itemize}

\Cref{table:BrainEigencompute}-\Cref{table:SwissRollEigencompute-2} show the efficiency of our algorithm. We can see that $\overline{\#}(B^{(k)})$ and $\overline{\#}(M^{(k)})$ are well bounded as expected, due to the artificial imposition of the condition bound $c$. $\widehat{\#}(A^{(k)})$ and the numerical condition number $\widehat{\#}(A^{(k)})/\log(1/\epsilon^{(k)})$ are also well controlled by choosing $\alpha$ properly to bound $\kappa=\alpha c/\eta$. It is worth mentioning that $\widehat{\#}(A^{(k)})/\log(1/\epsilon^{(k)})$ appears to be uniformly bounded for all levels, actually much smaller than $\kappa$, which reflects our uniform control on efficiency. \#Iter is well bounded due to the proper choice of $\beta$ for bounding $\gamma=(1+\beta)/\alpha$.

We may also compare the results for the same decomposition but from two different sets of parameters $(\alpha,\beta)$. For all three decompositions, the experiments with $(\alpha,\beta)=(5,2)$ have a smaller $\gamma=\frac{3}{5}$, and thus is more efficient in the refinement process (less \#Iter and less refinement Main Cost). While the experiments with $(\alpha,\beta)=(3,1)$ have a smaller $\kappa$ that leads to better efficiency in the extension process (smaller $\widehat{\#}(A^{(k)})/\log(1/\epsilon^{(k)})$ and less extension Main Cost). But since the dominant cost of the whole process comes from the extension process, thus the experiments with $(\alpha,\beta)=(3,1)$ have a smaller Total Main Cost.

We remark that the choice of $(\alpha,\beta)$ not only determines $(\kappa,\gamma)$ that will affect the algorithm efficiency, but also determines the segmentation of the target spectrum and its allocation towards different levels of the decomposition. Smaller values of $\alpha$ and $\beta$ means more eigenpairs being computed on coarser levels (larger $k$), which relieves the burden of the extension process for finer levels, but also increases the load of the refinement process. There could be an optimal choice of $(\alpha,\beta)$ that minimizes the total main cost, balancing between the refinement and the extension processes. However, without a priori knowledge of the distribution of the eigenvalues, which is the case in practice, a safe choice of $(\alpha,\beta)$ would be $\alpha,\beta=O(1)$.

\begin{table}[!h]
\centering
%\scriptsize %\tiny % \footnotesize, \small, \normalsize
\footnotesize
\renewcommand{\arraystretch}{1.5}
    \begin{tabular}{|c|c|c|c|c|c|c|c|}
    \hline
    Data & Decomposition & $(\alpha,\beta)$ & $(\eta,c)$ & $\kappa$ & $\gamma$ & Total \#Iter & Total Main Cost \\ \hline
    Brain & 4-level & $(5,2)$ & $(0.2,20)$ & $500$ & $3/5$ & 12 & $4.37\times 10^5 \cdot m$ \\
     & 4-level & $(3,1)$ & $(0.2,20)$ & $300$ & $2/3$ & 15 & $4.13\times 10^5 \cdot m$\\ \hline
    SwissRoll & 3-level & $(5,2)$ & $(0.1,20)$ & $1000$ & $3/5$ & 13 & $7.56\times 10^5 \cdot m$\\
     & 3-level & $(3,1)$ & $(0.1,20)$ & $600$ & $2/3$ & 16 & $7.17\times 10^5 \cdot m$\\ \hline
    SwissRoll & 4-level & $(5,2)$ & $(0.2,20)$ & $500$ & $3/5$ & 19 & $7.00\times 10^5 \cdot m$\\
     & 4-level & $(3,1)$ & $(0.2,20)$ & $300$ & $2/3$ & 28 & $5.86\times 10^5 \cdot m$\\
    \hline
    \end{tabular}
    \caption{Computation information. $m\triangleq nnz(A^{(0)})$.}
    \label{table:computationinformation}
\end{table}

\begin{table}[h!]
\centering
\footnotesize
\renewcommand{\arraystretch}{2}
\begin{tabular}{|c|c|c|c|c|c|c|c|c|c|}
\hline 
\multicolumn{10}{|c|}{$(\alpha, \beta) = (5,2)$} \\
\hline
\multirow{4}{*}{\rotatebox{-90}{Refinement}}
& Level $k$ & (\#I,\#O) & \#Iter & $\#_{\text{cg}}(B^{(k+1)})$ & $\overline{\#}(B^{(k+1)})$ & $\#_{\text{pcg}}(A^{(k)})$ & $\overline{\#}(A^{(k)})$ & $\overline{\#}(M^{(k)})$  & Main Cost \\ \cline{2-10} 
& 3 	& $(7, 4)$ 		& 4 	& 7 					& 24.43
& 28 	& 10.97 	& 6.10 	& $5.66\times 10^1 \cdot m$ \\ \cline{2-10}
& 2 	& $(41, 17)$ 	& 4 	& 41 					& 25.90 
& 164 	& 16.26 	& 6.12 	& $4.50\times 10^3\cdot m$ \\ \cline{2-10}
& 1 	& $(207, 84)$ 	& 4 	& 207 					& 23.44
& 828 	& 19.17 	& 4.64 	& $1.02\times 10^5 \cdot m$ \\ \cline{1-10}
\multirow{4}{*}{\rotatebox{-90}{Extension}} 
& Level $k$ & (\#I,\#O) & $\widehat{\#}(A^{(k)})$ & $\epsilon^{(k)}$ & $\frac{\widehat{\#}(A^{(k)})}{log(1/\epsilon^{(k)})}$ & $\#_{\text{pcg}}(A^{(k)})$ & $\overline{\#}(A^{(k)})$ & $\overline{\#}(M^{(k)})$ & Main Cost \\ \cline{2-10} 
& 3 	& $(4, 41)$ 		& 43 	& $2.5\times 10^{-4} $ 	& 5.18
& 175 	& 16.93 	& 5.39 	& $4.37\times 10^2 \cdot m$ \\ \cline{2-10} 
& 2 	& $(17, 207)$	& 75 	& $5.0\times 10^{-5}$ 	& 7.57
& 500 	& 32.27 	& 5.47 	& $2.27\times 10^4 \cdot m$ \\ \cline{2-10} 
& 1 	& $(84, 500)$	& 82 	& $10^{-5}$ 			& 7.12
& 1248 	& 44.23 	& 4.45 	& $3.07\times 10^5 \cdot m$ \\
\hline
\hline
\multicolumn{10}{|c|}{$(\alpha, \beta) = (3,1)$} \\
\hline
\multirow{4}{*}{\rotatebox{-90}{Refinement}}
& Level $k$ & (\#I,\#O) & \#Iter & $\#_{\text{cg}}(B^{(k+1)})$ & $\overline{\#}(B^{(k+1)})$ & $\#_{\text{pcg}}(A^{(k)})$ & $\overline{\#}(A^{(k)})$ & $\overline{\#}(M^{(k)})$  & Main Cost \\ \cline{2-10} 
& 3 	& $(15, 6)$ 		& 5 	& 15 					& 24.54
& 75 	& 7.74 		& 6.07 	& $1.08\times 10^2 \cdot m$ \\ \cline{2-10}
& 2 	& $(78, 28)$ 	& 5 	& 78 					& 25.85
& 390 	& 11.17 	& 6.01 	& $7.39\times 10^3 \cdot m$ \\ \cline{2-10}
& 1 	& $(276, 140)$ 	& 5 	& 276 					& 23.43
& 1380 	& 14.28 	& 4.67 	& $1.29\times 10^5 \cdot m$ \\ \cline{1-10}
\multirow{4}{*}{\rotatebox{-90}{Extension}} 
& Level $k$ & (\#I,\#O) & $\widehat{\#}(A^{(k)})$ & $\epsilon^{(k)}$ & $\frac{\widehat{\#}(A^{(k)})}{log(1/\epsilon^{(k)})}$ & $\#_{\text{pcg}}(A^{(k)})$ & $\overline{\#}(A^{(k)})$ & $\overline{\#}(M^{(k)})$ & Main Cost \\ \cline{2-10} 
& 3 	& $(6, 78)$ 		& 37 	& $2.5\times 10^{-4}$ 	& 4.46
& 225 	& 14.12 	& 5.41 	& $4.70\times 10^2 \cdot m$ \\ \cline{2-10}
& 2 	& $(28, 276)$ 	& 57 	& $5.0\times 10^{-5}$ 	& 5.75
& 600 	& 27.91 	& 5.43 	& $2.34\times 10^4 \cdot m$ \\ \cline{2-10}
& 1 	& $(140, 500)$ 	& 63 	& $10^{-5}$  			& 5.47
& 1080 	& 42.09 	& 4.46 	& $2.53\times 10^5 \cdot m$ \\
\hline
\end{tabular}
\caption{4-level eigenpairs computation of Brain data with $(\eta,c)=(0.2,20)$, $m\triangleq nnz(A^{(0)})$.}
\label{table:BrainEigencompute}
\end{table}

\begin{table}[h!]
\centering
\footnotesize
\renewcommand{\arraystretch}{2}
\begin{tabular}{|c|c|c|c|c|c|c|c|c|c|}
\hline
\multicolumn{10}{|c|}{$(\alpha, \beta) = (5,2)$} \\
\hline
\multirow{3}{*}{\rotatebox{-90}{Refinement}}
& Level $k$ & (\#I,\#O) & \#Iter & $\#_{\text{cg}}(B^{(k+1)})$ & $\overline{\#}(B^{(k+1)})$ & $\#_{\text{pcg}}(A^{(k)})$ & $\overline{\#}(A^{(k)})$ & $\overline{\#}(M^{(k)})$  & Main Cost \\ \cline{2-10} 
& 2 & $(21, 12)$ 	& 7 & 21 & 52.14 & 147 & 17.61 & 6.33 & $3.91 \times 10^3 \cdot m$ \\ \cline{2-10} 
& 1 & $(232, 100)$ 	& 6 & 232 & 47.23 & 1392 & 16.08 & 5.29 & $1.86 \times 10^5 \cdot m$ \\ \cline{1-10} 
\multirow{3}{*}{\rotatebox{-90}{Extension}} 
& Level $k$ & (\#I,\#O) & $\widehat{\#}(A^{(k)})$ & $\epsilon^{(k)}$ & $\frac{\widehat{\#}(A^{(k)})}{log(1/\epsilon^{(k)})}$ & $\#_{\text{pcg}}(A^{(k)})$ & $\overline{\#}(A^{(k)})$ & $\overline{\#}(M^{(k)})$ & Main Cost \\ \cline{2-10} 
& 2 & $(12, 232)$ 	& 94 & $10^{-5}$ & 8.16 & 650 & 28.20 & 7.25 & $2.67 \times 10^4 \cdot m$ \\ \cline{2-10} 
& 1 & $(100, 500)$ 	& 101 & $10^{-6}$ & 7.31 & 1200 & 59.44 & 6.10 & $5.42 \times 10^5 \cdot m$ \\
\hline
\hline
\multicolumn{10}{|c|}{$(\alpha, \beta) = (3,1)$} \\
\hline
\multirow{3}{*}{\rotatebox{-90}{Refinement}}
& Level $k$ & (\#I,\#O) & \#Iter & $\#_{\text{cg}}(B^{(k+1)})$ & $\overline{\#}(B^{(k+1)})$ & $\#_{\text{pcg}}(A^{(k)})$ & $\overline{\#}(A^{(k)})$ & $\overline{\#}(M^{(k)})$  & Main Cost \\ \cline{2-10} 
& 2 & $(35, 19)$ 	& 8 & 35 & 51.89 & 280 & 13.13 & 6.45 & $5.74 \times 10^3 \cdot m$ \\ \cline{2-10} 
& 1 & $(315, 165)$ 	& 8 & 315 & 46.85 & 2520 & 12.73 & 5.37 & $2.66 \times 10^5 \cdot m$ \\ \cline{1-10} 
\multirow{3}{*}{\rotatebox{-90}{Extension}}
& Level $k$ & (\#I,\#O) & $\widehat{\#}(A^{(k)})$ & $\epsilon^{(k)}$ & $\frac{\widehat{\#}(A^{(k)})}{log(1/\epsilon^{(k)})}$ & $\#_{\text{pcg}}(A^{(k)})$ & $\overline{\#}(A^{(k)})$ & $\overline{\#}(M^{(k)})$ & Main Cost \\ \cline{2-10} 
& 2 & $(19, 315)$ 	& 69 & $10^{-5}$ & 5.99 & 700 & 25.10 & 7.29 & $2.57 \times 10^4 \cdot m$ \\ \cline{2-10} 
& 1 & $(165, 500)$ 	& 78 & $10^{-6}$ & 5.65 & 1005 & 54.91 & 6.11 & $4.20 \times 10^5 \cdot m$ \\
\hline
\end{tabular}
\caption{3-level eigenpairs computation of SwissRoll data with $(\eta,c)=(0.1,20)$, $m\triangleq nnz(A^{(0)})$.}
\label{table:SwissRollEigencompute-1}
\end{table}

\begin{table}[h!]
\centering
\footnotesize
\renewcommand{\arraystretch}{2}
\begin{tabular}{|c|c|c|c|c|c|c|c|c|c|}
\hline
\multicolumn{10}{|c|}{$(\alpha, \beta) = (5,2)$} \\
\hline
\multirow{4}{*}{\rotatebox{-90}{Refinement}}
& Level $k$ & (\#I,\#O) & \#Iter & $\#_{\text{cg}}(B^{(k+1)})$ & $\overline{\#}(B^{(k+1)})$ & $\#_{\text{pcg}}(A^{(k)})$ & $\overline{\#}(A^{(k)})$ & $\overline{\#}(M^{(k)})$  & Main Cost \\ \cline{2-10} 
& 3 & $(18, 10)$ & 6 & 18 & 22.61 & 108 & 7.19 & 7.87 & $3.39 \times 10^2 \cdot m$ \\ \cline{2-10} 
& 2 & $(84, 44)$ & 8 & 84 & 43.45 & 672 & 10.42 & 6.49 & $2.11 \times 10^4 \cdot m$ \\ \cline{2-10} 
& 1 & $(390, 195)$ & 5 & 390 & 28.85 & 1950 & 11.68 & 5.42 & $1.92 \times 10^5 \cdot m$ \\ \cline{1-10} 
\multirow{4}{*}{\rotatebox{-90}{Extension}} 
& Level $k$ & (\#I,\#O) & $\widehat{\#}(A^{(k)})$ & $\epsilon^{(k)}$ & $\frac{\widehat{\#}(A^{(k)})}{log(1/\epsilon^{(k)})}$ & $\#_{\text{pcg}}(A^{(k)})$ & $\overline{\#}(A^{(k)})$ & $\overline{\#}(M^{(k)})$ & Main Cost \\ \cline{2-10} 
& 3 & $(10, 84)$ & 42 & $2.5\times 10^{-5}$ & 3.96 & 200 & 18.32 & 8.43 & $1.53 \times 10^3 \cdot m$ \\ \cline{2-10} 
& 2 & $(44, 390)$ & 63 & $5\times 10^{-6}$ & 5.16 & 1050 & 29.30 & 7.24 & $8.47 \times 10^4 \cdot m$ \\ \cline{2-10} 
& 1 & $(195, 50)$ & 71 & $10^{-6}$ & 5.13 & 915 & 57.47 & 6.10 & $4.00 \times 10^5 \cdot m$ \\
\hline
\hline
\multicolumn{10}{|c|}{$(\alpha, \beta) = (3,1)$} \\
\hline
\multirow{4}{*}{\rotatebox{-90}{Refinement}} 
& Level $k$ & (\#I,\#O) & \#Iter & $\#_{\text{cg}}(B^{(k+1)})$ & $\overline{\#}(B^{(k+1)})$ & $\#_{\text{pcg}}(A^{(k)})$ & $\overline{\#}(A^{(k)})$ & $\overline{\#}(M^{(k)})$  & Main Cost \\ \cline{2-10} 
& 3 & $(31, 16)$ & 7 & 31 & 22.45 & 217 & 6.09 & 8.09 & $5.89 \times 10^2 \cdot m$ \\ \cline{2-10} 
& 2 & $(95, 67)$ & 12 & 95 & 43.44 & 1140 & 7.66 & 6.66 & $2.63 \times 10^4 \cdot m$ \\ \cline{2-10} 
& 1 & $(459, 314)$ & 7 & 459 & 28.75 & 3656 & 8.71 & 5.56 & $2.65 \times 10^5 \cdot m$ \\ \cline{1-10} 
\multirow{4}{*}{\rotatebox{-90}{Extension}} 
& Level $k$ & (\#I,\#O) & $\widehat{\#}(A^{(k)})$ & $\epsilon^{(k)}$ & $\frac{\widehat{\#}(A^{(k)})}{log(1/\epsilon^{(k)})}$ & $\#_{\text{pcg}}(A^{(k)})$ & $\overline{\#}(A^{(k)})$ & $\overline{\#}(M^{(k)})$ & Main Cost \\ \cline{2-10} 
& 3 & $(16, 95)$ & 31 & $2.5\times 10^{-5}$ & 2.92 & 200 & 16.61 & 8.48 & $1.39 \times 10^3 \cdot m$ \\ \cline{2-10} 
& 2 & $(67, 459)$ & 49 & $5\times 10^{-6}$ & 4.01 & 1100 & 25.66 & 7.27 & $7.79 \times 10^4 \cdot m$ \\ \cline{2-10} 
& 1 & $(314, 500)$ & 55 & $10^{-6}$ & 3.98 & 558 & 50.61 & 6.12 & $2.15 \times 10^5 \cdot m$ \\
\hline
\end{tabular}
\caption{4-level eigenpairs computation of SwissRoll data: $(\eta,c)=(0.2,20)$, $m\triangleq nnz(A^{(0)})$.}
\label{table:SwissRollEigencompute-2}
\end{table}

\begin{figure}[!h]
\centering
    \begin{subfigure}[b]{0.43\textwidth}
        \includegraphics[width=\textwidth]{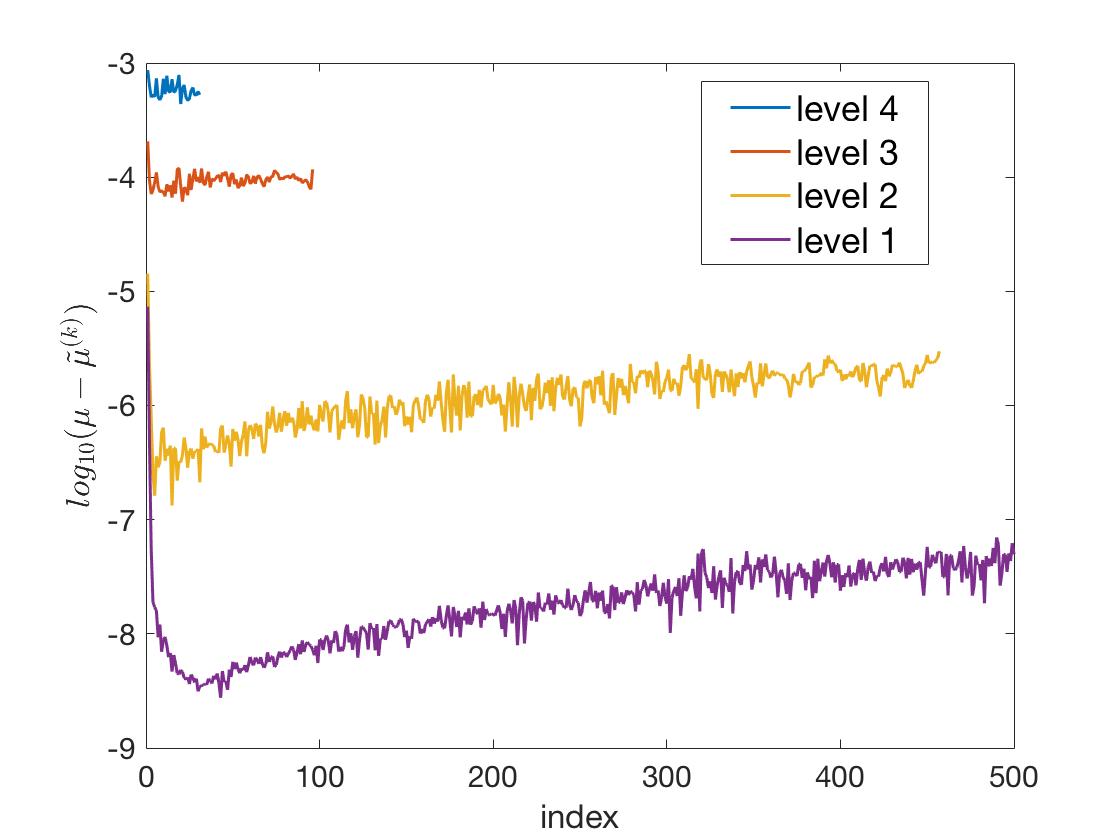}
        \caption{$\log_{10}(\mu_i-\tilde{\mu}_i^{(k)}),\ i=1,\cdots,m_k$}       
    \end{subfigure}
    \begin{subfigure}[b]{0.43\textwidth}
        \includegraphics[width=\textwidth]{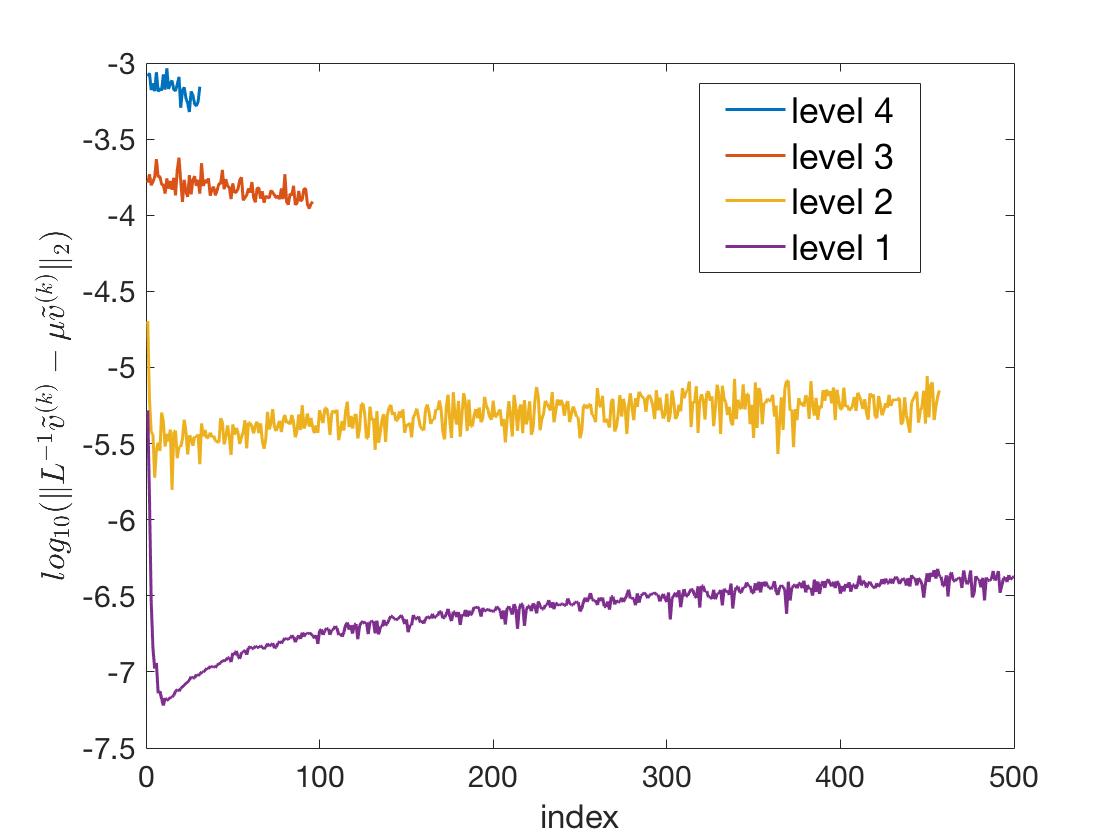}
        \caption{$\log_{10}(\|L^{-1}\tilde{v}_i^{(k)}-\mu_i\tilde{v}_i^{(k)}\|_2),\ i=1,\cdots,m_k$}      
    \end{subfigure}
    \begin{subfigure}[b]{0.43\textwidth}
        \includegraphics[width=\textwidth]{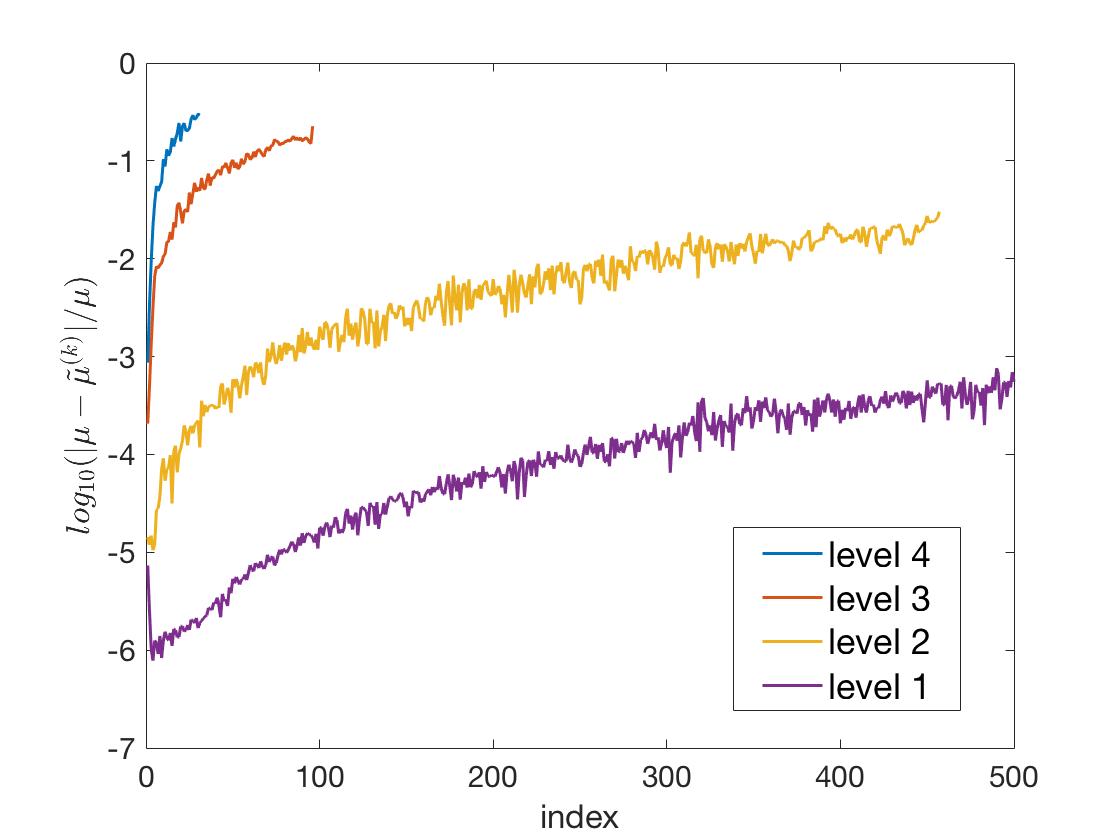}
        \caption{$\log_{10}\big((\mu_i-\tilde{\mu}_i^{(k)})/\mu_i\big),\ i=1,\cdots,m_k$}       
    \end{subfigure}
	\begin{subfigure}[b]{0.43\textwidth}
        \includegraphics[width=\textwidth]{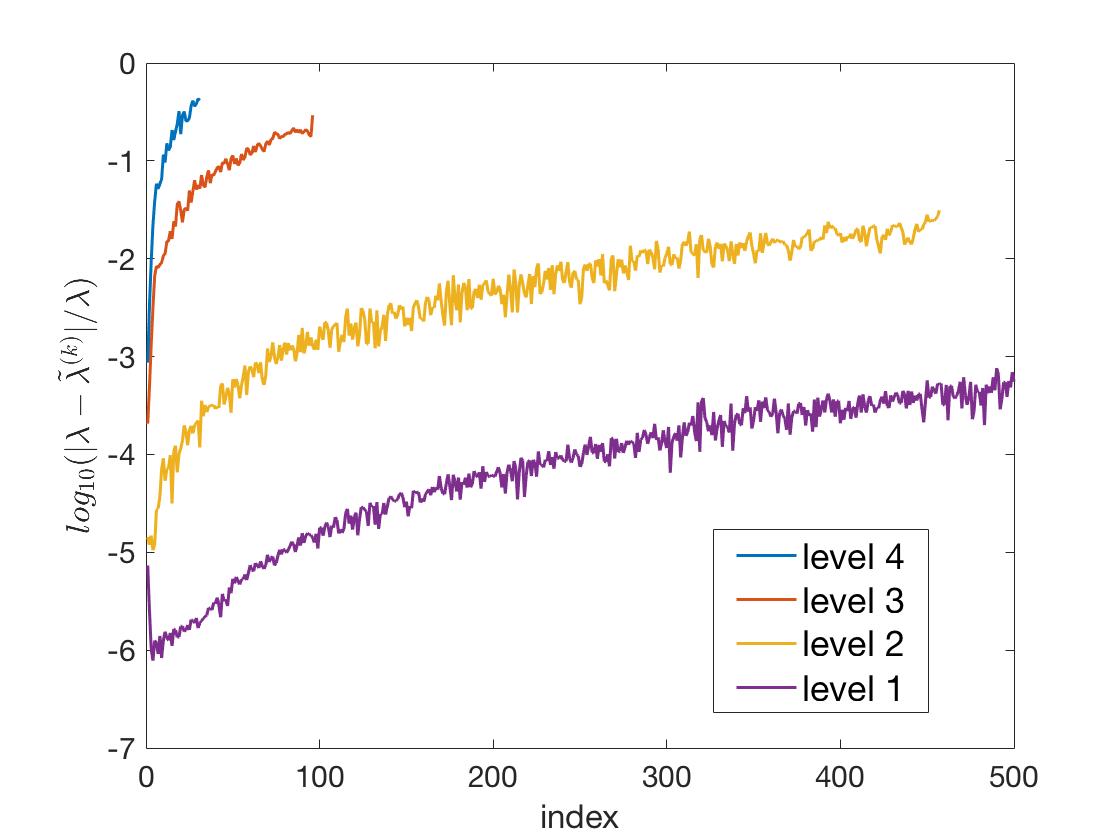}
        \caption{$\log_{10}\big((\lambda_i-\tilde{\lambda}_i^{(k)})/\lambda_i\big),\ i=1,\cdots,m_k$}
    \end{subfigure}
    \caption{Convergence of computed spectrum in different errors.}\label{fig:convergence}
\end{figure}

\begin{figure}[!]
\centering
    \begin{subfigure}[b]{0.32\textwidth}
        \includegraphics[width=\textwidth]{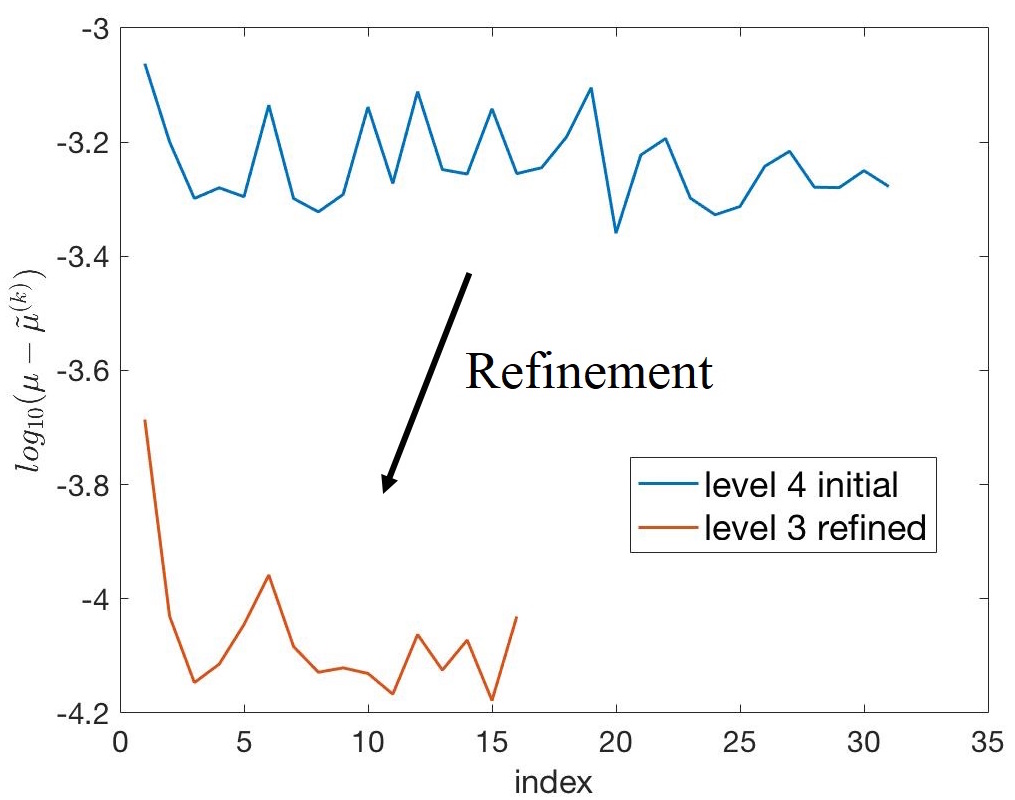}
        \caption{results after level 3 refinement}       
    \end{subfigure}
    \begin{subfigure}[b]{0.32\textwidth}
        \includegraphics[width=\textwidth]{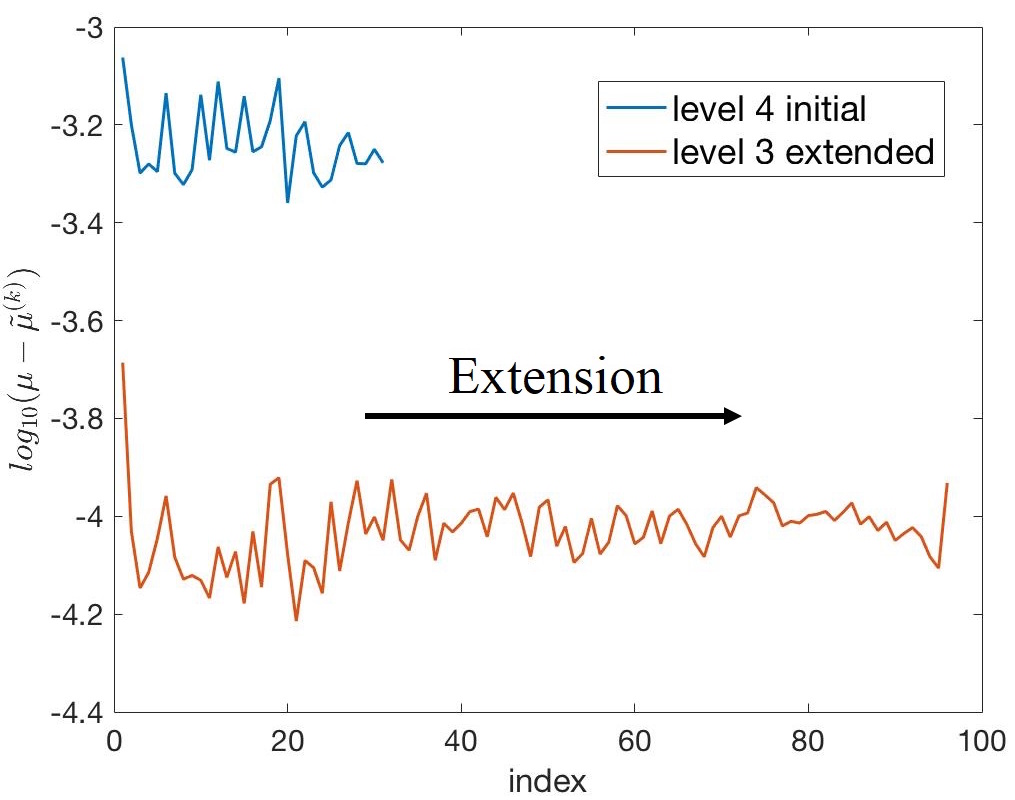}
        \caption{results after level 3 exension}      
    \end{subfigure}
    \begin{subfigure}[b]{0.32\textwidth}
        \includegraphics[width=\textwidth]{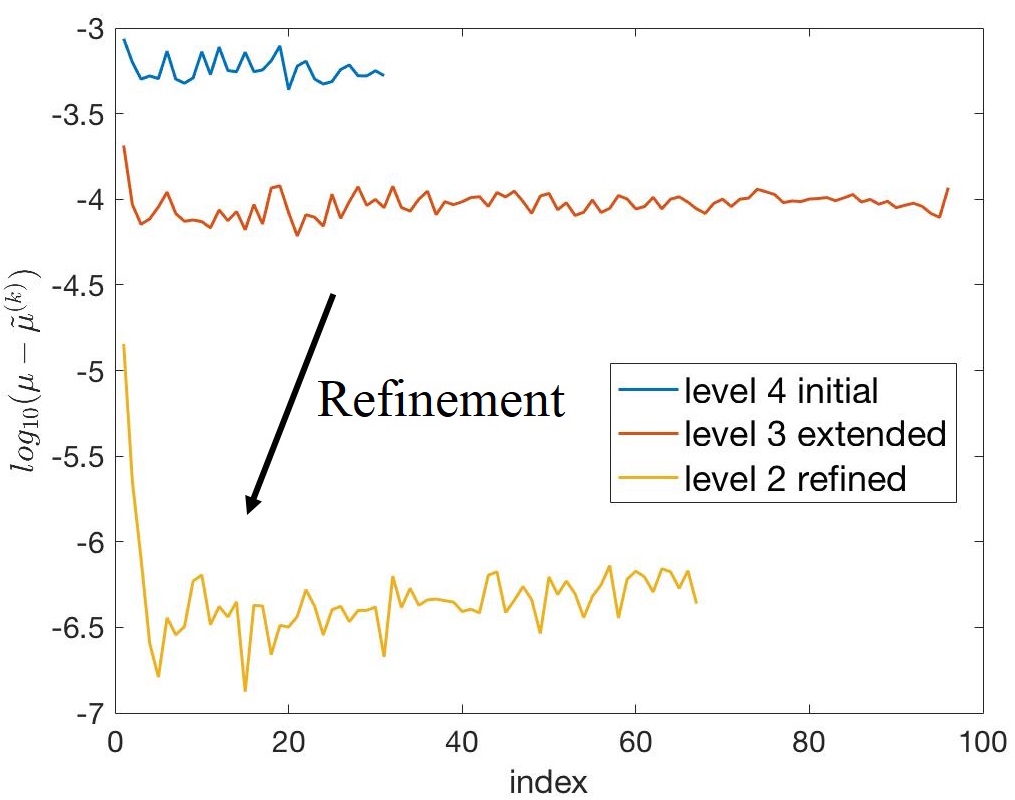}
        \caption{results after level 2 refinement}       
    \end{subfigure}
	\begin{subfigure}[b]{0.32\textwidth}
        \includegraphics[width=\textwidth]{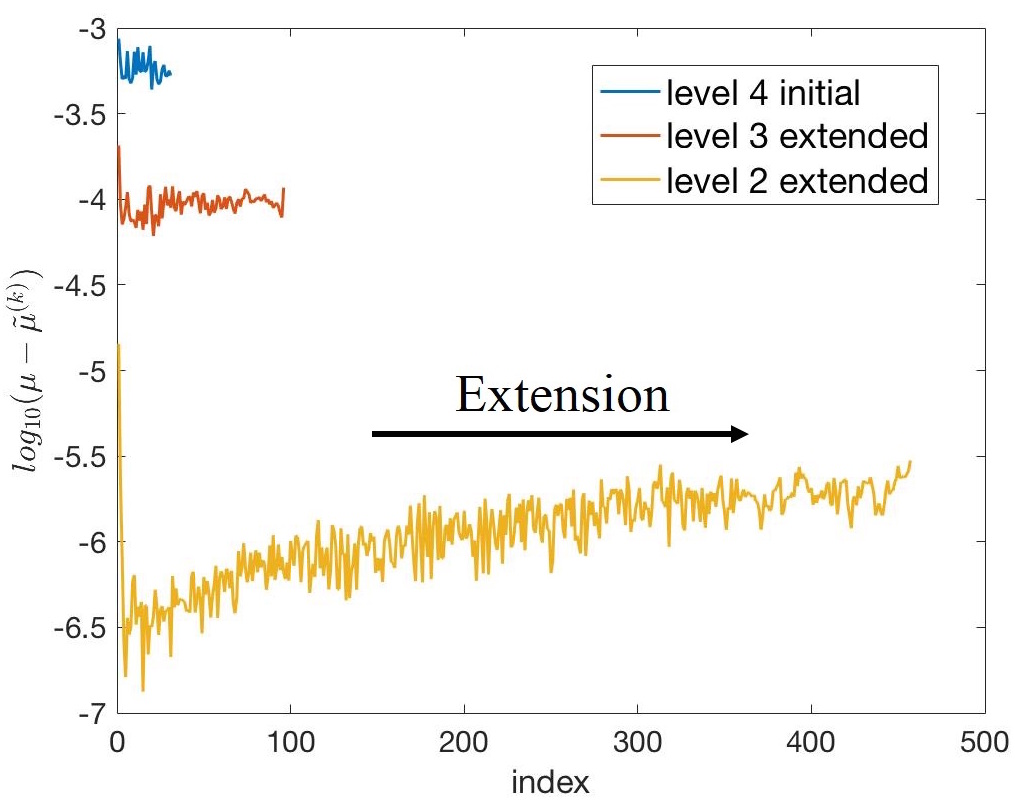}
        \caption{results after level 2 extension}
    \end{subfigure}
    \begin{subfigure}[b]{0.32\textwidth}
        \includegraphics[width=\textwidth]{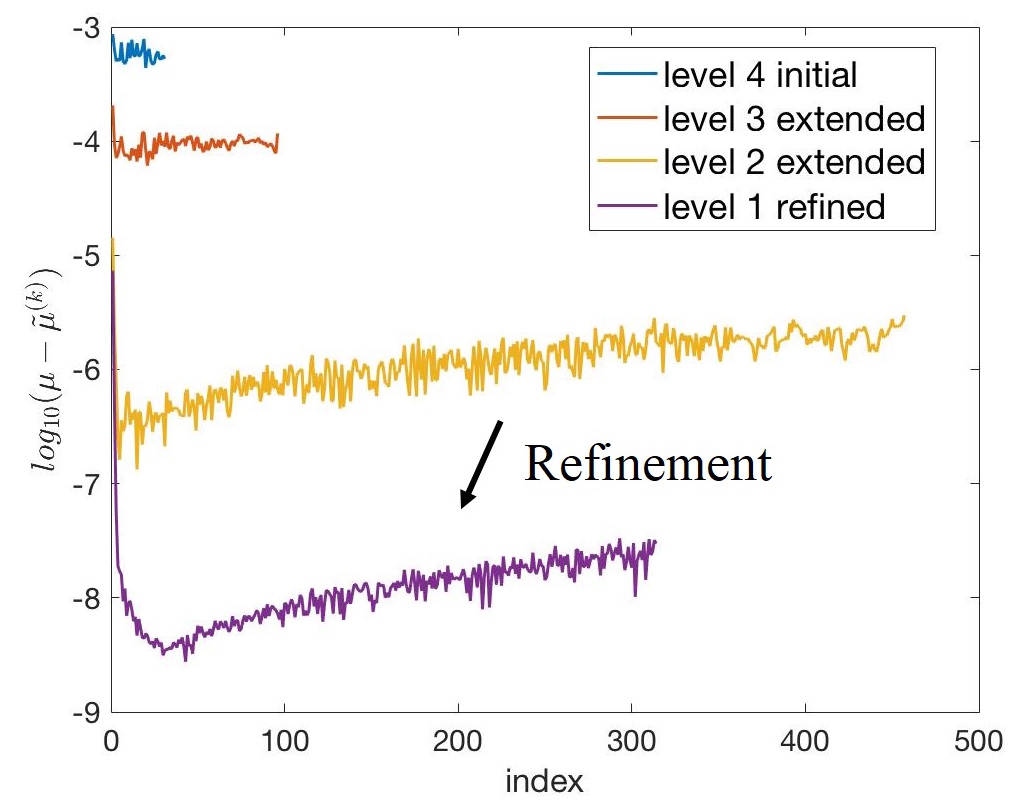}
        \caption{results after level 1 refinement}       
    \end{subfigure}
	\begin{subfigure}[b]{0.32\textwidth}
        \includegraphics[width=\textwidth]{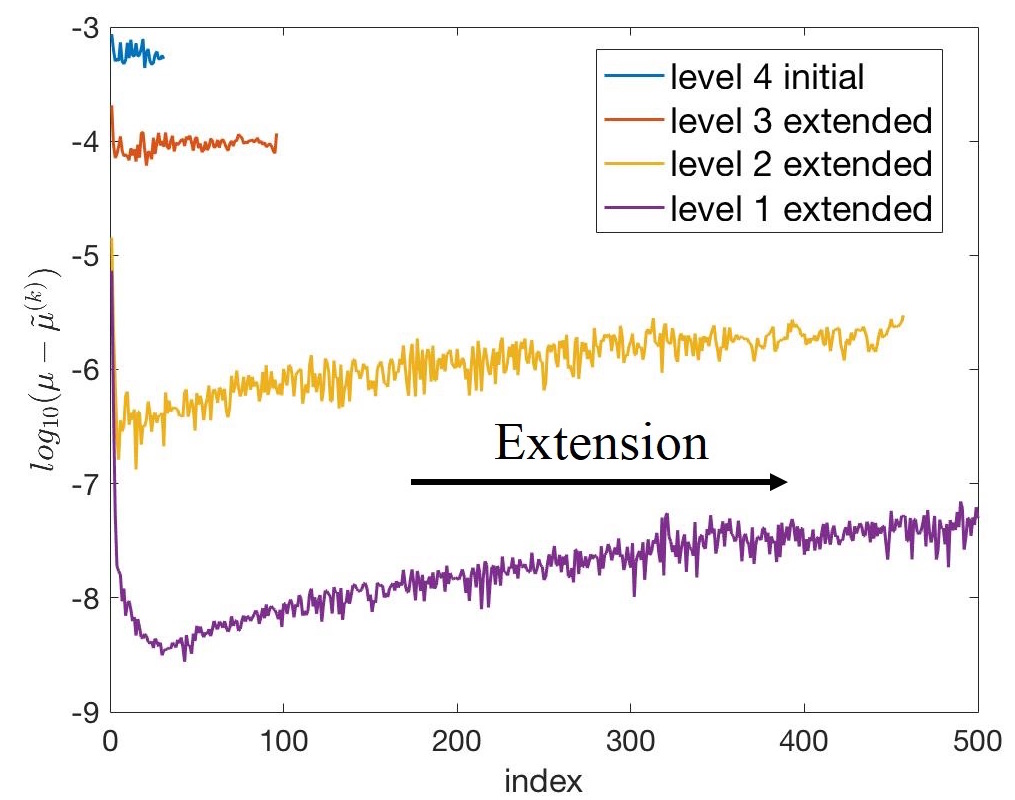}
        \caption{final results after level 1 extension}
    \end{subfigure}
    \caption{The completion and convergence process of the target spectrum. The refinement process retains part of the spectrum subject to threshold $\mu_{re}^{(k)}$ with improved accuracy, and the extension process extends the spectrum subject to threshold $\mu_{ex}^{(k)}$. The whole process is an iterative procedure that aims at improving the accuracy of the eigenvalue solver.}\label{fig:convergenceprocess}
\end{figure}

To further investigate the behavior of our algorithm, we focus on numerical experiments carried out on the 4-level decomposition of the SwissRoll data. \Cref{fig:convergence} shows the convergence of the computed spectrum in different errors. \Cref{fig:convergenceprocess} shows the completion and the convergence process of the target spectrum in the case of $(\alpha,\beta)=(3,1)$ (corresponding to \Cref{table:SwissRollEigencompute-2}). We use a log-scale plot to illustrate the error $|\mu_i-\tilde{\mu}^{(k)}|$ after we complete the refinement process and the extension process respectively on each level $k$. As we can see, each application of the refinement process improves the accuracy of the first $\hat{m}_k$ eigenvalues at least by a factor of $\eta=\frac{\varepsilon_k}{\varepsilon_{k+1}}$, but at the price of discarding the last $m_{k+1}-\hat{m}_k$ computed eigenvalues. So the computation of the last $m_{k+1}-\hat{m}_k$ computed eigenvalues on the coarser level $k+1$ actually serves as preconditioning to ensure the efficiency of the refinement process on level $k$. Then the extension process extends the spectrum to  $m_k$ that is determined by the threshold $\mu_{ex}^{(k)}$. The whole computation is an iterative process that improves the accuracy of the eigenvalues by applying the hierarchical Lanczos method to each eigenvalue at most twice.

\begin{figure}[h!]
\centering
\includegraphics[width=0.6\textwidth]{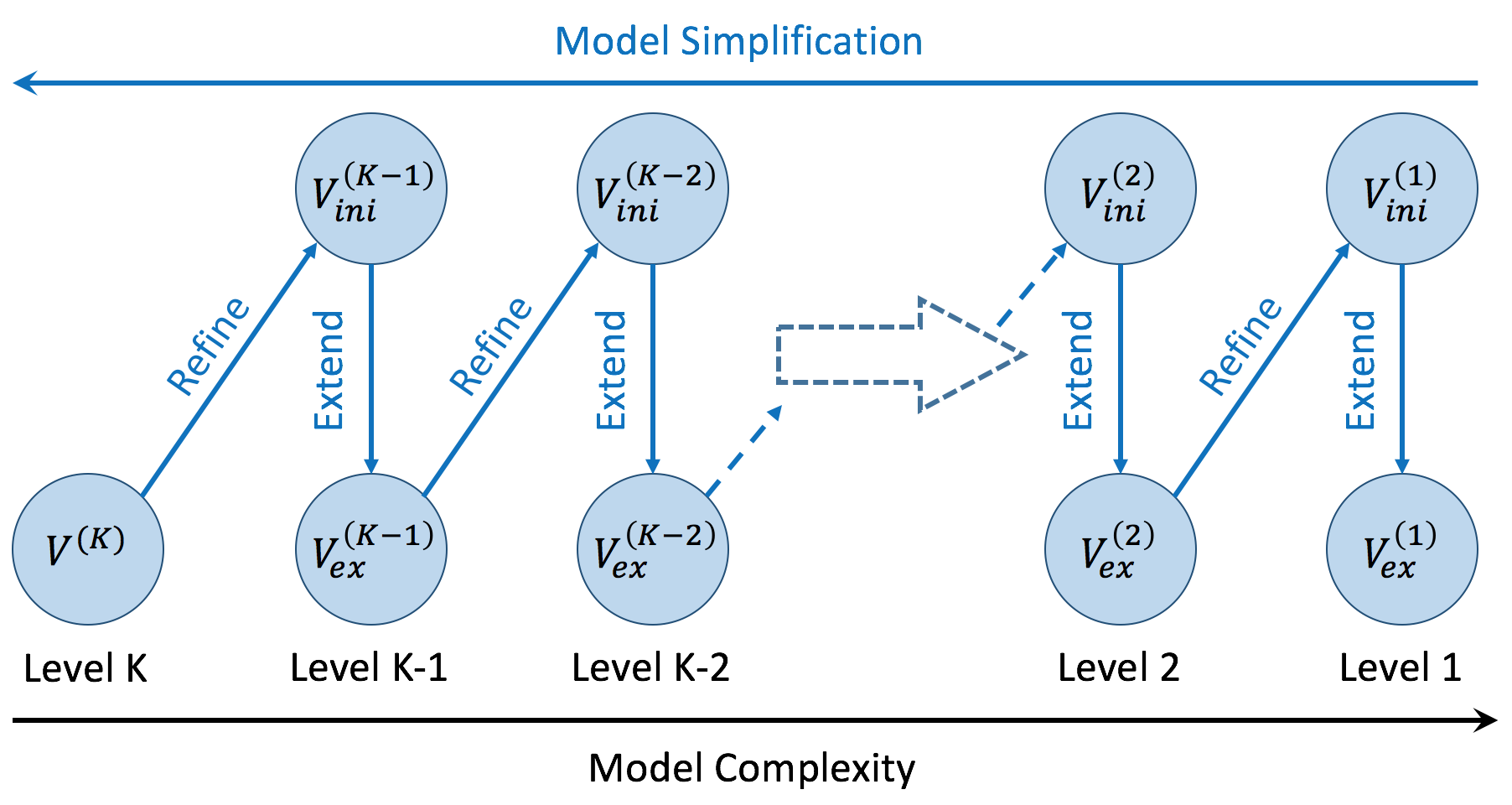} 
\caption{Flow Chart illustrating the procedure of \Cref{alg:hierarch_eigencompute}.}
\label{fig:flowchart}
\end{figure}

It could be clearer using a flow chart \Cref{fig:flowchart} to illustrate the procedure of our method. We can see the eigenproblem of the original matrix $A$ as a complicated model, and we are pursuing some solutions from this model. To resolve the complexity, we first use the multiresolution matrix decomposition to hierarchically simplify/coarsen the original model into a sequence of approximate models, so the model in each level $k$ is a simplification of the model in the higher level $k-1$. Then we start from the bottom level. Every time we obtain some partial solutions on an intermediate level, we feed them to the higher level through some correction process, and use the corrected ones to help us continue to complete the whole solution set.

\begin{figure}[h!]
\centering
    \begin{subfigure}[b]{0.43\textwidth}
        \includegraphics[width=\textwidth]{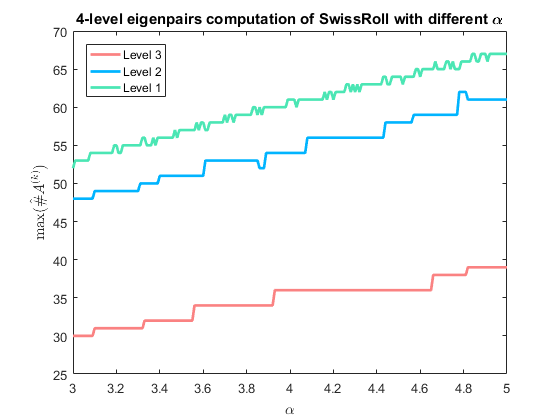}
        \caption{$\hat{\#} A^{(k)}$ versus $\alpha$}
        \label{fig:maxA_vs_alpha}
    \end{subfigure}
	\begin{subfigure}[b]{0.43\textwidth}
        \includegraphics[width=\textwidth]{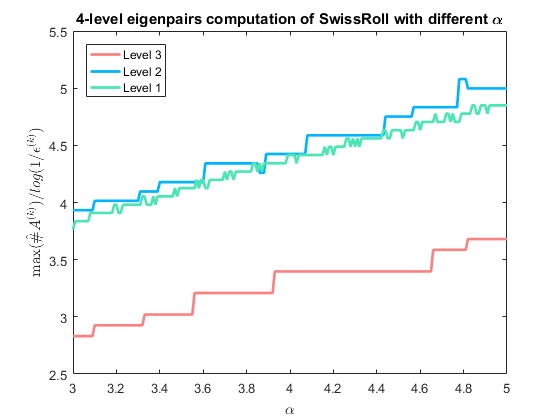}
        \caption{$\hat{\#} A^{(k)} / \log(\epsilon_{pcg}^{-1})$ versus $\alpha$}
        \label{fig:maxA_normalize_vs_alpha}
    \end{subfigure}
    \caption{$\hat{\#} A^{(k)}$ versus $\alpha$ in the 4-level SwissRoll example.}\label{fig:maxA}
\end{figure}

We also further verify our critical control on the restricted condition number $\kappa(A_{\bm{\Psi}}^{(k)},Z^{(k)}_{\hat{m}_k^+})$ by $\kappa=\alpha c/\eta$, by showing the dependence of $\widehat{\#}(A^{(k)})$(or $\widehat{\#}(A^{(k)})/\log(1/\epsilon^{(k)})$) on $\kappa$. Recall that $\widehat{\#}(A^{(k)})$ denotes the largest number of iterations in one single PCG call concerning $A^{(k)}$ on level k. Using the 4-level decomposition of the SwissRoll data with $(\eta,c)=(0.2,20)$, we perform \Cref{alg:hierarch_eigencompute} with fixed $\beta=1$ but different $\alpha\in[3,5]$. \Cref{fig:maxA} shows $\widehat{\#}(A^{(k)})$ versus $\alpha$ for all three levels. By \Cref{thm:multi_spectrum_completion_overall}, we expect that $\widehat{\#}(A^{(k)})\propto \kappa\cdot \log(1/\epsilon^{(k)}) \propto \alpha\cdot \log(1/\epsilon^{(k)})$. This linear dependence is confirmed in \Cref{fig:maxA}. It is also important to note that the curve(green) corresponding to level 1 is below the curve (blue) corresponding to level 2 in \Cref{fig:maxA_normalize_vs_alpha}, which again implies that $\widehat{\#}(A^{(k)})/\log(1/\epsilon^{(k)})$ is uniformly bounded for all levels.

\section{Comparison With The Implicit Restarted Lanczos Method (IRLM)}
\label{sec:compare_IRLM}

Owning to the observation in \cite{martinez2016tuned} that Implicit Restarted Lanczos Method (IRLM) is still one of the most performing and well-known algorithms for finding a large portion of smallest eigenpairs, in this section, we compare the computation complexity of our proposed algorithm with the IRLM.

To quantitatively compare the two methods, we record the computation time and the number of Conjugate gradient iterations as the benchmarks. The reasons for doing this are as follows: 
\begin{itemize}
\item In large-scale setting, direct method for solving sparse matrix $A^{-1}$ is general, not practical since large memory storage is required. Instead, iterative methods, especially the Conjugate gradient method (as $A$ is SPD in our case) is employed.
\item In both the IRLM and our proposed algorithm, the dominating complexity comes from the operator of solving for $A^{-1}$.
%\item Comparing the iteration number is also fairer than comparing the computation time directly as there are various public programming implementation for the IRLM, which are optimized in different scenarios. We also have to admit that the recent implementation of our proposed algorithm is not optimal. But we would like to emphasize that any (programming) techniques for optimizing the IRLM is also suitable for our implementation since the core computation are the conjugate gradient steps appear in the Hierarchical Eigenpair Computation Algorithm. 
\end{itemize}

\begin{remark}
For small-scale problems, a direct solver (such as sparse Cholesky factorization) for $A^{-1}$ is preferred in the IRLM. In this way, only one factorization step for $A$ is required prior to the IRLM. Moreover, solving for $A^{-1}$ in each iteration is replaced by solving two lower triangular matrix systems. This will bring a significant speedup for the IRLM. However, recall that we are aiming at understanding the asymptotic behavior and performance of these methods. Therefore, the IRLM discussed in this section employs the iterative solver instead of a direct solver.
\end{remark}

To be consistent, all the experiments are performed on a single machine equipped with Intel(R) Core(TM) i5-4460 CPU with 3.2GHz and 8GB DDR3 1600MHz RAM. Both the proposed algorithm and the IRLM are implemented using C++ with the Eigen Library for fairness. In particular, the built-in (Preconditioned) conjugate gradient solvers are used in the IRLM implementation, instead of implementing on our own.

\begin{table}[h!]
\centering
%\scriptsize %\tiny % \footnotesize, \small, \normalsize
\footnotesize
\renewcommand{\arraystretch}{1.2}
    \begin{tabular}{|c|C{1.5cm}|C{1.5cm}|C{2.5cm}|C{2.5cm}|C{2.5cm}|}
    \hline
\# Eigenpairs & \multicolumn{2}{c|}{Methods} & 4-level Brain & 4-level SwissRoll & 3-level SwissRoll \\ \hline
& \multicolumn{2}{c|}{Decomposition} & \cellcolor{Gray2} 34.589 & \cellcolor{Gray2} 8.124 & \cellcolor{Gray2} 9.430 \\ \hline
\multirow{7}{*}{300} 
& \multirow{5}{*}{Proposed} 
& Level-4 	& 0.010 & 0.011 & - \\ \cline{3-6} 
& & Level-3 & 0.841 & 0.560 & 0.083 \\ \cline{3-6} 
& & Level-2 & 29.122 & 40.796 & 18.729 \\ \cline{3-6} 
& & Level-1 & 61.286 & 18.846 & 22.440 \\ \cline{3-6} 
& & Total   & \cellcolor{Gray} 125.848 & \cellcolor{Gray} 68.337 & \cellcolor{Gray} 50.682 \\ \cline{2-6} 
& \multicolumn{2}{c|}{IRLM-CG} & 174.028 & \multicolumn{2}{c|}{81.005} \\ \cline{2-6} 
& \multicolumn{2}{c|}{IRLM-ICCG} & 525.73 & \multicolumn{2}{c|}{289.385} \\ \hline
    
\multirow{7}{*}{200} 
& \multirow{5}{*}{Proposed} 
& Level-4 	& 0.010 & 0.011 & - \\ \cline{3-6} 
& & Level-3 & 0.826 & 0.526 & 0.083 \\ \cline{3-6} 
& & Level-2 & 25.560 & 28.094 & 11.517 \\ \cline{3-6} 
& & Level-1 & 54.951 & 12.107 &  18.378 \\ \cline{3-6} 
& & Total   & \cellcolor{Gray} 115.936 & \cellcolor{Gray} 48.862 & \cellcolor{Gray} 39.408 \\ \cline{2-6} 
& \multicolumn{2}{c|}{IRLM-CG} & 124.871 & \multicolumn{2}{c|}{61.479} \\ \cline{2-6} 
& \multicolumn{2}{c|}{IRLM-ICCG} & 417.632 & \multicolumn{2}{c|}{196.217} \\ \hline

\multirow{7}{*}{100} 
& \multirow{5}{*}{Proposed} 
& Level-4 	& 0.010 & 0.011 & - \\ \cline{3-6} 
& & Level-3 & 0.831 & 0.531 & 0.083 \\ \cline{3-6} 
& & Level-2 & 25.056 & 22.062 & 9.883 \\ \cline{3-6} 
& & Level-1 & 31.882 & 8.066 & 12.029 \\ \cline{3-6} 
& & Total   & \cellcolor{Gray} 92.368 & \cellcolor{Gray} 38.794 & \cellcolor{Gray} 31.425 \\ \cline{2-6} 
& \multicolumn{2}{c|}{IRLM-CG} & 115.676 & \multicolumn{2}{c|}{48.713} \\ \cline{2-6} 
& \multicolumn{2}{c|}{IRLM-ICCG} & 324.648 & \multicolumn{2}{c|}{90.175} \\ \hline
    \end{tabular}
    \caption{Computation time (in seconds) for the 4-level Brain, 3-level SwissRoll and the 4-level SwissRoll examples using the proposed Hierarchical multi-level eigensolver; the IRLM with Conjugate Gradient solver and the IRLM with incomplete Cholesky preconditioned Conjugate Gradient solver.\label{table:computationtime}}
\end{table}

Table \ref{table:computationtime} shows the overall computation time for computing the leftmost (i) 300; (ii) 200 and (iii) 100 eigenpairs using (i) our proposed algorithm, (ii) the IRLM with incomplete Cholesky preconditioned Conjugate Gradient (IRLM-ICCG); and (iii) the IRLM with classical conjugate gradient method (IRLM-CG). In this numerical example, the error tolerance of the eigenvalues in all three cases are set to $10^{-5}$. Since the error for IRLM cannot be obtained a priori, we fine-tune the relative error tolerance for the (preconditioned) conjugate gradient solver such that eigenvalues error are of order $O(10^{-6})$. For the proposed algorithm, the time required for level-wise eigenpair computation is recorded. In the bottom level (level-4 or level-3 in these cases), we have used the built-in eigensolver function in the Eigen Library to obtain the full eigenpairs (corresponding to Line 1 in Algorithm 6). As the problem size is small, the time complexity is insignificant for all three examples. 

The total runtime of our proposed algorithm in each example is computed by summing up all levels' computation time, plus the operator decomposition time (which is the second row in \Cref{table:computationtime}). For all these examples, our proposed algorithm outperforms the IRLM. Although both the size of the matrices and their corresponding condition numbers are not extremely large, the numerical experiments already show a observable improvement. From the theoretical analysis discussed in the previous sections, this improvement will even be magnified if the SPD matrices are of larger scales and more ill conditioned. Indeed, we assert that our proposed algorithm cannot be fully utilized in these illustrations. Therefore, one of the main future works is to perform detailed numerical experiments in these cases. For instance, by considering the 3-level and 4-level SwissRoll examples, we observe that a 3-level decomposition is indeed sufficient for SwissRoll graph laplacian, where we recall the corresponding condition number is $\| A \|_2 = 1.15 \times 10^6$. Therefore, using a 3-level decomposition, the overall runtime reduction goes up to approximately 37\% if 300 eigenpairs are required.

\begin{figure}[!h]
\centering
    \begin{subfigure}[b]{0.43\textwidth}
        \includegraphics[width=\textwidth]{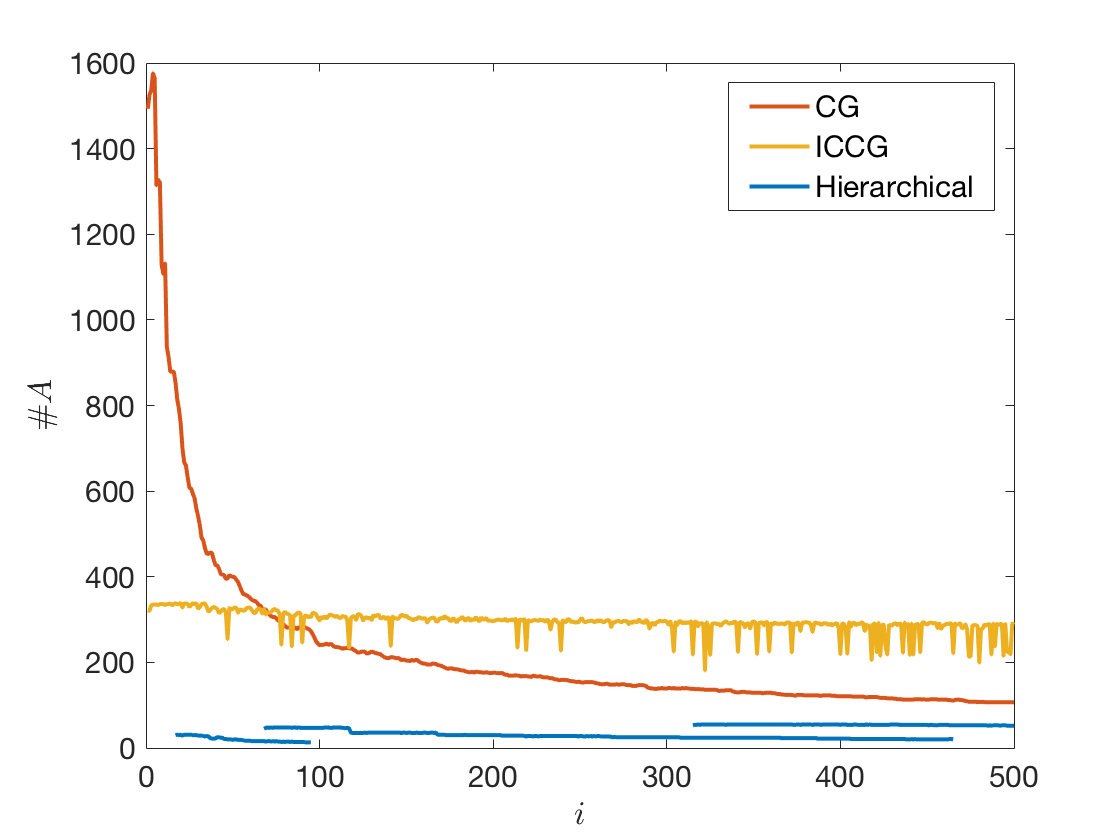}
        \caption{}
        \label{fig:CG_iter_number}
    \end{subfigure}
    \begin{subfigure}[b]{0.43\textwidth}
        \includegraphics[width=\textwidth]{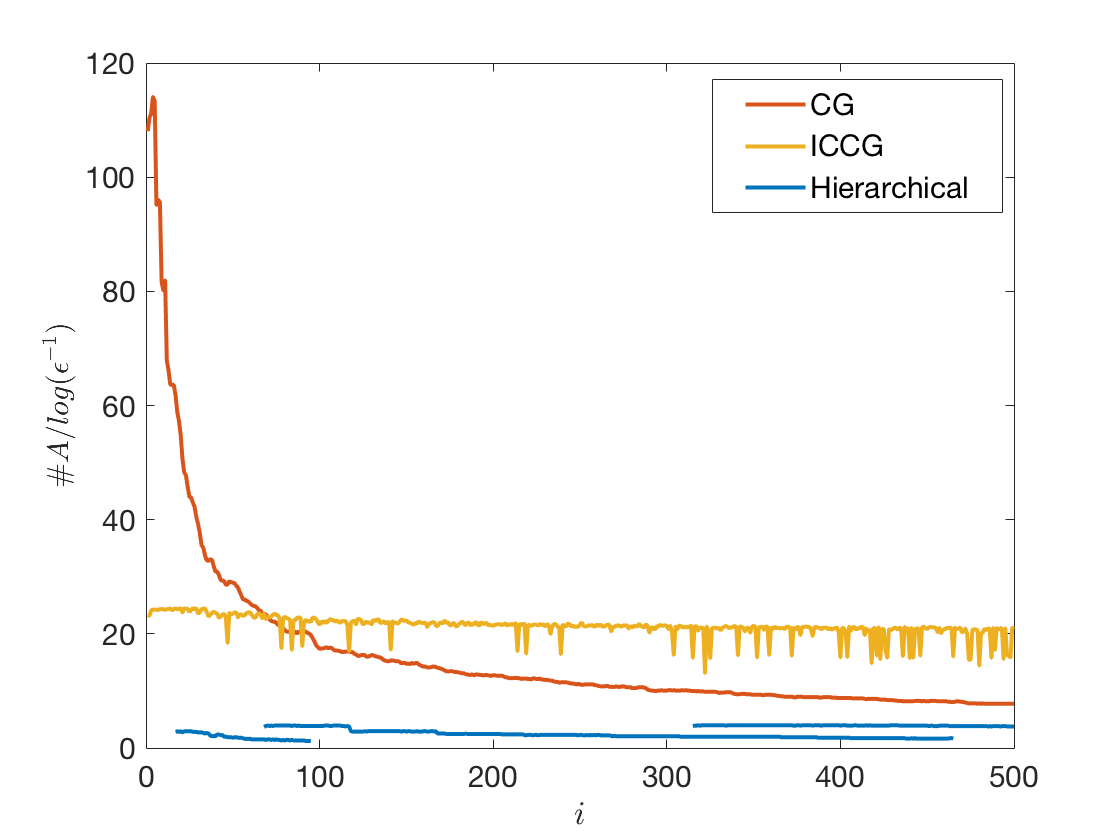}
        \caption{}
        \label{fig:log-CG_iter_number}
    \end{subfigure}
    \caption{(a) The PCG iteration number in the 4-level SwissRoll example. The IRLM-ICCG methods exhibits a uniform iteration number, while the IRLM-ID has an exponential decaying iteration number. For our proposed algorithm, since the spectrum-preserving hierarchical preconditioner $M$ is employed, the CG iteration number is minimum. This is also consistent to the time complexity shown in \Cref{table:computationtime}. (b) The corresponding normalized plot, where the iteration number is normalized by $\log (\epsilon)$.}
\label{fig:CG_iter-number}
\end{figure}

Notice that the time required for the IRLM-ICCG is notably much more than that of the IRLM-CG, which contradicts to our usual experience regarding preconditioning. In fact, such phenomenon can be explained as follows: In the early stage of the IRLM, preconditioning with incomplete Cholesky factorization helps reducing the iteration number of the CG. However, when the eigen-subspace are gradually projected away throughout the IRLM process, the spectrum of the remaining subspace reduces and therefore CG iteration numbers also drops significantly. On the contrary, preconditioning with incomplete Cholesky ignores such update in spectrum and therefore the CG iteration number is uniform throughout the whole Lanczos iteration. Hence, the classical CG method is preferred if a large number of leftmost eigenpairs are required. \Cref{fig:CG_iter_number} shows the CG iteration numbers in the IRLM-ICCG, IRLM-CG and respectively, our proposed hierarchical eigensolver versus the Lanczos iteration. More precisely, if we call $V_k$ in \cref{eqt:lanczos_vector} to be the {\it Lanczos vector}, the x-axis in the figure then corresponds to the first time we generate the $i$-th column of the Lanczos vector. For IRLM methods, it is equivalent to the extension procedure for the $i$-th column of the Lanczos vector, which corresponds to Line 6 -- 8 in Algorithm 1. In particular, the CG iteration number recorded in this figure corresponds to the operation $op$ in Line 7 of Algorithm 1. For our proposed algorithm, there are three separate sections, each section's CG iteration numbers correspond to the formation of Lanczos vectors in the $3^{\text{rd}}$-, $2^{\text{nd}}$- and the $1^{\text{st}}$-level respectively. Since we may also update some of these Lanczos vector during the refinement process, therefore some overlaps in the recording of CG iteration numbers corresponding to those Lanczos vector are observed. With the spectrum-preserving hierarchical preconditioner $M$ introduced in our algorithm, the CG iteration number for solving $A^{-1}$ is tremendously reduced. In contrast, the CG iteration number for IRLM-CG is the largest at the beginning but decreases exponentially and asymptotically converges to our proposed result. For IRLM-ICCG, the incomplete Cholesky factorization does not capture the spectrum update and therefore the iteration numbers is uniform throughout the computation. This observation is also consistent to the time complexity as shown in \Cref{table:computationtime}. \Cref{fig:log-CG_iter_number} shows the corresponding normalized plot, where the iteration number is normalized by $\log (\frac{1}{\epsilon})$.

Similar results can also be plotted for the 4-level Brain and the 3-level SwissRoll examples. We therefore skip those plots to avoid repetition.

\section{Conclusion And Future Works}
\label{sec:conclusion}
In this work, we propose a spectrum preserving preconditioned hierarchical eigensolver to compute a large number of leftmost eigenpairs of a sparse symmetric positive definite matrix. This eigensolver exploits the well-conditioned property of the decomposition components obtained through the MMD, the nice spectral property Lanczos procedure and also the preconditioning characteristics of the CG method. In particular, we propose an extension-refinement iterative scheme, in which eigenpairs are hierarchically extended and refined from the ones obtained from the previous level up to the desired amount. A specially designed spectrum-preserving preconditioner is also introduced for the PCG method to solve for $A^{-1}$ during the iterations. Theoretical analysis on the runtime complexity and the asymptotic behavior of our proposed algorithm are reported. Quantitative numerical experiments and comparison with the IRLM are also reported to demonstrate the efficiency and effectiveness of our proposed algorithm.

We would like to remark that the proposed algorithm and its implementation are still in the early stage as the main purpose of this work is to explore the possibility of integrating the multiresolution operator compression framework with the Krylov-type iterative eigensolver. Therefore, one of the future topics is to conduct a comprehensive numerical studies of our algorithm to various large-scale, real data such as graph Laplacians of real network data, or stiffness matrices stemmed from the discretization of high-contrasted elliptic PDEs. These studies will help numerically confirm the asymptotic behavior of the relative condition numbers of $M$ and $A_{st}$, especially when we need to compute a large number of leftmost eigenpairs from large-scale operators. Another possible research direction is to investigate the parallelization of this algorithm. This is important when we solve a large scale eigenvalue problem.

\section*{Acknowledgment}
The research was in part supported by the NSF Grants DMS‐1318377 and DMS-1613861. Ziyun Zhang would like to acknowledge ACM, Caltech and SMS, PKU for supporting her research visit to Caltech in 2017 summer. She would also like to thank ACM's staff for their hospitality during her visit.

\appendix
\section*{Appendix A}
In this section, we compare the our method for compressed eigenproblem and the method proposed by Ozoliņ\v{s} et al. \cite{ozolicnvs2013compressed}. We start with the straightforward compression directly using the eigenvectors corresponding to smallest eigenvalues, which can be obtained by solving the following optimization problem:
\begin{equation}
\begin{array}{ll}
 \Psi\ =& \underset{\widehat{\Psi}}{\argmin}\ \sum_{i=1}^N \hat{\psi}_i^TA\hat{\psi}_i,\\ 
 & \text{s.t.}\quad \hat{\psi}_i^T\hat{\psi}_j=\delta_{ij},\ i,j=1,2,\cdots,N.
\end{array}
\label{eqt:PCA_op}
\end{equation}
The compression using eigenvectors is well known as the PCA method is optimal in 2-norm sense for fixed compressed dimension $N$. However, computing a large number of eigenvectors is a hard problem itself, not to mention that we actually intend to approximate eigenpairs using the compressed operator. Also the spatially extended profiles of exact eigenvectors make them less favorable in many fields of researches. Then as modification, Ozoli\v{s} et al. \cite{ozolicnvs2013compressed} added a $L_1$ regularization term to impose the desired locality on $\Psi$. They modified the optimization problem \cref{eqt:PCA_op} as
\begin{equation}
\begin{array}{ll}
 \Psi\ =& \underset{\widehat{\Psi}}{\argmin}\ \sum_{i=1}^N \Big(\hat{\psi}_i^TA\hat{\psi}_i+\frac{1}{\mu}\|\hat\psi_i\|_1\Big),\\ 
 & \text{s.t.}\quad \hat{\psi}_i^T\hat{\psi}_j=\delta_{ij},\ i,j=1,2,\cdots,N.
\end{array}
\label{eqt:Osher_op}
\end{equation}
The $L_1$ regularization, as widely used in many optimization problems for sparsity pursuit, effectively ensures each output $\psi_i$ to have spatially compact support, at the cost of compromising the approximation accuracy compared to PCA. The factor $\mu$ controls the locality of $\Psi$. A smaller $\mu$ gives more localized profiles of $\Psi$, which, however, results in larger compression error for a fixed $N$. The loss of approximation accuracy can be compensated by increasing, yet not significantly, the basis number $N$. An algorithm based on the split Bregman iteration was also proposed in \cite{ozolicnvs2013compressed} to effectively solve the problem \cref{eqt:Osher_op}. In summary, their work provides an effective method to find a bunch of localized basis functions that can approximately span the eigenspace of smallest eigenvalues of $A$.

Although our approach to operator compression is originally developed from a different perspective based on Finite Element Method (FEM), it can be reformulated as an optimization problem similar to \Cref{eqt:PCA_op}. In fact, to obtain the basis $\Psi$ used in our method, we can simply replace the nonlinear constraints $\psi_i^T\psi_j=\delta_{ij},\ i,j=1,2,\cdots,N,$ by linear constraints $\psi_i^T\phi_j=\delta_{ij},\ i,j=1,2,\cdots,N,$ to get 
\begin{equation}
\begin{array}{ll}
 \Psi\ =& \underset{\widehat{\Psi}}{\argmin}\ \sum_{i=1}^N \hat{\psi}_i^TA\hat{\psi}_i,\\ 
 & \text{s.t.}\quad \hat{\psi}_i^T\phi_j=\delta_{ij},\ i,j=1,2,\cdots,N.
\end{array}
\label{eqt:Ourmethod_op}
\end{equation}
Here $\Phi=[\phi_1,\phi_2,\cdots,\phi_N]$ is a dual basis that we construct ahead of $\Psi$ to provide a priori compression error estimate as stated in \Cref{thm:compression_thm1}. As the constraints become linear, problem \cref{eqt:Ourmethod_op} can be solved explicitly by $\Psi=A^{-1}\Phi(\Phi^TA^{-1}\Phi)^{-1}$ as mentioned in \cref{eqt:psi_closedform}. Instead of imposing locality by adding $L_1$ regularization as in \cref{eqt:Osher_op}, we obtain the exponential decaying feature of $\Psi$ by constructing each dual basis function $\phi_i$ locally. That is the locality of $\Phi$ and the strong correlation $\Psi^T\Phi=I$ automatically give us the locality of $\Psi$ under energy minimizing property. The optimization form \cref{eqt:Ourmethod_op} was derived by Owhadi in \cite{owhadi2017multigrid} where $\Psi$ was used as the FEM basis to solve second-order elliptic equations with rough coefficients. This methodology was then generalized to problems on higher order elliptic equations \cite{hou2016sparse}, general Banach space \cite{owhadi2017universal} and general sparse SPD matrix \cite{hou2017adaptive}. In all previous works the nice spectral property of $\Psi$ has been observed and in particular the eigenspace corresponding to the smallest $M$ eigenvalues of $A$ can be well approximately spanned by $\Psi$ of a relative larger dimension $N=O(M)$.

To further compare the problems \cref{eqt:Osher_op} and \cref{eqt:Ourmethod_op}, we test both of them on the one-dimensional Kronig--Penney (KP) model studied in \cite{ozolicnvs2013compressed} with rectangular
potential wells replaced by inverted Gaussian potentials. In this example, the matrix $A$ comes from discretization of the PDE operator $-\frac{1}{2}\Delta+V(x)$ defined on the domain $\Omega$ with periodic boundary condition. In particular, $\Omega=[0,50]$, and $V(x)=-V_0\sum_{j=1}^{N_{el}}\exp\big(-\frac{(x-x_j)^2}{2\delta^2}\big)$. As in \cite{ozolicnvs2013compressed}, we discretize $\Omega$ with 512 equally spaced nodes, and we choose $N_{el}=5$, $V_0=1$, $\delta=3$, and $x_j=10j-5$(instead of $x_j=10j$ in \cite{ozolicnvs2013compressed}, which essentially changes nothing).

For problem \cref{eqt:Ourmethod_op}, we divide $\Omega$ into $N$ equal-length intervals $\{\Omega_i\}_{i=1}^N$, and choose the dual basis $\Phi=[\phi_1,\phi_2,\cdots,\phi_N]$ such that $\phi_i$ is the discretization of the indicator function $\mathbf{1}(\Omega_i)$($\mathbf{1}(\Omega_i)(x)=1$ for $x\in\Omega_i$, otherwise $\mathbf{1}(\Omega_i)(x)=0$). We use $\Psi_{o}$ to denote the exact result of problem \cref{eqt:Ourmethod_op}, namely $\Psi_{o}=A^{-1}\Phi(\Phi^TA^{-1}\Phi)^{-1}$. Since $\Psi_{o}$ is not orthogonal, we should compute the eigenvalues from the general eigenvalue problem $\Psi_o^TA\Psi_o v=\lambda\Psi_o^T\Psi_ov$(\Cref{lemma:generaleigenproblem}) as approximations of the eigenvalues of $A$. We use $\lambda_o$ to denote these approximate eigenvalues.

For problem \cref{eqt:Osher_op}, we use Algorithm 1 and exactly the same parameters provided in \cite{ozolicnvs2013compressed}, which means we are simply reproducing their results, except that we use a finer discretization (512 rather than 128) and we shift the potential $V(x)$. We have used normalized $\Phi$ as the initial guess for Algorithm 1 in \cite{ozolicnvs2013compressed}, and choose $\mu=10$. We use $\Psi_{cm}$ to denote the result of problem \cref{eqt:Osher_op}. We use $\lambda_{cm}$ to denote the eigenvalues of $\Psi_{cm}^TA\Psi_{cm}$. 

\begin{figure}[!h]
\centering
		$N=50$\qquad\qquad\qquad\qquad\quad $N=75$\qquad\qquad\qquad\qquad\quad $N=100$\\
        \includegraphics[width=0.25\textwidth]{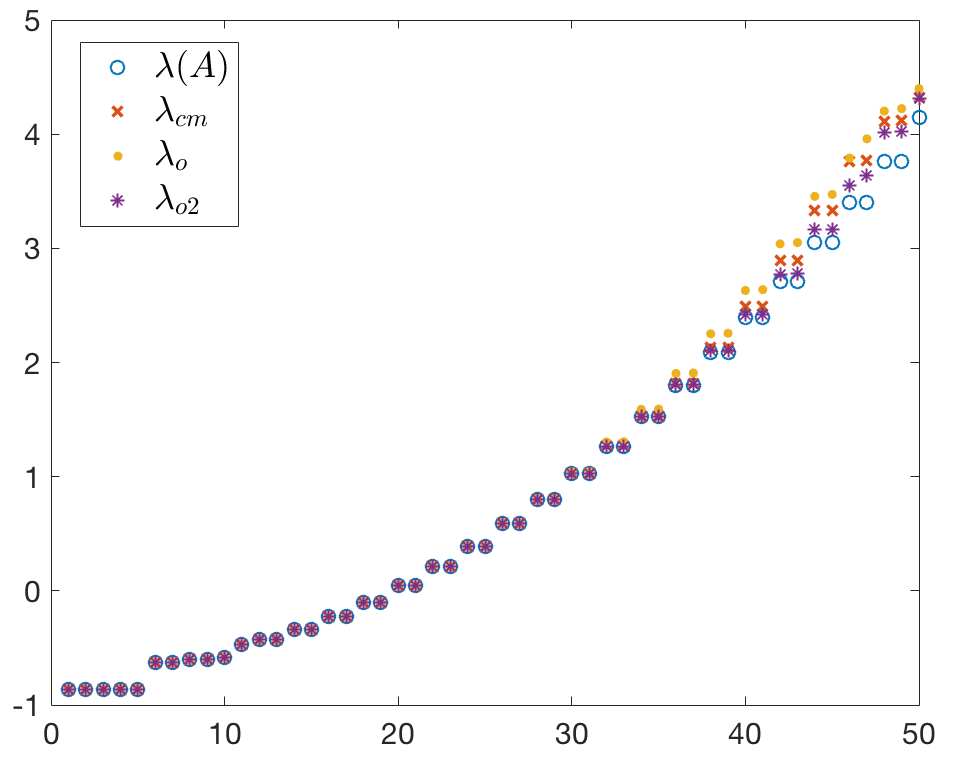}
        \includegraphics[width=0.25\textwidth]{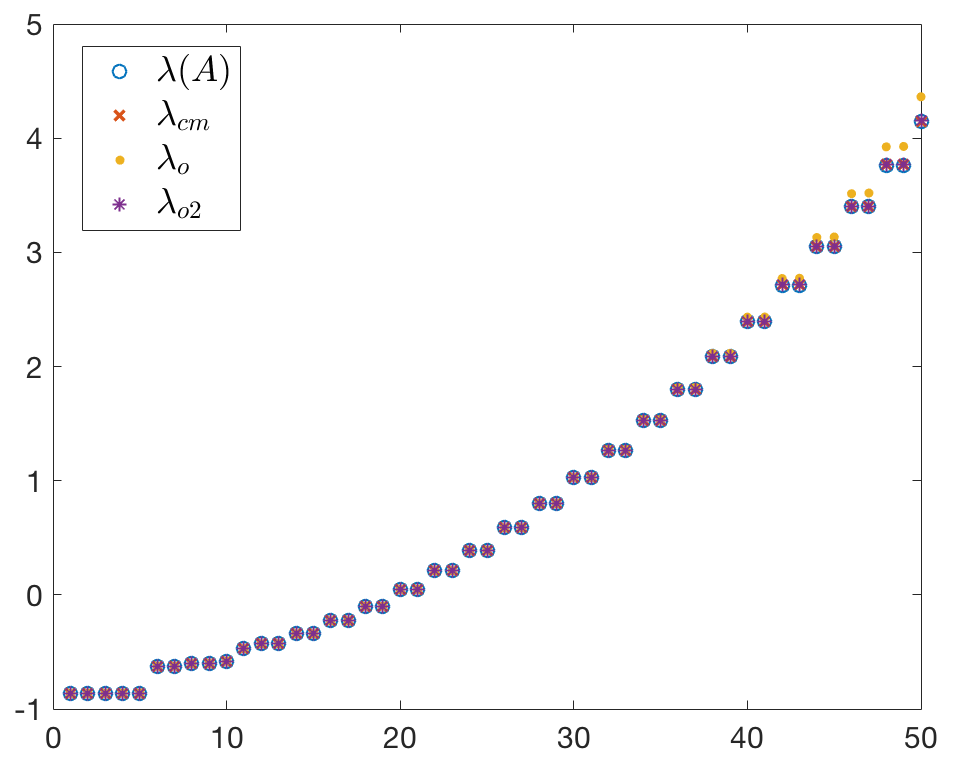}      
        \includegraphics[width=0.25\textwidth]{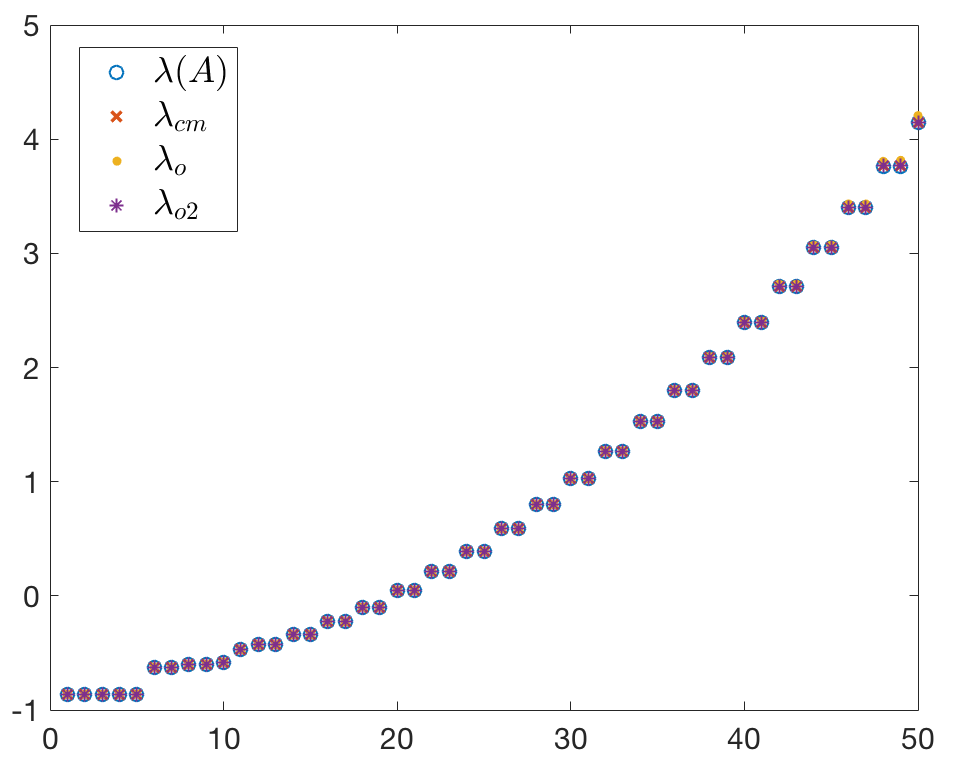}
        \includegraphics[width=0.25\textwidth]{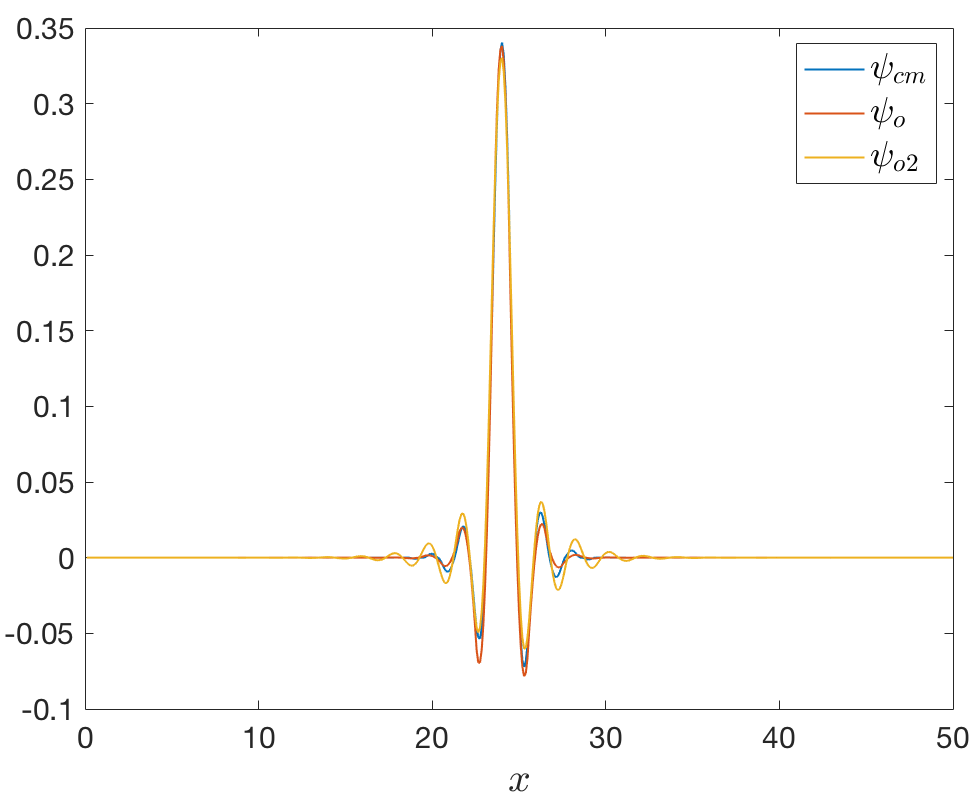}
        \includegraphics[width=0.25\textwidth]{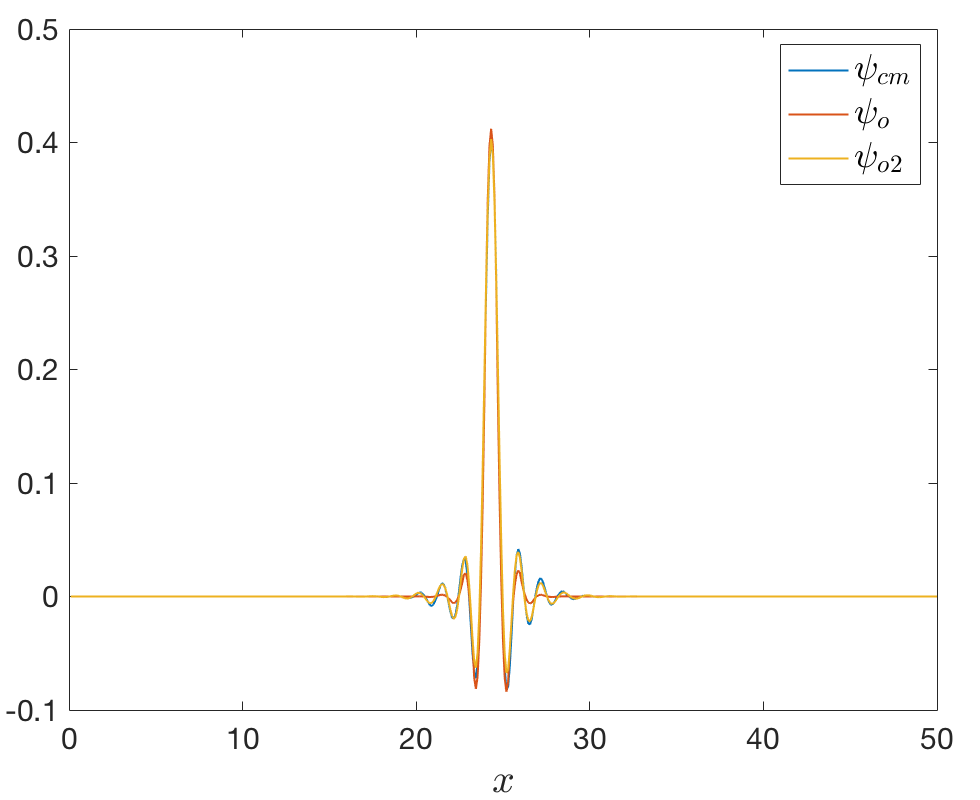}      
        \includegraphics[width=0.25\textwidth]{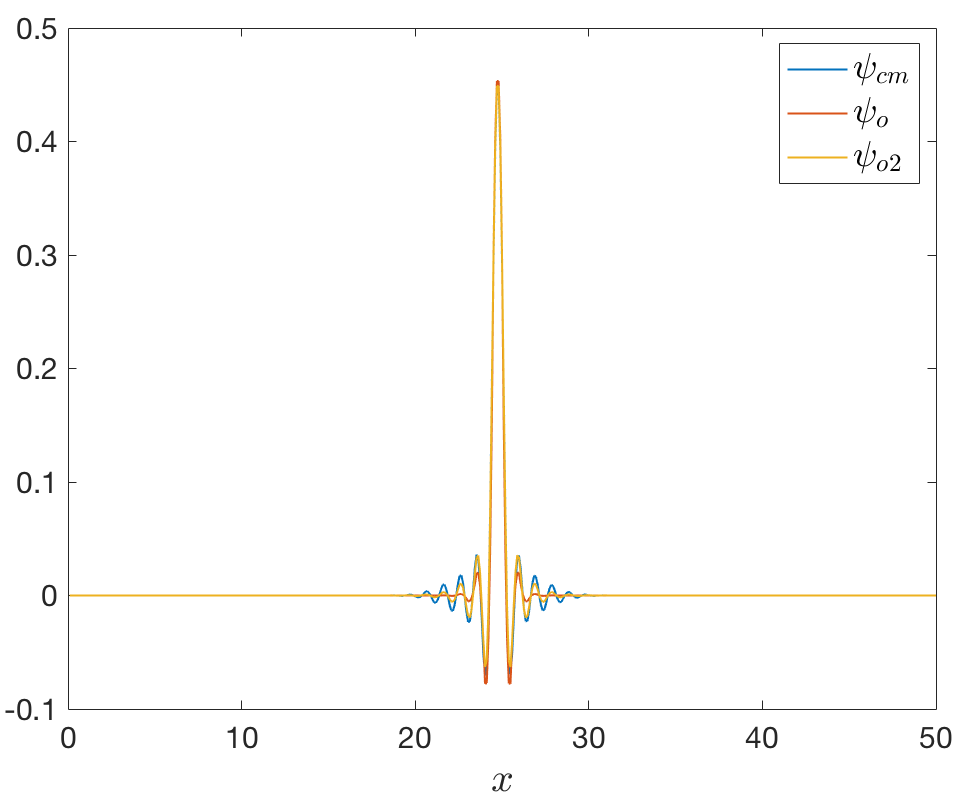}
        \includegraphics[width=0.25\textwidth]{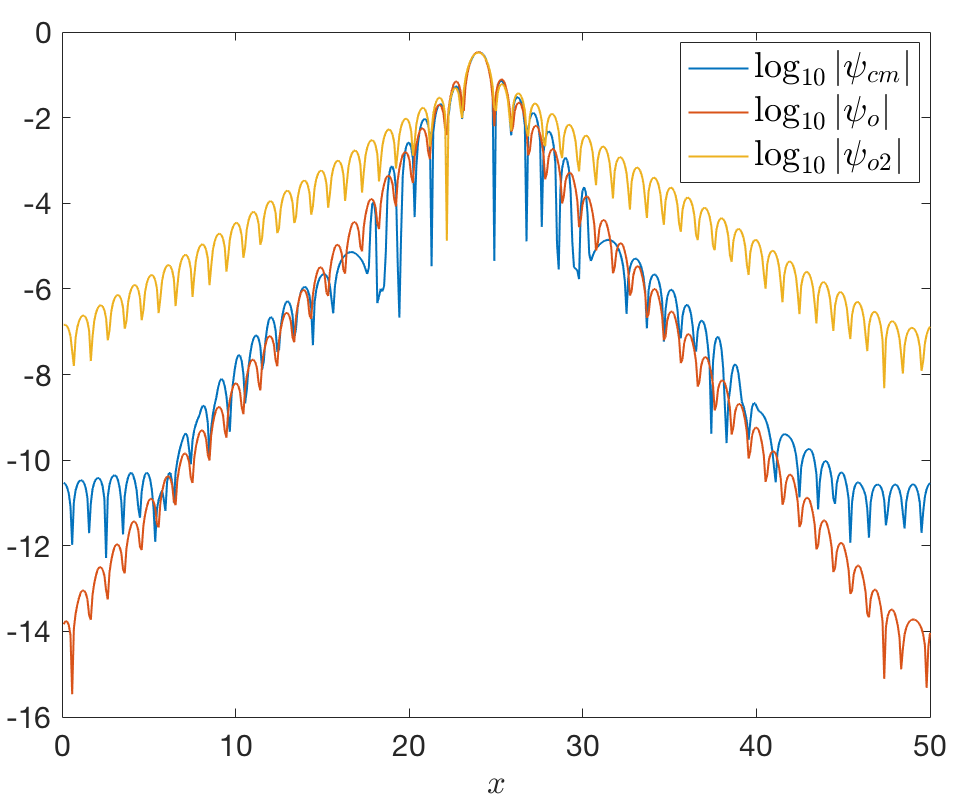}
        \includegraphics[width=0.25\textwidth]{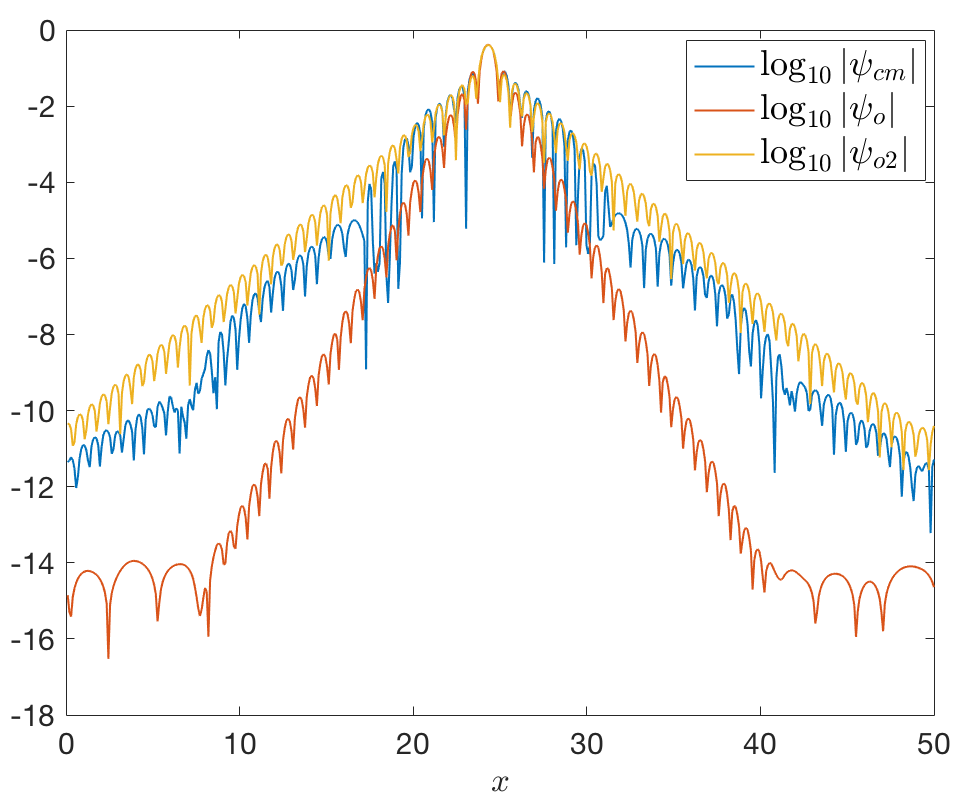}      
        \includegraphics[width=0.25\textwidth]{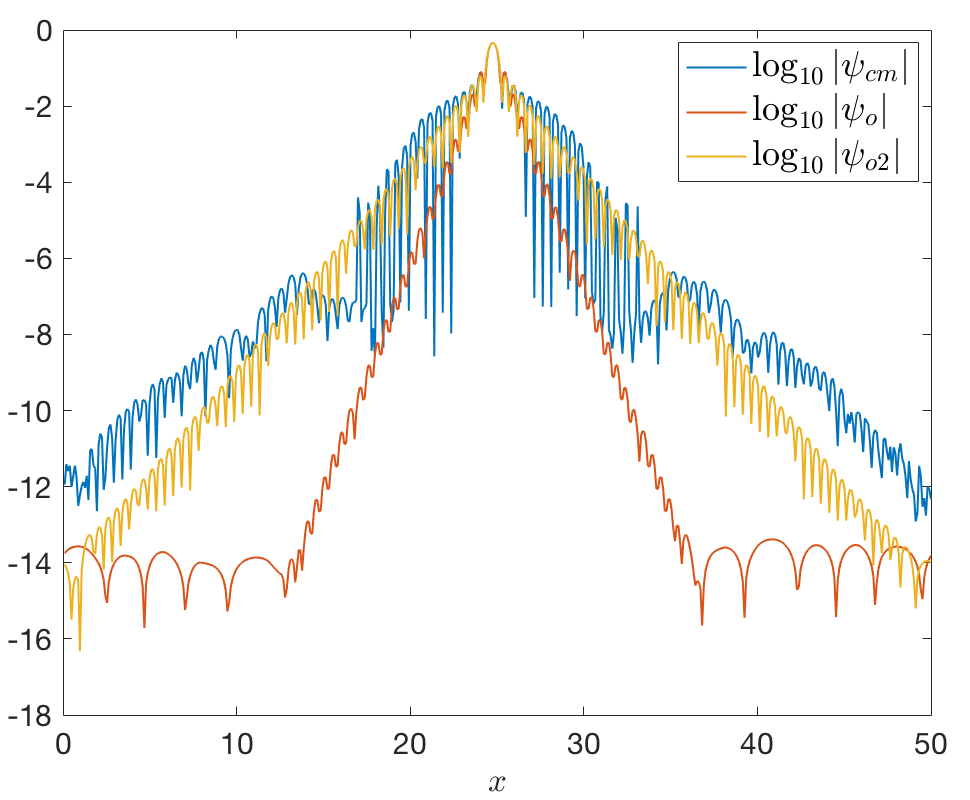}
\caption{Results of problems \cref{eqt:Osher_op} and \cref{eqt:Ourmethod_op} for $N=50$(first column), $N=75$(second column) and $N=100$(third column). First row: the first 50 eigenvalues of $A$ and those of the compressed problems. Second row: examples of local basis functions. Third row: examples of local basis functions in log scale. } 
\label{fig:compare_osher}     
\end{figure}

We compare the approximate eigenvalues to the first 50 eigenvalues of $A$. The first row of \Cref{fig:compare_osher} shows that both methods give very good approximations of $\lambda(A)$. And when $N$ increases, the approximations become better. But relatively, the results $\lambda_{cm}$ from \Cref{fig:compare_osher} is closer to the ground truth than our results $\lambda_o$ from \cref{eqt:Ourmethod_op}. To improve our results, we simply solve problem \cref{eqt:Ourmethod_op} again, but this time using previous result $\Psi_o$ as the dual basis. That is we compute $\Psi_{o2}=A^{-1}\Psi_o(\Psi_o^TA^{-1}\Psi_o)$, and compute eigenvalues $\lambda_{o2}$ from the general eigenvalue problem $\Psi_{o2}^TA\Psi_{o2} v=\lambda\Psi_{o2}^T\Psi_{o2}v$. We can see that the approximate eigenvalues $\lambda_{o2}$ are even closer to the ground truth. An interpretation of this improvement is that if we see $\Psi_o=A^{-1}\Phi(\Phi^TA^{-1}\Phi)^{-1}$ as a transformation from $\Phi$ to $\Psi_o$, then the part $A^{-1}\Phi$ is equivalent to applying inverse power method to make $\Psi_o$ more aligned to the eigenspace of the smallest eigenvalues, while the part $(\Phi^TA^{-1}\Phi)^{-1}$ is to force $\Psi_o^T\Phi=I$ so $\Psi_o$ inherits some weakened locality from $\Phi$. So if we apply this transformation to $\Psi_o$ again to obtain $\Psi_{o2}$, $\Psi_{o2}$ will approximate the eigenspace of the smallest eigenvalues better, but with more loss of locality.

In the second row and third row of \Cref{fig:compare_osher}, we show some examples of the local basis functions $\psi_{cm}$, $\psi_o$ and $\psi_{o2}$ (all are normalized to have unit $l_2$ norm). Interestingly, these basis functions are not just localized as expected, but indeed they have very similar profiles. One can see that for $N=75$, the basis functions $\psi_{cm}$ and $\psi_o$ are almost identical. So it seems that in spite of how we impose locality (either the $L_1$ minimization approach, or the construction of the dual basis $\Phi$), the local behaviors of the basis functions are determined by the operator $A$ itself. We believe that this ``coincidence'' is governed by some intrinsic property of $A$, which may be worth further exploring and studying. If we can understand a higher level, unified mechanism that results in the locality of the basis, we may be able to extend these methods to a more general class of operators. We also observed that as $N$ goes large, $\psi_o$ and $\psi_{o2}$ become more and more localized since the support of the dual basis functions are smaller and smaller. However the locality of $\psi_{cm}$ doesn't change much as $N$ increases, since we use the same penalty parameter $\mu=10$ for \cref{eqt:Osher_op} in this experiment. 

We would like to remark that, though these two problems result in local basis functions with similar profiles, problem \cref{eqt:Osher_op} requires to use the split Bregman iteration to obtain the $N$ basis functions simultaneously. In our problem \cref{eqt:Ourmethod_op}, since the constraints are linear and separable, the basis functions can be obtained separately and directly without iteration. Furthermore, thanks to the exponential decay of the basis functions, each subproblem for obtaining one basis function can be restricted to a local domain without significant loss of accuracy, and the resulting local problem can be solved very efficiently. For definitions and detailed properties of these local problems for obtaining localized basis, please refer to section 3 in \cite{hou2017adaptive}. Therefore the algorithm for solving problem \cref{eqt:Ourmethod_op} can be highly localized and embarrassingly parallel.

\section*{Appendix B} 
In this section, we qualitatively examine the accuracy of the approximate eigenvectors of the compressed operators by comparing their behaviors in image segmentation to those of the true eigenvectors of the original Laplacian operators. In the image segmentation, the eigenvectors of graph Laplacian provide a solution to graph partitioning problem. Namely, for a partition $(A,B)$ that satisfies $A \cup B = V$ and $A \cap B = \emptyset$, a measure of their disassociation called the normalized cut ($Ncut$) is defined as \cite{shi2000normalized}
% * <jefferykclam@gmail.com> 2017-11-28T00:19:30.527Z:
%
% > normalized cut
%
% ^.
\begin{equation}
Ncut(A,B) = \frac{cut(A,B)}{assoc(A,V)} + \frac{cut(A,B)}{assoc(B,V)} ,
\end{equation}
where
\begin{equation*}
cut(A,B) = \sum_{u \in A, v \in B} w(u,v), \qquad assoc(A,V) = \sum_{u \in A, t \in V} w(u,t).
\end{equation*}
Shi and Malik \cite{shi2000normalized} shows that, for a connected graph, minimizing $Ncut$ can be rephrased as finding the eigenvector $v_2$ that corresponds to the second smallest eigenvalue $\lambda_2$ of the graph Laplacian (since we always have $\lambda_1=0$ and $v_1$ a uniform vector). Taking $sign(v_2)$ transforms it into a binary vector which gives a satisfactory cut. Moreover, the next few eigenvectors provide further cuts of the previously partitioned fractions. Therefore, our eigensolver may serve as a powerful tool for graph partitioning, as well as its applications including image segmentation and manifold learning. 

\begin{figure}[!h]
\centering
\includegraphics[width=0.7\textwidth]{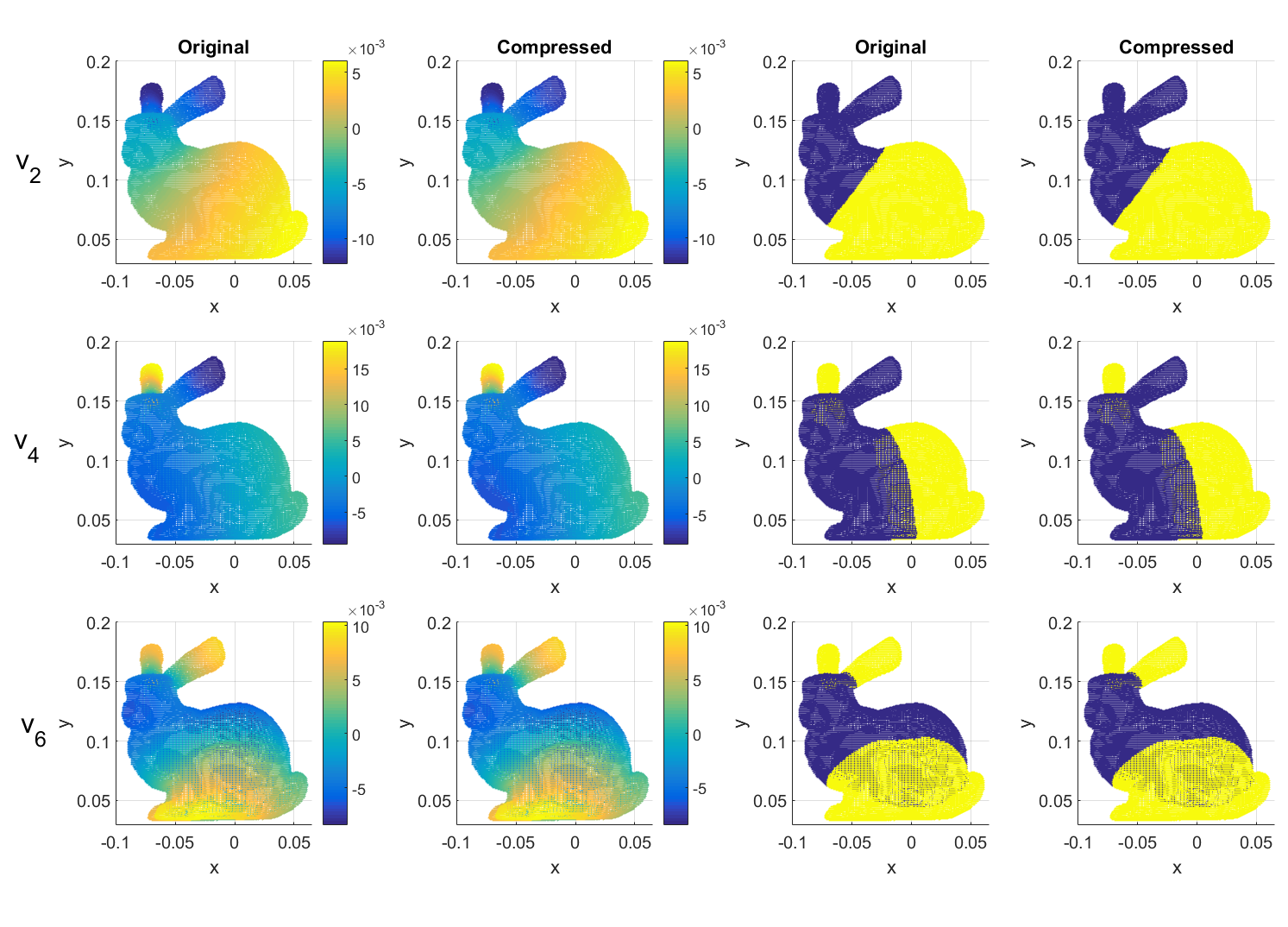} 
\caption{Colormap (left) and partition (right) using the $2^{\text{nd}}$, $4^{\text{th}}$ and $6^{\text{th}}$ eigenvectors of the original/compressed operator}
\label{fig:bunnycpr}
\end{figure}

\begin{figure}[!h]
\centering
\includegraphics[width=0.7\textwidth]{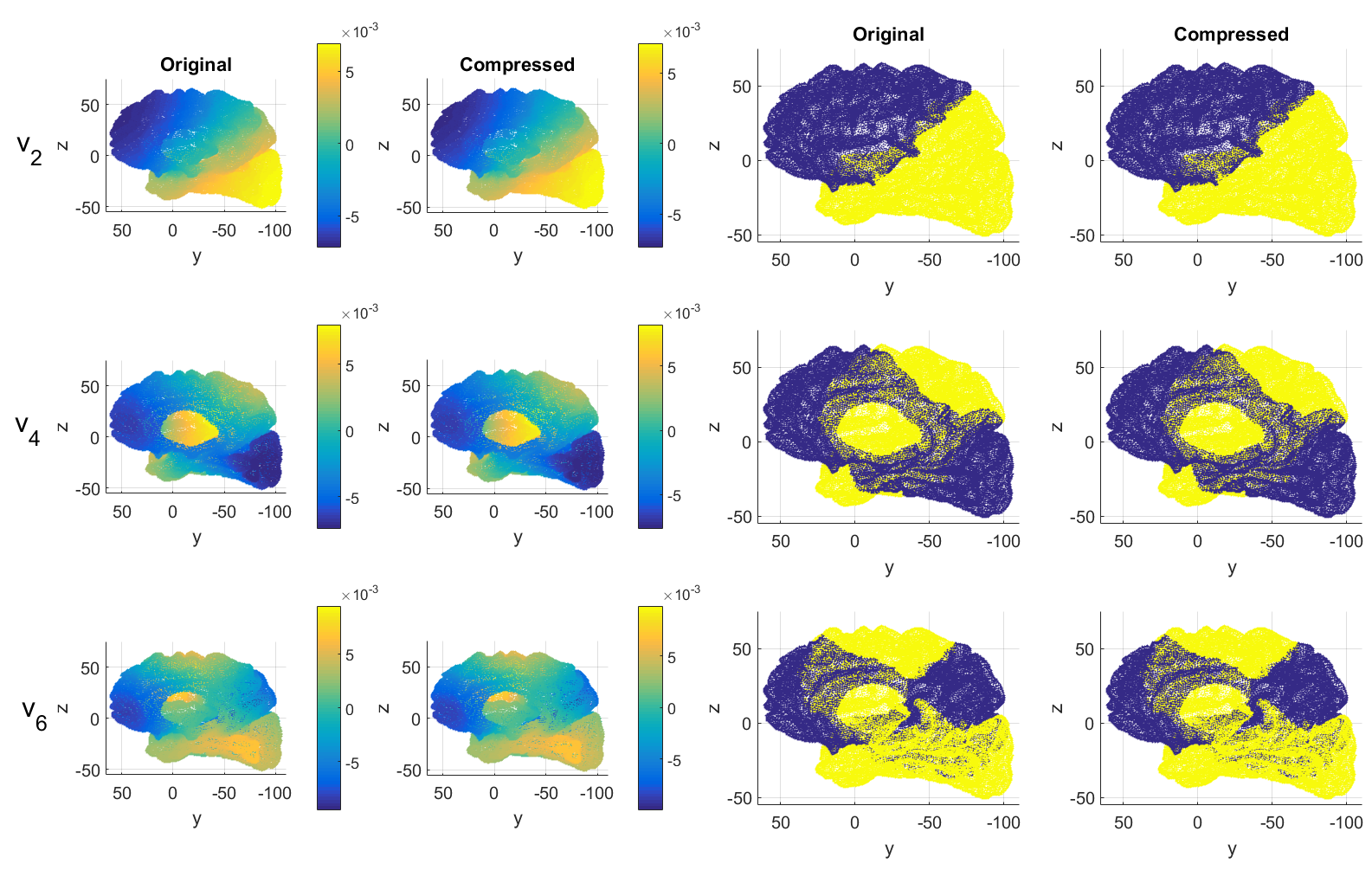} 
\caption{Colormap (left) and partition (right) using the $2^{\text{nd}}$, $4^{\text{th}}$ and $6^{\text{th}}$ eigenvectors of the original/compressed operator}
\label{fig:braincpr}
\end{figure}

We test graph partitioning on bunny and brain datasets using the eigenvectors of both original and compressed operators. Figures \ref{fig:bunnycpr} and \ref{fig:braincpr} shows the colormap and the partition generated by some selected eigenvectors. From the pictures we can see that the original and the compressed operators give very similar results when it comes to graph partitioning. The compressed operator is not only easier to compute, but also gives a satisfactory partition in practical settings.

\begin{figure}[!h]
\centering
\includegraphics[width=0.7\textwidth]{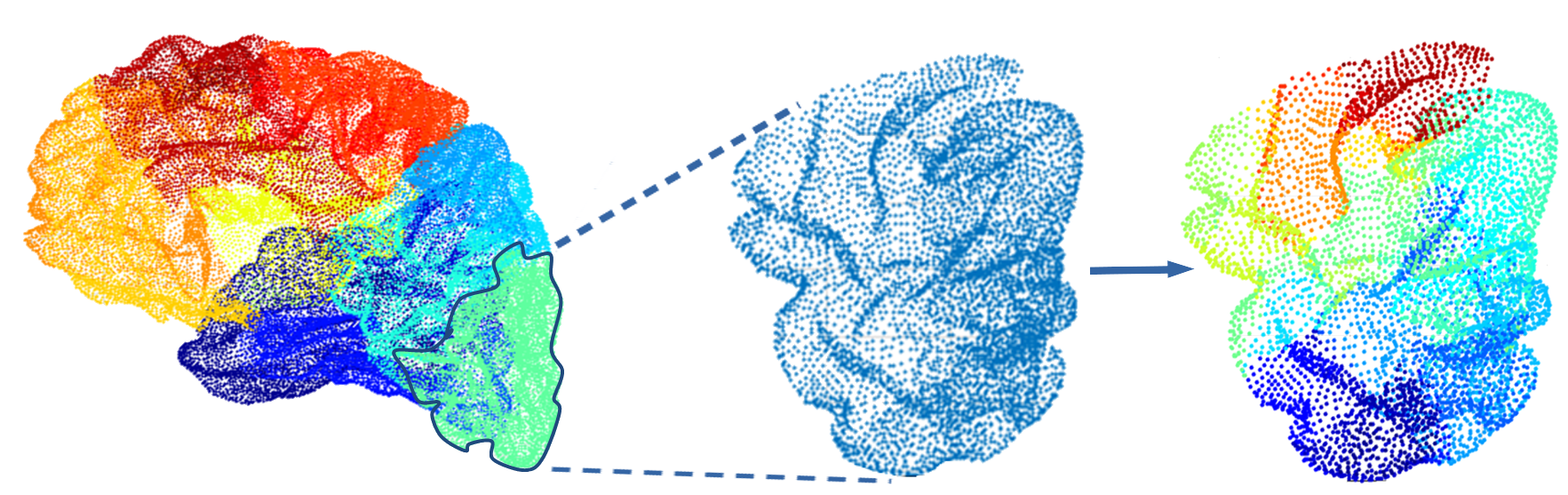} 
\caption{Heaping up more eigenvectors leads to finer partition. Left: partition using the first 5 eigenvectors. Middle: a uniform fraction from the previous partition. Right: further partition using the next 10 eigenvectors.}
\label{fig:multi}
\end{figure}

Figure \ref{fig:multi} gives an example of refining the partition with more eigenvectors. In the brain data, a fraction that is left intact in the first 5 eigenvectors (the light green part on the left) is divided into a lot more fractions when eigenvectors pile up to 15.

\bibliographystyle{siamplain}
\bibliography{reference}

\begin{thebibliography}{10}

\bibitem{bergamaschi2015computing}
{\sc L.~Bergamaschi and E.~Bozzo}, {\em Computing the smallest eigenpairs of
  the graph laplacian}, SeMA Journal,  (2015), pp.~1--16.

\bibitem{bergamaschi1997asymptotic}
{\sc L.~Bergamaschi, G.~Gambolati, and G.~Pini}, {\em Asymptotic convergence of
  conjugate gradient methods for the partial symmetric eigenproblem}, Numerical
  linear algebra with applications, 4 (1997), pp.~69--84.

\bibitem{bozzo2012effective}
{\sc E.~Bozzo and M.~Franceschet}, {\em Effective and efficient approximations
  of the generalized inverse of the graph laplacian matrix with an application
  to current-flow betweenness centrality}, arXiv preprint arXiv:1205.4894,
  (2012).

\bibitem{bozzo2013resistance}
{\sc E.~Bozzo and M.~Franceschet}, {\em Resistance distance, closeness, and
  betweenness}, Social Networks, 35 (2013), pp.~460--469.

\bibitem{calvetti1994implicitly}
{\sc D.~Calvetti, L.~Reichel, and D.~C. Sorensen}, {\em An implicitly restarted
  lanczos method for large symmetric eigenvalue problems}, Electronic
  Transactions on Numerical Analysis, 2 (1994), p.~21.

\bibitem{chung1997spectral}
{\sc F.~R. Chung}, {\em Spectral graph theory}, vol.~92, American Mathematical
  Soc., 1997.

\bibitem{cocco2013principal}
{\sc S.~Cocco, R.~Monasson, and M.~Weigt}, {\em From principal component to
  direct coupling analysis of coevolution in proteins: Low-eigenvalue modes are
  needed for structure prediction}, PLoS computational biology, 9 (2013),
  p.~e1003176.

\bibitem{francis1961}
{\sc J.~Francis}, {\em The transformation: a unitary analogue to the
  transformation. i}, Comput. J., 4 (1961), pp.~265--271.

\bibitem{goedecker1995low}
{\sc S.~Goedecker}, {\em Low complexity algorithms for electronic structure
  calculations}, Journal of Computational Physics, 118 (1995), pp.~261--268.

\bibitem{hou2017adaptive}
{\sc Y.~T. Hou, D.~Huang, K.~C. Lam, and P.~Zhang}, {\em An adaptive fast
  solver for a general class of positive definite matrices via energy
  decomposition}, Preprint: arXiv:1707.08277v2 [math.NA].,  (2017).

\bibitem{hou2016sparse}
{\sc Y.~T. Hou and P.~Zhang}, {\em Sparse operator compression of higher-order
  elliptic operators with rough coefficients}, Research in Mathematical
  Sciences, in press,  (2017).

\bibitem{lehoucq1996deflation}
{\sc R.~B. Lehoucq and D.~C. Sorensen}, {\em Deflation techniques for an
  implicitly restarted arnoldi iteration}, SIAM Journal on Matrix Analysis and
  Applications, 17 (1996), pp.~789--821.

\bibitem{lehoucq1998arpack}
{\sc R.~B. Lehoucq, D.~C. Sorensen, and C.~Yang}, {\em ARPACK users' guide:
  solution of large-scale eigenvalue problems with implicitly restarted Arnoldi
  methods}, SIAM, 1998.

\bibitem{maalqvist2014localization}
{\sc A.~M{\aa}lqvist and D.~Peterseim}, {\em Localization of elliptic
  multiscale problems}, Mathematics of Computation, 83 (2014), pp.~2583--2603.

\bibitem{martinez2016tuned}
{\sc {\'A}.~Mart{\'\i}nez}, {\em Tuned preconditioners for the eigensolution of
  large spd matrices arising in engineering problems}, Numerical Linear Algebra
  with Applications, 23 (2016), pp.~427--443.

\bibitem{meirovitch1975elements}
{\sc L.~Meirovitch}, {\em Elements of vibration analysis}, McGraw-Hill, 1975.

\bibitem{newman2010networks}
{\sc M.~Newman}, {\em Networks: an introduction}, Oxford university press,
  2010.

\bibitem{ng2001spectral}
{\sc A.~Y. Ng, M.~I. Jordan, Y.~Weiss, et~al.}, {\em On spectral clustering:
  {A}nalysis and an algorithm}, in NIPS, vol.~14, 2001, pp.~849--856.

\bibitem{owhadi2017multigrid}
{\sc H.~Owhadi}, {\em Multigrid with rough coefficients and multiresolution
  operator decomposition from hierarchical information games}, SIAM Review, 59
  (2017), pp.~99--149.

\bibitem{owhadi2017universal}
{\sc H.~Owhadi and C.~Scovel}, {\em Universal scalable robust solvers from
  computational information games and fast eigenspace adapted multiresolution
  analysis}, arXiv preprint arXiv:1703.10761,  (2017).

\bibitem{ozolicnvs2013compressed}
{\sc V.~Ozoli{\c{n}}{\v{s}}, R.~Lai, R.~Caflisch, and S.~Osher}, {\em
  Compressed modes for variational problems in mathematics and physics},
  Proceedings of the National Academy of Sciences, 110 (2013),
  pp.~18368--18373.

\bibitem{romero2010parallel}
{\sc E.~Romero, M.~B. Cruz, J.~E. Roman, and P.~B. Vasconcelos}, {\em A
  parallel implementation of the jacobi-davidson eigensolver for unsymmetric
  matrices.}, in VECPAR, Springer, 2010, pp.~380--393.

\bibitem{schafer2017compression}
{\sc F.~Sch{\"a}fer, T.~Sullivan, and H.~Owhadi}, {\em Compression, inversion,
  and approximate pca of dense kernel matrices at near-linear computational
  complexity}, arXiv preprint arXiv:1706.02205,  (2017).

\bibitem{shi2000normalized}
{\sc J.~Shi and J.~Malik}, {\em Normalized cuts and image segmentation}, IEEE
  Transactions on pattern analysis and machine intelligence, 22 (2000),
  pp.~888--905.

\bibitem{sleijpen2000jacobi}
{\sc G.~L. Sleijpen and H.~A. Van~der Vorst}, {\em A jacobi--davidson iteration
  method for linear eigenvalue problems}, SIAM review, 42 (2000), pp.~267--293.

\bibitem{sorensen1992implicit}
{\sc D.~C. Sorensen}, {\em Implicit application of polynomial filters in
  ak-step arnoldi method}, Siam journal on matrix analysis and applications, 13
  (1992), pp.~357--385.

\bibitem{sorensen1997implicitly}
{\sc D.~C. Sorensen}, {\em Implicitly restarted arnoldi/lanczos methods for
  large scale eigenvalue calculations}, in Parallel Numerical Algorithms,
  Springer, 1997, pp.~119--165.

\bibitem{stewart1969accelerating}
{\sc G.~Stewart}, {\em Accelerating the orthogonal iteration for the
  eigenvectors of a hermitian matrix}, Numerische Mathematik, 13 (1969),
  pp.~362--376.

\end{thebibliography}
\end{document}